\documentclass[11pt,reqno]{amsart}
\usepackage{amsfonts}
\usepackage{amsmath}
\usepackage{amssymb}
\usepackage{caption}
\usepackage[colorlinks=true, pdfstartview=FitV, linkcolor=blue,  citecolor=blue, urlcolor=blue]
{hyperref}   
\usepackage{textcase}

\usepackage{tikz,float}
\usepackage{tikz-cd}
\usetikzlibrary{decorations.markings}
\tikzset{-{latex[length=3mm, width=2mm]}-/.style={decoration={
  markings,
  mark=at position #1 with {\arrow{>}}},postaction={decorate}}}
  \tikzset{middlearrow/.style={
        decoration={markings,
            mark= at position 0.55 with {\arrow{#1}} ,
        },
        postaction={decorate}
    }
}
\usetikzlibrary{patterns}

\usepackage{color}

\definecolor{DarkGreen}{rgb}{0.00, 0.39, 0.00}
\newcommand{\nn}{\nonumber}
\newcommand{\p}{\partial}
\newcommand{\supp}{\text{supp}}
\newcommand{\rr}{\mathbb{R}}

\newcommand{\zz}{\mathbb{Z}}

\newcommand{\NN}{\mathbb{N}}

\newcommand{\nin}{\noindent}

\def\dott
{\tikz[]
{\filldraw [black] (-0.1,-0.1) circle (0.8pt);}
}

\theoremstyle{plain}  % default
\newtheorem{theorem}{Theorem}[section]
\newtheorem{proposition}{Proposition}[section]
\newtheorem{lemma}{Lemma}[section]

\newtheorem{remark}{Remark}[section]
\numberwithin{figure}{section}
\numberwithin{equation}{section}

\begin{document}
\title{
The Nonlinear Schr\"odinger equation on the half-space
}
\author{A. Alexandrou Himonas* \&  Fangchi Yan}
\date{June 19, 2024
 \mbox{}$^*$\!\textit{Corresponding author}:
himonas.1@nd.edu}

\keywords{
Cubic Nonlinear Schr\"odinger equations, 
initial-boundary value problem,
Fokas Unified Transform Method, 
well-posedness in Sobolev spaces,  
linear space-time estimates,
 trilinear estimates in Bourgain spaces}

\subjclass[2020]{Primary: 35Q55, 35G31, 35G16, 37K10}

\begin{abstract} 
This work studies the initial-boundary value problem for both the 
linear Schr\"odinger equation
and the
cubic nonlinear Schr\"odinger equation on the half-space
in  higher dimensions ($n\ge 2$).
 First, the forced  linear problem is solved on the half-space via the Fokas method and then using the obtained
solution formula new and interesting linear estimates are
derived with data and forcing in appropriate spaces.
Second, the well-posedness of the  nonlinear problem 
on the half-space is proved with initial data in 
Sobolev spaces $H^s(\rr^n_+)$, with $s>\frac{n}{2}-1$,
and boundary data in natural
 Bourgain spaces $\mathcal{B}^s$ that reflect the 
 boundary regularity of the linear problem.
The proof method consists of showing that the iteration 
map defined via the Fokas solution formula is a contraction
by establishing sharper trilinear estimates.
The presence of the boundary introduces solution spaces that 
involve   temporal Bourgain spaces.

\end{abstract}

\maketitle

\markboth{
The Nonlinear Schr\"odinger equation on the half-space
}{
 A. Himonas \& F. Yan }
 
 \setcounter{tocdepth}{1}

%\tableofcontents

%
%
%
%%%%%%%%%%%%%%%%%%%%%
%
%	Introduction
%
%%%%%%%%%%%%%%%%%%%%%
%
%
%
\section{Introduction  and Results}
\label{sec:Intro}
\setcounter{equation}{0}
In this paper we study the 
cubic nonlinear Schr\"odinger (NLS) equation  on the half-space
\begin{subequations}
\label{NLS-ibvp}
\begin{alignat}{2}
\label{NLS-eqn}
&iu_t + \Delta u \pm |u|^2u=0,
&& 
\quad
(x, t)\in  \rr^{n-1} \times\rr^+ \times (0, T), 
\, \, T< 1/2,
\\
\label{NLS-ic}
&u(x, 0) = u_0(x),
&&
\quad
x\in\rr^{n-1}\times\rr^+,
\\
\label{NLS-bc}
&u(x', 0 ,t)
=
g_0(x', t),
&&
\quad
(x',t)\in\rr^{n-1}
\times(0,T),
\end{alignat}
\end{subequations}
and prove its local well-posedness in  natural data spaces
for all dimensions $n\ge 2$.
To state our results precisely we  introduce the 
data and  solution spaces,  starting with 
the following notation for a vector
$
x=(x_1,\cdots,x_n)=(x',x_n)$ in $\rr^n
$
 and its dual  $\xi=(\xi_1,\cdots,\xi_n)=(\xi',\xi_n)$,
$$
\xi x
=
\xi\cdot x
=
\langle \xi, x\rangle,
\quad
\xi^2
\doteq
\xi_1^2+\cdots+\xi_n^2
=
\xi'^2+\xi_n^2
\quad
\text{and}
\quad
|\xi|
\doteq
\sqrt{\xi_1^2+\cdots+\xi_n^2}
=
\sqrt{\xi'^2+\xi_n^2}.
$$
The initial data $u_0(x)$ belong in 
Sobolev spaces  $H^s(\rr^{n-1}\times\rr^+)$, $s>s_n$,
with  critical exponent
\begin{equation}
\label{critical-points-def}
s_n=\frac{n}{2}-1,
\quad
n\ge 2.
\end{equation}
These spaces are restrictions of the Sobolev spaces  $H^s(\rr^n)$ consisting of all  temperate distributions $F$
with finite norm
$
\big\| F \big\|_{H^s(\mathbb R^n)}
\doteq
\big(
\int_{\mathbb R^n} \hskip-0.02in (1+ |\xi|)^{2s}
|\widehat F(\xi)|^2 d\xi \big)^{1/2},
$
where
$\widehat F(\xi)$ is the  Fourier transform defined by 
$
\widehat F(\xi)
\doteq
\int_{\rr^n} e^{-i\xi\cdot x}$ $F(x)dx.
$
We recall that for  an open set $\Omega\subset \rr^n$
the space $H^s(\Omega)$  is defined by 
\begin{equation}
H^s(\Omega) \doteq \Big\{f: f= F\big|_{\Omega}\
 \mathrm{where}\ F\in H^s(\rr^n)\
  \mathrm{and}\
\| f \|_{H^s(\Omega)}\doteq \inf_{F \in H^s(\mathbb R^n)} 
\| F \|_{H^s(\rr^n)} <\infty
   \Big\}.
\end{equation}
The boundary data  $g_0(x', t)$ belong in 
the  {\it boundary Bourgain spaces $\mathcal{B}^s_T$},
which are the restriction of the spaces 
$\mathcal{B}^s$ on the interval $(0,T)$.  
For any $s\in\rr$ the  space
$\mathcal{B}^s=\mathcal{B}^s(\rr^{n-1}_{x'}\times \rr_t)$ is defined by the norm
\begin{equation}
\label{sharp-Bs}
\|g\|_{\mathcal{B}^s}^2
\doteq
\int_{\rr^n}
(1+\xi'^2+|\tau+\xi'^2|)^{s}
(1+|\tau+\xi'^2|)^{1/2}
|\widehat{g}(\xi',\tau)|^2
d\xi'
d\tau.
\end{equation}
Note that in the case $n=1$ these are  the familiar spaces
$H_t^{\frac{2s+1}{4}}$, which express the time regularity 
 for the linear Schr\"odinger equation.
Furthermore,   $\mathcal{B}^s_T=\mathcal{B}^s_{\rr^{n-1}\times(0,T)}$ is the restriction space 
\begin{equation}
\label{B-restrict}
\mathcal{B}^s_T
\doteq
\{
g:
g(x',t)
=
\tilde{g}(x',t)
\;\;\mbox{on}\;\; \rr^{n-1}\times(0,T)
\,\,
\text{with}
\,\,
\tilde{g}\in \mathcal{B}^s
\},
\end{equation}
with norm defined by
\begin{equation}
\label{sharp-Bs-restrict}
\|g\|_{\mathcal{B}^s_T}
\doteq
\inf\limits_{\tilde{g}\in \mathcal{B}^s}\{
\| \tilde{g}\|_{\mathcal{B}^s}:
\;\tilde{g}(x',t)
=
g(x',t)
\;\;\mbox{on}\;\; \rr^{n-1}\times(0,T)
\}.
\end{equation}
Finally, the solution/forcing spaces 
$X^{s,b}_{\rr^{n-1}\times\rr^+\times(0, T)}$
are restrictions of the well-known  Bourgain spaces $X^{s,b}(\rr^{n+1})$ associated to the linear  Schr\"odinger
equation, which for any real numbers $s$, $b$,  
are defined by the norm (see \cite{b1993-nls, kpv1996-nls})
\begin{equation}
\label{bourgain-nls}
\|u\|_{X^{s,b}}^2
=
\|u\|_{X^{s,b}(\rr^{n+1})}^2
\doteq
\int_{\rr^{n+1}}
(1+|\xi|)^{2s}
(1+|\tau+\xi^2|)^{2b}
|\widehat{u}(\xi,\tau)|^2
d\xi
d\tau.
\end{equation}
More precisely, the spaces 
$X^{s,b}_{\rr^{n-1}\times\rr^+\times(0, T)}$,
are defined as follows
\begin{equation}
\label{sb-restrict}
X^{s,b}_{\rr^{n-1}\times\rr^+\times(0, T)}
\doteq
\{
u:
u(x',x_n,t)
=
\tilde{u}(x',x_n,t)
\;\;\mbox{on}\;\; \rr^{n-1}\times\mathbb{R^+}\times(0,T)
\,\,
\text{with}
\,\,
\tilde{u}\in X^{s,b}(\rr^{n+1})
\},
\end{equation}
and are equipped  with the norm
\begin{equation}
\label{sb-restrict-norm}
\| u\|_{X^{s,b}_{\rr^{n-1}\times\rr^+\times(0, T)}}
\doteq
\inf\limits_{\tilde{u}\in X^{s,b}}\big\{
\| \tilde{u}\|_{X^{s,b}}:
\;\tilde{u}(x',x_n,t)
=
u(x',x_n,t)
\;\;\mbox{on}\;\; \rr^{n-1}\times\rr^+\times(0,T)
\big\}.
\end{equation}
In addition, for $s>1/2$ we utilize
 the following temporal Bourgain spaces
\begin{align}
\label{temporal-B}
\|u\|_{Y^{s,b}}^2
=
\|u\|_{Y^{s,b}(\rr^{n+1})}^2
\doteq
\int_{\rr^{n+1}}
(1+|\tau|)^{s}
(1+|\tau+\xi^2|)^{2b}
|\widehat{u}(\xi,\tau)|^2
d\xi
d\tau,
\,\, \, 
s, b \in \rr,
\end{align}
which in the one-dimensional space situation
were introduced in
%by Faminskii 
\cite{Fa2004, Fa2007}.
Their restriction is 
\begin{equation}
\label{Y-sb-restrict}
Y^{s,b}_{\rr^{n-1}\times\rr^+\times(0, T)}
\doteq
\{
u:
u(x',x_n,t)
=
\tilde{u}(x',x_n,t)
\;\;\mbox{on}\;\; \rr^{n-1}\times\mathbb{R^+}\times(0,T)
\,\,
\text{with}
\,\,
\tilde{u}\in Y^{s,b}(\rr^{n+1})
\},
\end{equation}
and is equipped with the norm
\begin{equation}
\label{Y-sb-restrict-norm}
\| u\|_{Y^{s,b}_{\rr^{n-1}\times\rr^+\times(0, T)}}
\doteq
\inf\limits_{\tilde{u}\in Y^{s,b}}\big\{
\| \tilde{u}\|_{Y^{s,b}}:
\;\tilde{u}(x',x_n,t)
=
u(x',x_n,t)
\;\;\mbox{on}\;\; \rr^{n-1}\times\rr^+\times(0,T)
\big\}.
\end{equation}

Using the above definitions, we state our main result as follows.
\begin{theorem}[Local well-posedness]
%[\textcolor{blue}{(small data) 
%Well-posedness for NLS ibvp in 
%high dimensions}]
\label{wp-thm-small}
If $n\ge 2$ and  $s> s_n\doteq\frac{n}{2}-1$ 
with $\frac{2s-1}{4}\not\in \NN_0$,
%(or $s\neq \frac12+2k$, $k=0,1,2,\cdots$),
%
then for sufficiently small initial data $u_0\in H^s(\rr_n^+)$
and  boundary data $g_0\in \mathcal{B}^s_T$
 %satisfying a  smallness condition
 the NLS ibvp \eqref{NLS-ibvp}, subject to the compatibility conditions \eqref{compa-condi}, 
 %(placed in Extensions lemma Section), 
 has a unique solution $u\in X^{s,b}_{\rr^{n-1}\times\rr^+\times(0,T)}\cap Y^{s,b}_{\rr^{n-1}\times\rr^+\times(0,T)}$,
 for some $b\in(0,\frac12)$,
 with lifespan 
 $T_0=T<1/2$.  
Moreover, the data-to-solution map  $\{u_0,g_0\}\mapsto u$ is locally Lipschitz continuous.
\end{theorem}

In dimension $n=2$ we have any size data. 
More precisely, we have the following result.
\begin{theorem}
%[\textcolor{blue}{Any data}]
\label{wp-thm}
If $n= 2$ and $0< s<\frac12$, 
then for initial data $u_0 \in H^s(\rr\times\rr^+)$ 
and  boundary data  $g_0\in \mathcal{B}_T^s$
there  is a lifespan  $0 < T_0\leq T<\frac12$ such that the 
NLS ibvp \eqref{NLS-ibvp}
has a unique solution $u\in X^{s,b}_{\rr\times\rr^+\times(0,T_0)}$ satisfying the size estimate
\begin{align}
\label{X-norm-def}
\|u\|_{X^{s,b}_{\rr\times\rr^+\times(0,T_0)}}
\lesssim
\|
u_0
\|_{H^{s}(\rr\times\rr^+)}
+
\|g_0\|_{\mathcal{B}_T^s},
\end{align}
for some $b\in(0,\frac12)$. Also, 
we have the following estimate for the lifespan 
\begin{align}
\label{mnls-lifespan}
T_0
=c_0
\Big[
1
+
\|
u_0
\|_{H^{s}(\rr\times\rr^+)}
+
\|g_0\|_{\mathcal{B}_T^s}
\Big]^{-32/s}.
\end{align}
Moreover, the data-to-solution map  $\{u_0,g_0\}\mapsto u$ is locally Lipschitz continuous.
\end{theorem}
We mention that the local well-posedness of the NLS on the half-plane 
%for $s>0$ and for any size data 
has been studied in 
%Ran, Sun and Zhang 
\cite{RSZ2018} via the Laplace transform 
and in 
%Himonas and Mantzavinos
\cite{hm2020}  via the Fokas method.
Concerning the dimension $n=1$ the well-posedness of 
the NLS initial-boundary value problem \eqref{NLS-ibvp} 
for $s\ge 0$
%and for arbitrary size data
has been studied earlier in several works.
%
% Holmer
In \cite{h2005} it was studied by adapting the method 
developed in 
%Colliander-Kenig
\cite{ck2002} for the gKdV, which is based on reducing the 
initial-boundary value problem to simpler initial value problems
and employing the Bourgain spaces theory of dispersive equations
(see also \cite{h2006}).
%
%J.L. Bona, S.M. Sun and B.-Y. Zhang
In \cite{bsz2018} it is studied via the  Laplace transform,
which was also used in earlier works for the KdV
(see \cite{bsz2002}, \cite{bsz2003}).
%
%M.B. Erdo\v gan and N. Tzirakis,
Furthermore,  following \cite{bsz2018},  in \cite{et2016} 
regularity properties and well-posedness of the cubic nonlinear Schrödinger equation on the half-line with rough initial data
are studied. Also, it is proved that the nonlinear part of the cubic NLS on the half-line is smoother than the initial data, like in the case 
of initial value problem.
In \cite{fhm2017} the Fokas method for the study
of the NLS ibvp \eqref{NLS-ibvp} on the half-line 
was implemented for the first time and local well-posedness 
was proved for $s>1/2$ in Hadamard solution spaces.
Finally, in \cite{hy2024-JMA-cNLS} the NLS ibvp  \eqref{NLS-ibvp}
on the half-line in low regularity spaces
 is studied as part the  initial-boundary value problem for the Schr\"odinger equation with cubic nonlinearities of the form $u^{3-k}\bar{u}^{k}$, $k=0,1, 2, 3$, using the Fokas method and  trilinear estimates in Bourgain spaces. It is interesting that for
the nonlinearities  $u^3$,  
$|u|^2u$, $|u|^2\bar{u}$  well-posedness holds when 
$s\ge 0$ while for  $\bar{u}^3$
 it  holds when $s>-1/3$.

 Our method for proving Theorem \ref{wp-thm-small}  
is analogous to that  used for the well-posedness
in Sobolev spaces $H^s$ of the corresponding initial value problem (ivp), where the iteration map is obtained via the Fourier transform and a fixed point is 
established by deriving linear and trilinear estimates in  spatial
Bourgain spaces $X^{s,b}$ with $b>1/2$.
For our NLS ibvp \eqref{NLS-ibvp}
 this method consists of three steps. 
First, using the Fokas unified transform method we solve the forced linear Schr\"odinger  ibvp.
Second, using the obtained solution formula we derive 
linear estimates and uncover the boundary data space
$\mathcal{B}^s$ and the natural solution/forcing spaces.
Third, we prove the trilinear estimates suggested by 
the nonlinearity and which are needed for making 
the iteration map defined by the Fokas solution formula 
a contraction.
As we will see below the novelties in ibvp 
arise from the boundary, which introduces 
new  temporal Bourgain spaces $Y^{s,b}$, 
in addition to the spaces $X^{s,b}$, both with $b<1/2$.
This requires sharper trilinear estimates for 
establishing well-posedness. 
%

%
%%%%%%%%%%%%%%
%
%
%     The Linear Problem
%
%
%%%%%%%%%%%%%%%
%
\vskip0.05in
\nin
{\bf Solving the Linear Problem.}
Replacing in our NLS ibvp \eqref{NLS-ibvp} the nonlinearity 
by a forcing $f(x, t)$ we obtain the following  forced linear Schr\"odinger ibvp
\begin{subequations}
\label{lnls-ibvp}
\begin{align}
\label{lnls-sym:eqn}
&i\partial_t u
+
\Delta_{x'} u
+
\partial^2_{x_n} u
=
f,
\hskip-0.7in
&& 
(x', x_n, t)\in \rr^{n-1} \times \rr^+ \times (0, T),
\\
\label{lnls-sym:ic}
&u(x',x_n,0) 
= 
u_0 (x',x_n),
\hskip-0.7in
&& (x', x_n)\in \rr^{n-1} \times \rr^+,
\\
\label{lnls-sym:bc}
&u(x',0, t) 
= 
g_0(x',t),
\hskip-0.7in
&& (x', t)\in \rr^{n-1} \times (0, T)
.
\end{align}
\end{subequations}
Then, applying the Fokas method, also referred as unified transform method (an outline of which we present 
in Section  \ref{sec:Fokas-solution}), we find that the
solution to the linear  problem \eqref{lnls-ibvp} is given by the 
formula:
\begin{align}
\label{nd-lnls-utm-sln}
u(x,t)
=
S[u_0,g_0;f](x,t)
=&
\frac{1}{(2\pi)^{n}}
\int_{\rr^{n}}
e^{i\xi\cdot x}
\cdot
e^{-i\xi^2 t}
\cdot
[\widehat{u}_{0}(\xi',\xi_n)-i\widetilde{F}(\xi',\xi_n,t)]
d\xi
\\
-&
\frac{1}{(2\pi)^{n}}
\int_{\rr^{n-1}}
\int_{\p D^+}
e^{i\xi'\cdot x'+i\xi_n x_n}
\cdot
e^{-i\xi'^2t-i\xi_n^{2}t}
[\widehat{u}_{0}(\xi',-\xi_n)-i\widetilde{F}(\xi',-\xi_n,t)]
d\xi_n
d\xi'
\nonumber
\\
+&
\frac{2}{(2\pi)^{n}}
\int_{\rr^{n-1}}
\int_{\p D^+}
e^{i\xi'\cdot x'+i\xi_n x_n}
\cdot
e^{-i\xi'^2t-i\xi_n^{2}t}
\cdot
\xi_n
\cdot
\tilde{g}_{0}(\xi',-\xi'^2-\xi_n^{2},T)
d\xi_n
d\xi',
\nn
\end{align}
where
\begin{equation}
\label{initial-FT-half}
\widehat {u}_{0}(\xi',\xi_n)
\doteq
\int_0^\infty
\int_{\rr^{n-1}}
e^{-i\xi'\cdot x'-ix_n\xi_n}
u_{0}(x',x_n)
dx'
dx_n,
\quad 
\xi'\in\rr^{n-1},
\,\,
\text{Im}\, \xi_n \le 0,
\end{equation}
\begin{align}
\label{forcing-time-trans}
&\widetilde{F}(\xi',\xi_n, t)
\\
\doteq&
\int_0^t \int_0^\infty  \int_{\rr^{n-1}}
e^{i\xi'^2\tau+i\xi_n^2\tau}
e^{-i\xi'\cdot x'-i\xi_n x_n}
f(x',x_n,\tau)
dx'
dx_nd\tau,
\quad
\xi'\in\rr^{n-1},
\,\,
\text{Im} \xi_n \le 0,
\,\,
0<t<T,
\nonumber
\end{align}
and
\begin{align}
\label{bc-time-trans}
\tilde{g}_{0}(\xi',\tau,T)
\doteq
\int_0^T
\int_{\rr^{n-1}}
e^{-i\tau t} 
e^{-i\xi'\cdot x'}
\cdot
g_0(x',t)
dx'
dt,
\quad 
\,\,
\xi'\in\rr,
\,\,
\tau\in \mathbb{C}. 
\end{align}
\begin{minipage}{1\linewidth}
\begin{center}
\begin{tikzpicture}[scale=0.8]
%%%%%%%%%%%%%%%%%%%
%
%Variables defined
%
%%%%%%%%%%%%%%%%%%%%
\newcommand\X{0};
\newcommand\Y{0};
\newcommand\FX{10};
\newcommand\FY{10};
\newcommand\R{0.6};
%%%%%%%%%%%%%%%%%%
%
%End
%
%%%%%%%%%%%%%%%%%%%%%
%
\filldraw [fill=gray, fill opacity=0.2, draw opacity=0.01,variable=\x,domain=0:pi/2] (0,0)--
plot({2.7*cos(deg(\x))},{2.7*sin(deg(\x))})
--cycle;

\draw [] (-3,0)--(-1.5,0);
\draw [] (0,0)--(-1.5,0);
\draw [blue] (3,0)--(1.5,0);
\draw [blue,-{Stealth[scale=1.3]}] (0,0)--(1.5,0);
\draw [blue,-{Stealth[scale=1.3]}] (0,3)--(0,1.5);
\draw [blue] (0,1.5)--(0,0);
\draw [] (1.5,1.5) node {\fontsize{\FX}{\FY}$D^+$};

\draw[line width=0.5pt, black] 
(2.1,2.9)--(2.1,2.4) (2.6,2.4)--(2.1,2.4);
\node[] at ({2.4},{2.7}) {\fontsize{\FX}{\FY}$\xi_n$}
;

\end{tikzpicture}
\end{center}
\vskip-0.1in
\captionof{figure}{Domain $D^+$}
\label{nls-domain}
\end{minipage}

%
%%%%%%%%%%%%%%
%
%
%     The Linear Problem
%
%
%%%%%%%%%%%%%%%
%
\vskip0.1in
\nin
{\bf Linear Estimates.}  Next, we estimate 
 the Fokas solution $S\big[u_0,g_0;f\big]$ defined
 by  formula \eqref{nd-lnls-utm-sln}  in terms 
 of the initial, boundary data 
 $u_0$, $g_0$, and the forcing $f$.
 Choosing the initial data $u_0$ in Sobolev spaces $H^s$
 leads in having the boundary data belonging in 
  natural
  boundary Bourgain spaces $\mathcal{B}^s$, when 
  the solution is estimated in spatial or temporal 
  Bourgain spaces. The spaces for the forcing depend
  on the range of the Sobolev space $s$.
  More precisely we have the following result.
\begin{theorem}
[Linear estimates]
\label{forced-linear-nls-thm}   
Let  $0<T<1/2$ and $0\le b<1/2$.
Then the solution $S\big[u_0,g_0;f\big]$
of the forced linear NLS ibvp  \eqref{lnls-ibvp},
defined by the Fokas formula \eqref{nd-lnls-utm-sln}, 
satisfies the following  spatial and temporal estimates,
under the compatibility condition \eqref{compa-condi}.

\vskip0.03in
\nin
$\bullet$
 If  $-\frac 12 \le s<\frac 12$, then
$S\big[u_0,g_0;f\big]$ satisfies the following estimate 
in spatial Bourgain spaces
\begin{align}
\label{noY-forced-linear-nls-est}
\|S\big[u_0,g_0;f\big]\|_{X_{\rr^{n-1}\times\rr^+\times(0,T)}^{s,b}}
\hskip-0.03in
\le
\hskip-0.03in
c_1\Big[
\|u_0\|_{H_x^s(\rr^{n-1}\times\rr^+)}
+
\|g_0\|_{\mathcal{B}_T^s}
+
 \|f\|_{X^{s,-b}_{\rr^{n-1}\times\rr^+\times(0,T)}}
 \Big].
\end{align}
\vskip0.01in
\nin
$\bullet$
If 
$s>-\frac32$ and 
 $\frac{2s-1}{4}\not\in \NN_0$,
then  the spatial $X^{s,b}$ norm of $S\big[u_0,g_0;f\big]$  is bounded as follows
\begin{align}
\label{X-forced-linear-nls-est}
 \hskip-0.15in
\|S\big[u_0,g_0;f\big]\|_{X_{\rr^{n-1}\times\rr^+\times(0,T)}^{s,b}}
 \hskip-0.25in
\le
c_1\Big[
\|u_0\|_{H_x^s(\rr^{n-1}\times\rr^+)}
\hskip-0.02in
+
\hskip-0.02in
\|g_0\|_{\mathcal{B}_T^s}
\hskip-0.02in
+
\hskip-0.02in
 \|f\|_{X^{s,-b}_{\rr^{n-1}\times\rr^+\times(0,T)}}
 \hskip-0.15in
+
\hskip-0.01in
 \|f\|_{Y^{s,-b}_{\rr^{n-1}\times\rr^+\times(0,T)}}
\Big].
\end{align}
\vskip0.01in
\nin
$\bullet$
If $s\ge 0$
and
%$s\neq\frac12+2k$,  
$\frac{2s-1}{4}\not\in \NN_0$, 
then  the temporal Bourgain $Y^{s,b}$ norm of $S\big[u_0,g_0;f\big]$  is bounded  as follows
\begin{align}
\label{Y-forced-linear-nls-est}
 \hskip-0.07in
\|S\big[u_0,g_0;f\big]\|_{Y_{\rr^{n-1}\times\rr^+\times(0,T)}^{s,b}}
 \hskip-0.2in
\le
c_1\Big[
\|u_0\|_{H_x^s(\rr^{n-1}\times\rr^+)}
\hskip-0.02in
+
\hskip-0.02in
\|g_0\|_{\mathcal{B}_T^s}
\hskip-0.02in
+
\hskip-0.02in
 \|f\|_{X^{s,-b}_{\rr^{n-1}\times\rr^+\times(0,T)}}
\hskip-0.08in
+
\hskip-0.02in
 \|f\|_{Y^{s,-b}_{\rr^{n-1}\times\rr^+\times(0,T)}}
\Big].
\end{align}
\end{theorem}
Observe that in linear estimate \eqref{noY-forced-linear-nls-est},
which holds for $-\frac 12 \le s<\frac 12$,
the solution $S\big[u_0,g_0;f\big]$ and the forcing 
$f$ are both in spatial Bourgain spaces.
Considering that the trilinear estimates stated below hold 
for $s>\frac{n}{2}-1$, we see that we can utilize spatial estimate 
 \eqref{noY-forced-linear-nls-est} 
 only for $0 \le s<\frac 12$ (low regularity data)
 and  prove well-posedness 
  for arbitrary size data in this 
 range of the Sobolev exponent $s$ (see Theorem \ref{wp-thm}).
 For data of higher smoothness ($s>1/2$) estimating 
 the spatial norm $X^{s,b}$ of the solution requires 
 the temporal norm $Y^{s,-b}$ of the forcing, in addition 
 to its  $X^{s,-b}$ norm.
 This leads to the companion estimate 
 \eqref{Y-forced-linear-nls-est} in the temporal $Y^{s,b}$ norm.
 Combining these norms enables us to prove 
 well-posedness for higher smoothness data
 and in all dimensions.
Related 
%Strichartz 
estimates for the Schr\"odinger equation on the half-space and 
in any dimension were derived in \cite{a2018}.
In the two-dimensional case and on the half-plane related linear estimates were derived in \cite{RSZ2018} and \cite{hm2020}.

%
%%%%%%%%%%%%%%
%
%
%    Trilinear Estimates
%
%
%%%%%%%%%%%%%%%
%
\vskip0.1in
\nin
{\bf Trilinear Estimates.} Next, we state the trilinear estimates suggested by the $X^{s,-b}$ and $Y^{s,-b}$ norms of the forcing $f$ appearing in linear estimates of Theorem \ref{forced-linear-nls-thm}, when $f$ is replaced by the NLS nonlinearity.
These  estimates help us show that 
the iteration map defined by the Fokas solution formula $S\big[u_0,g_0;\mp|u|^2u\big]$ is a contraction 
on a ball of the solution space 
$
X^{s,b}_{\rr^{n-1}\times\rr^+\times(0,T)}
\cap
Y^{s,b}_{\rr^{n-1}\times\rr^+\times(0,T)}
$
in Theorem \ref{wp-thm-small}, and the  space 
$
X^{s,b}_{\rr\times\rr^+\times(0,T)}
$
in Theorem \ref{wp-thm}.
More precisely, we have the following result.

\begin{theorem} [Trilinear estimates in Bourgain type spaces]
\label{Trilinear-estimate-thm}  
If $s>s_n\doteq \frac{n}{2}-1$, $n\ge 2$,  then
there exists $b$ and $b'$ with $\frac12-\beta_{n}\le b'\le b<\frac12$ 
such that the following trilinear estimates hold
\begin{align}
\label{trilinear-est-B}
&\| f\bar{g}h \|_{X^{s,-b}}
\le
c_2
\| f \|_{X^{s,b'}} \| g \|_{X^{s,b'}}\| h \|_{X^{s,b'}},
\\
\label{trilinear-est-Y}
&\|f\bar{g}h\|_{Y^{s,-b}}
\le
c_2
(\|f\|_{X^{s,b'}}+\|f\|_{Y^{s,b'}})
(\|g\|_{X^{s,b'}}+\|g\|_{Y^{s,b'}})
(\|h\|_{X^{s,b'}}+\|h\|_{Y^{s,b'}}),
\end{align}
where $\beta_n>0$ is defined as follows
\begin{align}
\label{beta-choice}
\beta_n
=
\min\{\frac{1}{16},\frac18+\frac18s-\frac1{16}n\},
\quad
n\ge 2.
\end{align}
\end{theorem}
Observe that the trilinear estimates above are  sharper
than the ones used for the well-posedness of NLS ivp,
where it is assumed that $b>1/2$, while here
the linear estimates require that $b<1/2$.
Also, the result below states that these estimates are
 optimal.

\begin{theorem}
\label{optimal-thm}
The trilinear estimates \eqref{trilinear-est-B} in Theorem \ref{Trilinear-estimate-thm} 
 are optimal, 
 i.e., they fail for
  $s\le \frac{n}{2}-1$.
\end{theorem}

Concerning the initial value problem (ivp) of the nonlinear Schr\"odinger equation on the whole space $\rr^n$, 
 it is well-posed for 
 initial data in Sobolev spaces $H^s(\rr^n)$ when $s>\frac{n}{2}-1$.
 This follows from a more general result in   \cite{cw-nls1990},
 where the higher order nonlinearity $|u|^{\alpha-1}u$, 
 $\alpha>1$, is considered and where  Sobolev and 
 Besov spaces are utilized for their study.
 Related earlier results are proved in 
  \cite{gv-1979,gv-1978,gv-1985,ka-nls1987,Ts-nls1987}.
 In the framework of Bourgain spaces $X^{s,b}(\rr^{n+1})$,
 one can prove this ivp result in a way which is 
  analogues to the proof of our Theorem \ref{wp-thm-small}.
For the forced linear ivp, the analogous to 
the forced linear ibvp estimate 
\eqref{noY-forced-linear-nls-est} is the following 
 much simpler estimate
\begin{align}
\label{forced-linear-nls-est}
\|S\big[u_0;f\big]\|_{X_{\rr^{n}\times(0,T)}^{s,b}}
\le
c_1\big[
\|u_0\|_{H_x^s(\rr^{n})}
+
 \|f\|_{X^{s,b-1}_{\rr^{n}\times(0,T)}}
 \big],
\end{align}
which holds for any $s$ and $\frac 12<b<1$.
The solution to the forced linear ivp $S\big[u_0;f\big]$
is obtained via the Fourier transform in a 
straightforward manner.
And, the well-posedness proof consists of
combining the
linear estimate \eqref{forced-linear-nls-est},
with the ivp trilinear estimate
\begin{align}
\label{trilinear-est-B-ivp}
&\| f\bar{g}h \|_{X^{s,b-1}}
\le
c_2
\| f \|_{X^{s,b'}} \| g \|_{X^{s,b'}}\| h \|_{X^{s,b'}},
\end{align}
and a basic multiplier lemma \cite{tao-book},
to deal with large data.
The presence of the boundary in the ibvp case requires 
the more sophisticated Fokas method in order 
to produce a formula for the solution of the forced 
linear ibvp and this becomes possible by using 
 integration over appropriate contours in the complex plane.
Furthermore, due to the boundary, numerous compatibility 
conditions arise (see \eqref{compa-condi}) that depend 
on the regularity of the data, and new temporal Bourgain spaces appear in the derivation of the linear estimates.
All this makes the analysis of the NLS ibvp much more delicate than the one for the NLS ivp
although, thanks to the Fokas solution formula, 
there is a certain analogy between them.
We note that there is an extensive literature on the well-posedness 
of the initial value problem for the NLS, KdV and related dispersive 
equations for which 
we refer the reader to the works
\cite{
Bourgain-book,
kpv1989,kpv1991-reg,kpv1993, kpv1996,
kv2019,
lpbook, 
ss1999, tao2001}
and the references therein.
Concerning results on initial-boundary value problem for NLS
and related equations using  the Fokas method we refer to the works
\cite{dtv2014, f1997,f2002, fhm2016, fis2005, fl2012, fpbook2015, fs2012,hm2015,hm2021, hy2022-KdVm,
hy2023-NA-mKdV, lf2012a,lenells2013}
and the references therein.

\vskip0.05in
\noindent
The nonlinear Schr\"odinger equation is one of the two 
 leading models in the theory of nonlinear dispersive equations (the other being the Korteweg–de Vries (KdV) equation). It  arises in many physical situations.
 For example, it appears  in nonlinear optics \cite{t1964},   plasmas \cite{ww1977}, Bose-Einstein condensates \cite{ps2003}, and water waves \cite{p1983,css1992}. 
Moreover, when $n=1$ it is an integrable equation
(like the KdV) and can be studied via the inverse scattering transform  \cite{zs1972}.

%
%%%%%%%%%%%%%%%%%%%%
%
%
%	 Structure
%
%
%%%%%%%%%%%%%%%%%%%%
%

\vskip0.05in
\noindent
{\bf \large  Structure.}  
The paper is organized as follows.
In section 2, we derive the linear estimates
for the reduced pure linear ibvp and
 uncover the natural spaces $\mathcal{B}^s$
for the boundary data.
Then, in section 3, using this estimate 
and decomposing the forced linear ibvp 
into simpler problems (two of which are 
initial value problems) we prove the linear estimates stated in Theorem \ref{forced-linear-nls-thm}.
In Section 4, we prove the estimates for homogeneous ivp
while in Section 5 we prove 
the estimates
for 
non-homogeneous ivp.
Then, in section 6 we prove the needed trilinear estimates in spatial Bourgain spaces.
In section 7, we prove the companion temporal trilinear estimates.
In section 8, we  demonstrate their optimality.
 In Section 9, we conclude by proving our main 
well-posedness result.
Finally,
in Section 10, we outline the derivation
of  the Fokas solution formula for the forced linear 
Schr\"odinger ibvp.

%
%
%%%%%%%%%%%%%%%%%%%%
%
%
%             Reduced pure ibvp
%
%
%%%%%%%%%%%%%%%%%%%% 
%
%
\section{Linear Estimates -- The Reduced Pure ibvp
 of the Schr\"odinger equation
 }
\label{sec:reduced-pure-ibvp}
\setcounter{equation}{0}
The basic step in our derivation of  the linear estimates for the  Schr\"odinger equation, stated in  Theorem \ref{forced-linear-nls-thm},  is to first derive such estimates for its reduced pure ibvp,
 which reads as follows
\begin{subequations}
\label{nd-nls-pure}
\begin{alignat}{2}
&iv_t + \Delta_{x'} v+\p_{x_n}^2v=0,
\quad
&&
(x', x_n)\in \rr^{n-1}_{x'} \times \rr^+_{x_n},\,\,
 t\in (0, 2), 
\\
&v(x', x_n, 0) = 0, 
\\
&v(x', 0, t) = h(x', t),
\quad
&& \supp (h) \subset \rr^{n-1}\times (0,2).
\end{alignat}
\end{subequations}
Using  Fokas  formula \eqref{nd-lnls-utm-sln} with $v_0=0$ and $f=0$, we get its solution 
\begin{align}
\label{pure-ibvp-utm}
v(x',x_n,t)
=
\frac{2}{(2\pi)^n}
\int_{\xi'\in{\rr^{n-1}}} 
\int_{\xi_n\in\p D^+} 
e^{i\xi'\cdot x'+i\xi_nx_n -i(\xi'^2+\xi_n^2)t}
\cdot
\xi_n 
\cdot
\tilde{h}(\xi',-\xi'^2-\xi_n^2,2)
d\xi_n
d\xi',
\end{align}
where $\tilde{h}(\xi',-\xi'^2-\xi_n^2,2)$ is given by \eqref{bc-time-trans}. 
Furthermore, since $\supp\, h\in \rr^{n-1}\times (0,2)$ we have
\begin{align*}
\tilde{h}(\xi',-\xi'^2-\xi_n^2,2)
=&
\int_0^2
\int_{\rr^{n-1}}
e^{-i\xi'\cdot x'}
e^{i(\xi'^2+\xi_n^2)\tau}
h(x',\tau)
\, 
dx'
d\tau
\\
=&
\int_\rr
e^{i(\xi'^2+\xi_n^2)t}
\int_{\rr^{n-1}}
e^{-i\xi'\cdot x'}
h(x',t)
\, 
dx'
dt
=
\widehat{h}
\big(
\xi',-\xi'^2-\xi_n^2
\big),
\nonumber
\end{align*}
where  $\widehat{h}=\widehat{h}^{x',t}$ is the full Fourier transform of $h$.

Next, we state the estimates for the solution of reduced pure ibvp \eqref{nd-nls-pure}.
\begin{theorem}
\label{pure-ibvp-thm}
For a test function 
$h(x',t)$ supported in  $\rr_{x'}^{n-1} \times (0, 2)$,
Fokas  formula \eqref{nd-lnls-utm-sln} defines a solution $v$ to the reduced pure ibvp \eqref{nd-nls-pure} that satisfies the following Bourgain spaces  estimates
\begin{align}
\label{B-est-redu-pure}
\|
S[0,h;0]
\|_{X^{s, b}_{\rr^{n-1}\times\rr^+\times (0,2)}} 
\le
c_{s,b}
\|h\|_{\mathcal{B}^s},
\quad s>-\frac32,
\quad
0\le b<\frac12,
\\
\label{Y-B-est-redu-pure}
\|
S[0,h;0]
\|_{Y^{s, b}_{\rr^{n-1}\times\rr^+\times (0,2)}} 
\le
c_{s,b}
\|h\|_{\mathcal{B}^s},
\quad s\ge 0,
\quad
0\le b<\frac12.  
\end{align}
\end{theorem}

\begin{remark}
\label{optimal-homogeneous-est}
As we will see in the proof below,
combining inequalities \eqref{optimal-hom-est}, \eqref{Bs-est-v0},  \eqref{Bs-est-v1}, \eqref{V0-Y-est} with \eqref{Y-norm-v1-no2}, 
 we obtain the following optimal homogeneous estimates
\begin{align}
\label{B-est-redu-pure-optimal}
\|
S[0,h;0]
\|_{X^{s, b}_{\rr^{n-1}\times\rr^+\times (0,2)}} 
\le
c_{s,b}
\|h\|_{\overset{\dott}{\mathcal{B}^s}},
\quad s>-\frac32,
\quad
0\le b<\frac12,
\\
\label{Y-B-est-redu-pure-optimal}
\|
S[0,h;0]
\|_{Y^{s, b}_{\rr^{n-1}\times\rr^+\times (0,2)}} 
\le
c_{s,b}
\|h\|_{\overset{\dott}{\mathcal{B}^s}},
\quad s\ge 0,
\quad
0\le b<\frac12.  
\end{align}
where  $\overset{\dott}{\mathcal{B}^s}=\overset{\dott}{\mathcal{B}^s}(\rr_{x'}^{n-1}\times\rr_t)$ is the following   homogeneous
version of $\mathcal{B}^s$  \cite{a2018}
\begin{equation}
\label{homo-sharp-Bs}
\|h\|_{\overset{\dott}{\mathcal{B}^s}}^2
\doteq
\int_{\rr^n}
(1+\xi'^2+|\tau+\xi'^2|)^{s}
|\tau+\xi'^2|^{1/2}
|\widehat{h}(\xi',\tau)|^2
d\xi'
d\tau.
\end{equation}
\end{remark}

In the proof of Theorem \ref{pure-ibvp-thm} and throughout this work, 
we will need a standard  {\bf time localizer} 
 $\psi(t)$, which  is defined as follows:
\begin{equation}
\label{time-localizer}
\psi \in C^{\infty}_0(-1, 1), 
\,\,\,
0\le \psi \le 1
\,\,
\text{ and }
\,\,
\psi(t)=1
\,\,
\text{ for }
\,\,
|t|\le 1/2.
\end{equation}
Also, for two quantities $Q_1$ and $Q_2$
depending on several variables, we  write
$Q_1 \lesssim Q_2$ if there is a constant 
$c>0$ such that $Q_1 \le c\, Q_2$.  If $Q_1 \lesssim Q_2$ and $Q_2\lesssim Q_1$,
then we write $Q_1\simeq Q_2$.
\vskip0.05in
\nin
{\bf Proof of  Theorem \ref{pure-ibvp-thm}.}
For the integration over $\p D^+$, using the parametrization $[0,\infty)\ni \xi_n\rightarrow i  \xi_n$ over the imaginary axis,
we rewrite  Fokas formula \eqref{pure-ibvp-utm} 
as follows  (see Figure \ref{domain-D-plus})
\begin{equation}
\label{sln-rpibvp-split}
v(x',x_n,t)
=
\frac{2}{(2\pi)^n}
v_r(x',x_n,t)
-
\frac{2}{(2\pi)^n}
v_i(x',x_n,t),
\end{equation}
where $v_r$ and $v_i$ are the integrals of $\xi_n$ over the real and imaginary axis respectively, that is
\noindent

\hskip-0.2in
\begin{minipage}{0.75\linewidth}
\begin{align}
\label{vr-def}
&v_r
\doteq
\int_{\rr^{n-1}} 
\int_{\xi_n=0}^\infty 
e^{i\xi'\cdot x'+i \xi_nx_n -i(\xi'^2+\xi_n^2)t}
\cdot
\xi_n
\widehat{h}^{x',t}
\big(
\xi',-\xi'^2-\xi_n^2
\big)
d\xi_n
d\xi',
\\
\label{vi-def}
&v_i
\doteq
\int_{\rr^{n-1}} 
\int_{\xi_n=0}^\infty 
e^{i\xi'\cdot x'-\xi_nx_n -i(\xi'^2-\xi_n^2)t}
\cdot
\xi_n 
\widehat{h}^{x',t}
\big(
\xi',-\xi'^2+\xi_n^2
\big)
d\xi_n
d\xi'.
\end{align}
\end{minipage}
\hskip-0.2in
\begin{minipage}{0.35\linewidth}
\begin{center}
\begin{tikzpicture}[scale=0.7]
%%%%%%%%%%%%%%%%%%%
%
%Variables defined
%
%%%%%%%%%%%%%%%%%%%%
\newcommand\X{0};
\newcommand\Y{0};
\newcommand\FX{9};
\newcommand\FY{9};
\newcommand\R{0.6};
%%%%%%%%%%%%%%%%%%
%
%End
%
%%%%%%%%%%%%%%%%%%%%%
%
\filldraw [fill=gray, fill opacity=0.2, draw opacity=0.01,variable=\x,domain=0:pi/2] (0,0)--
plot({2.7*cos(deg(\x))},{2.7*sin(deg(\x))})
--cycle;

\draw [] (-2,0)--(-1.5,0);
\draw [] (0,0)--(-1.5,0);
\draw [blue] (3,0)--(1.5,0);
\draw [blue,-{Stealth[scale=1.3]}] (0,0)--(1.5,0);
\draw [blue,-{Stealth[scale=1.3]}] 
(0,3)
%node[yshift=0.2cm]
%{\fontsize{\FX}{\FY}Re$\xi_n$}
--(0,1.5);
\draw [blue] (0,1.5)--(0,0);
\draw [] (1.5,1.5) node {\fontsize{\FX}{\FY}$D^+$};

\draw[line width=0.5pt, black] 
(2.1,2.9)--(2.1,2.4) (2.6,2.4)--(2.1,2.4);
\node[] at ({2.4},{2.7}) {\fontsize{\FX}{\FY}$\xi_n$}
;

\draw[line width=1pt,black]
(1,0)
node[yshift=-0.2cm]
{\fontsize{\FX}{\FY}$v_r$}
(0,1.5)
node[xshift=-0.2cm]
{\fontsize{\FX}{\FY}$v_i$}

;

\end{tikzpicture}
\end{center}
\vskip-0.2in
\captionof{figure}{
%$D^+$
}
\label{domain-D-plus}
\end{minipage}

\vskip0.05in
\noindent
Now, we  prove Bourgain space estimate \eqref{B-est-redu-pure} for $v_r$ and $v_i$.

\vskip.05in
\noindent
{\bf Estimate \eqref{B-est-redu-pure} for $v_r$.} We begin by extending  $v_r$ 
in $x_n$ from $\rr^+$ to $\rr$, and in 
$t$  from $(0, 2)$ to $\rr$,  as follows
(since they appear in an oscillatory way)
\begin{equation*}
%\label{vr-ext}
V_r
=
\int_{\rr^{n-1}}
\int_{\xi_n=0}^\infty
e^{i\xi'\cdot x'+i\xi_nx_n -i(\xi'^2+\xi_n^2)t}
\cdot
\xi_n 
\cdot
\widehat{h}^{x',t}
\big(
\xi',-\xi'^2-\xi_n^2
\big)
d\xi_n
d\xi',
\quad
(x',x_n,t)
\in
\rr^{n+1}.
\end{equation*} 
In order  to prove estimate \eqref{B-est-redu-pure} for $v_r$, it suffices to show that 
\begin{equation}
\label{Vr-B-est}
\|\psi_4(t)V_r\|_{X^{s,b}(\rr^{n+1})}
\lesssim
\|h\|_{\mathcal{B}^s},
\quad
\text{for all}
\,\,
s\in \rr
\,\,
\text{and}
\,\,
b\ge 0,
\end{equation}
where $\psi_4=\psi(t/4)$ with $\psi$ being the time localizer given by \eqref{time-localizer}. 

\vskip0.05in
\noindent
\underline{Proof of inequality \eqref{Vr-B-est}.}
Taking Fourier transform, we obtain
\begin{align*}
%\label{Vr-full-FT}
\widehat{\psi_4(t) V_r}(\xi',\xi_n,\tau)
\simeq
\chi_{\xi_n>0}
(\xi_n)
\cdot
\widehat{\psi}_4(\tau+\xi^2)
\cdot
\xi_n 
\cdot
\widehat{h}^{x',t}
\big(
\xi',-\xi'^2-\xi_n^2
\big).
\end{align*}
Hence, we have
\begin{align*}
\|
\psi_4
\cdot
V_r
\|_{X^{s,b}}^2
\lesssim
\int_{\rr^{n-1}}
\int_{\xi_n=0}^\infty
\Big[
\int_\rr (1+|\tau+\xi^2|)^{2b}|\widehat{\psi}_4(\tau+\xi^2)|^2
d\tau
\Big]
(1+|\xi|)^{2s}
\xi_n^2
\big|
\widehat{h}^{x',t}
(
\xi',-\xi'^2-\xi_n^2)
\big|^2
d\xi_n
d\xi'.
\end{align*}
Since
$
\int_\rr (1+|\tau+\xi^2|)^{2b}|\widehat{\psi}_4(\tau+\xi^2)|^2
d\tau
=
\|(1+|\cdot|)^{2b}|\widehat{\psi}_4(\cdot)|^2\|_{L^1}
=
\|\psi_4\|_{H^{b}}^2
\lesssim
1,
$
the last inequality gives
\begin{align*}
%\label{pure-ibvp-Vr-ine}
\|
\psi_4(t)
V_r
\|_{X^{s,b}(\rr^{n+1})}^2
\lesssim
\int_{\rr^{n-1}}
\int_{\xi_n=0}^\infty
(1+|\xi|)^{2s}
\xi_n^2
\big|
\widehat{h}^{x',t}
(
\xi',-\xi'^2-\xi_n^2
)
\big|^2
d\xi_n
d\xi'.
\end{align*}
Now, doing the change of variables  
$
\xi'=\xi'$ and $\tau=-\xi'^2-\xi_n^2
$
 with Jacobian  
$
J=\frac{\partial (\xi', \tau)}{\partial (\xi', \xi_n)}=-2\xi_n
$
we get
\begin{align*}
%\label{almost-optimal-hom-est}
\|
\psi_4
\cdot
V_r
\|_{X^{s,b}(\rr^{n+1})}^2
\lesssim&
\int_{\rr^{n-1}}
\int_{\tau=-\infty}^{-\xi'^2}
\big(
1+[\xi'^2+|\tau+\xi'^2|]^{1/2}
\big)^{2s}
|\tau+\xi'^2|^{1/2}
\big|
\widehat{h}^{x',t}
\big(
\xi',\tau
\big)
\big|^2
d\tau
d\xi'.
\end{align*}
Furthermore,  using the equivalence 
$
1+[\xi'^2+|\tau+\xi'^2|]^{1/2}
\simeq
[1+\xi'^2+|\tau+\xi'^2|]^{1/2}
$, which follows from the relation
$\sqrt{1+a^2}\le 1+|a| \le  \sqrt{2}\sqrt{1+a^2}$,
and extending the $\tau$-integration over the whole $\rr$ we get
\begin{align}
\label{optimal-hom-est}
\|
\psi_4
\cdot
V_r
\|_{X^{s,b}(\rr^{n+1})}^2
\lesssim&
\int_{\rr^{n-1}}
\int_{\tau\in\rr}
\big(
1+\xi'^2+|\tau+\xi'^2|
\big)^{s}
|\tau+\xi'^2|^{1/2}
\big|
\widehat{h}^{x',t}
\big(
\xi',\tau
\big)
\big|^2
d\tau
d\xi'
\doteq
\|h\|_{\overset{\dott}{\mathcal{B}^s}}^2,
\end{align}
which is  the optimal homogeneous estimate \eqref{B-est-redu-pure-optimal} for $v_r$ stated in Remark \ref{optimal-homogeneous-est}.
Finally, using $ |\tau+\xi'^2|^{\frac{1}{2}} \le (1+|\tau+\xi'^2|)^{\frac{1}{2}}$
we get  the desired inhomogeneous  bound \eqref{Vr-B-est}  for $V_r$
\begin{align*}
%\label{inhom-Bs-est-vr}
\|
\psi_4
\cdot
V_r
\|_{X^{s,b}(\rr^{n+1})}^2
\lesssim
\int_{\rr^{n}}
\big(
1+\xi'^2+|\tau+\xi'^2|
\big)^{s}
\big(1+|\tau+\xi'^2|\big)^{1/2}
\big|
\widehat{h}^{x',t}
\big(
\xi',\tau
\big)
\big|^2
d\tau
d\xi'
=
\|h\|_{\mathcal{B}^s}^2.
\end{align*}
{\bf Estimate \eqref{B-est-redu-pure} for $v_i$.}
We split $v_i(x',x_n,t)$ as the sum of two integrals, 
one for $\xi_n$ near $0$ and the other away from $0$,
that is we write
\begin{align}
\label{pure-sln-vi-split}
v_i(x',x_n,t)
\simeq
v_0(x',x_n,t)
+
v_1(x',x_n,t),
\quad
x'\in\rr^{n-1},
\,\,
x_n\in\rr^+,
\quad
t\in[0,2],
\end{align}
where the functions $v_0$ and $v_1$ are defined as follows
\begin{align}
\label{def-v0}
v_0(x',x_n,t)
\doteq&
\int_{\rr^{n-1}} 
\int_{\xi_n=0}^1
e^{i\xi'\cdot x'-\xi_nx_n -i(\xi'^2-\xi_n^2)t}
\cdot
\xi_n 
\cdot
\widehat{h}^{x',t}
\big(
\xi',-\xi'^2+\xi_n^2
\big)
d\xi_n
d\xi',
\\
\label{def-v1}
v_1(x',x_n,t)
\doteq&
\int_{\rr^{n-1}} 
\int_{\xi_n=1}^\infty
e^{i\xi'\cdot x'-\xi_nx_n -i(\xi'^2-\xi_n^2)t}
\cdot
\xi_n 
\cdot
\widehat{h}^{x',t}
\big(
\xi',-\xi'^2+\xi_n^2
\big)
d\xi_n
d\xi'.
\end{align}
\underline{\it Proof of  estimate \eqref{B-est-redu-pure} for $v_0$.}
This proof
follows from the $L^2$-boundedness of the Laplace transform.
We begin with extending $v_0$ from $\rr^{n-1}\times\rr^+\times[0,2]$ to $\rr^{n-1}\times\rr\times\rr$ and we denote this extension by
\begin{equation}
\label{v0-extension}
V_0(x',x_n,t)
\doteq
\int_{\rr^{n-1}} 
\int_{\xi_n=0}^1
e^{i\xi'\cdot x'-\xi_n\varphi_1(x_n) -i(\xi'^2-\xi_n^2)t}
\xi_n 
\widehat{h}^{x',t}
\big(
\xi',-\xi'^2+\xi_n^2
\big)
d\xi_n
d\xi',
\,\,
(x',x_n,t)\in\rr^{n+1},
\end{equation}
where the function  $\varphi_1(x_n)$ is a smooth version of $|x_n|$ and its graph is shown in Figure \ref{fig:phi1}.
\vskip-0.05in
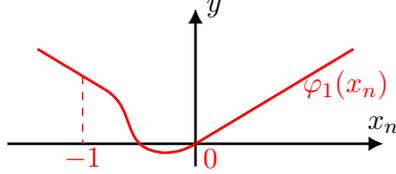
\begin{minipage}{0.49\linewidth}
\begin{center}
\begin{tikzpicture}[yscale=0.6, xscale=1]
%%%%%%%%%%%%%%%%%%%
%
%Variables defined
%
%%%%%%%%%%%%%%%%%%%%
\newcommand\X{0};
\newcommand\Y{0};
\newcommand\FX{11};
\newcommand\FY{11};
\newcommand\R{0.6};
\newcommand*{\TickSize}{2pt};
%%%%%%%%%%%%%%%%%%
%
%End
%
%%%%%%%%%%%%%%%%%%%%%
\draw[black,line width=1pt,-{Latex[black,length=2mm,width=2mm]}]
(-2.5,0)
--
(2.5,0)
node[above]
{\fontsize{\FX}{\FY}\bf \textcolor{black}{$x_n$}};

\draw[black,line width=1pt,-{Latex[black,length=2mm,width=2mm]}]
(0,-0.5)
--
(0,3)
node[right]
{\fontsize{\FX}{\FY}\bf \textcolor{black}{$y$}};

\draw[line width=1pt, yscale=1,domain=-2.1:-1.2,smooth,variable=\x,red]  plot ({\x},{-\x});

\draw[line width=1pt, yscale=1,domain=0:2.1,smooth,variable=\x,red]  plot ({\x},{\x});

\draw[smooth,line width=1pt, red]
(-0.4,-0.2)
to[out=5,in=-135]
(0,0)
;

\draw[smooth,line width=1pt, red]
(-1.2,1.2)
to[out=-45,in=175]
(-0.4,-0.2)
;

\draw[red,dashed, line width=0.5pt]
(2,1.3)
node[]
{\fontsize{\FX}{\FY}$\varphi_1(x_n)$}

(0,0)
node[yshift=-0.2cm,xshift=0.2cm]
{\fontsize{\FX}{\FY}$0$}

(-1.5,1.5)
--
(-1.5,0)
node[yshift=-0.2cm]
{\fontsize{\FX}{\FY}$-1$};

\end{tikzpicture}

\vskip-0.15in
\captionof{figure}{}
\label{fig:phi1}
\end{center}
\end{minipage}
\begin{minipage}{0.48\linewidth}
\begin{equation*}
%\label{phi1-x-even-ext}
\varphi_1(x_n)
=
\begin{cases}
x_n,
\quad
\,\,\,\,
x_n\geq
0,
\\
-x_n,
\quad
x_n\leq
-1,
\\
\text{ smooth on } \rr.
\quad
\end{cases}
\end{equation*}
\end{minipage}

\vskip0.05in
\noindent
Since $\psi_4(t)V_0=v_0$ for $(x',x_n,t)\in \rr^{n-1}\times\rr^+\times[0,2]$, where $\psi_4=\psi(t/4)$ and $\psi$ is the time localizer given by \eqref{time-localizer},
 in order to show estimate \eqref{B-est-redu-pure} for $v_0$, it suffices to show that 
\begin{equation}
\label{V0-B-est}
\|\psi_4(t)V_0\|_{X^{s,b}(\rr^{n+1})}
\lesssim
\|h\|_{\mathcal{B}^s(\rr^{n})},
\quad
s\in\rr
\quad
\text{and}
\quad
b\ge 0.
\end{equation}
For this we  separate  $\xi'$ and $\xi_n$ 
in the Bourgain norm multipliers via the basic inequality
\begin{align*}
%\label{sb-term-est}
(1+|\xi|)^{2s}
(1+|\tau+\xi^2|)^{2b}
\lesssim
(1+|\xi'|)^{2s}
(1+|\tau+\xi'^2|)^{2b}
(1+|\xi_n|)^{2|s|+4b},
\end{align*}
and thus we bound $\|\psi_4(t)V_0\|_{X^{s,b}}$
as follows
\begin{align*}
%\label{V0-B-norm-separated}
\|\psi_4(t)V_0\|^2_{X^{s,b}}
\lesssim
\int_{\rr^{n+1}}
(1+|\xi'|)^{2s}
(1+|\tau+\xi'^2|)^{2b}
(1+|\xi_n|)^{2|s|+4b}
|\widehat{\psi_4(t)V_0}
(\xi',\xi_n,\tau)|^2
d\xi_n
d\xi'
d\tau.
\end{align*}
Furthermore, we split the $\xi_n$-integration by using the inequality
\begin{align*}
%\label{sb-term-est-1}
(1+|\xi_n|)^{2|s|+4b}
\lesssim
\chi_{|\xi_n|\le 1} (\xi_n)
+
\chi_{|\xi_n|> 1} (\xi_n)
\cdot
|\xi_n|^{2|s|+4b}
\lesssim
1+|\xi_n|^{2(\lfloor |s|+2b\rfloor+1)},
\end{align*}
and obtain the following bound
\begin{align*}
\|\psi_4(t)V_0\|^2_{X^{s,b}(\rr^{n+1})}
\lesssim&
\int_{\rr^n}
(1+|\xi'|)^{2s}
(1+|\tau+\xi'^2|)^{2b}
\Big(
\int_\rr
|\widehat{\psi_4(t)V_0}
(\xi',\xi_n,\tau)|^2
d\xi_n
\Big)
d\xi'
d\tau 
\\
+&
\int_{\rr^n}
(1+|\xi'|)^{2s}
(1+|\tau+\xi'^2|)^{2b}
\Big(
\int_\rr
|\xi_n|^{2(\lfloor |s|+2b\rfloor+1)}
|\widehat{\psi_4(t)V_0}
(\xi',\xi_n,\tau)|^2
d\xi_n
\Big)
d\xi'
d\tau.
\end{align*}
Also, using Plancherel's theorem in the $\xi_n$-integral, 
we go to the physical variable $x_n$ and 
we get  
\begin{align}
\label{V0-L2-bound}
\hskip-0.1in
\|\psi_4V_0\|^2_{X^{s,b}}
\lesssim&
\int_{\rr^n}
\hskip-0.08in
(1
\hskip-0.01in+\hskip-0.01in
|\xi'|)^{2s}
(1
\hskip-0.01in+\hskip-0.01in
|\tau
\hskip-0.01in+\hskip-0.01in
\xi'^2|)^{2b}
\Big(
\|
\widehat{\psi_4V_0}^{x',t}
(\xi',\cdot,\tau)\|_{L_{x_n}^2}^2
\hskip-0.04in+\hskip-0.01in
\|
\p_{x_n}^{k}
\widehat{\psi_4V_0}^{x',t}
(\xi',\cdot,\tau)\|_{L_{x_n}^2}^2
\Big)
d\xi'
d\tau,
\end{align}
where $k=\lfloor |s|+2b\rfloor+1$. For $\|
\p_{x_n}^{k}
\widehat{\psi_4(t)V_0}^{x',t}
(\xi',\cdot,\tau)\|_{L_{x_n}^2}^2$, $k\in\NN_0$, we have the following result.
\begin{lemma}
\label{near0-multiplier-lem}
For any $k\in\NN_0$, we have 
\begin{align}
\label{near0-multiplier-est}
\|
\p_{x_n}^{k}
\widehat{\psi_4(t)V_0}^{x',t}
(\xi',\cdot,\tau)\|_{L_{x_n}^2}^2
\lesssim
\int_0^1
\Big|
\widehat{\psi}_4(\tau+\xi'^2-\xi_n^2)
\xi_n 
\widehat{h}^{x',t}
\big(
\xi',-\xi'^2+\xi_n^2
\big)
\Big|^2
d\xi_n.
\end{align}
\end{lemma}
\noindent
We prove Lemma \ref{near0-multiplier-lem}  later.
%using the Laplace transform $L^2$-boundedness.
Now, 
using it, from inequality \eqref{V0-L2-bound} we get
\begin{align*}
\|\psi_4(t)V_0\|^2_{X^{s,b}}
\lesssim
\int_{\rr^n}
(1+|\xi'|)^{2s}
(1+|\tau+\xi'^2|)^{2b}
\Big(
\int_0^1
\Big|
\widehat{\psi}_4(\tau+\xi'^2-\xi_n^2)
\xi_n 
\widehat{h}^{x',t}
\big(
\xi',-\xi'^2+\xi_n^2
\big)
\Big|^2
d\xi_n
\Big)
d\xi'
d\tau.
\end{align*}
Furthermore, applying Fubini's theorem 
we compute the $d\tau$-integral first,
and since  $0\le \xi_n\le 1$, we have 
$
\int_{\rr}
(1+|\tau+\xi'^2|)^{2b}
|\widehat{\psi}_4(\tau+\xi'^2-\xi_n^2)|^2
d\tau
\lesssim
\|\psi_4\|_{H^b}^2
\lesssim
1.
$
This gives us that 
\begin{align*}
\|\psi_4(t)V_0\|^2_{X^{s,b}(\rr^{n+1})}
\lesssim
\int_{\rr^{n-1}}
\int_{\xi_n=0}^1
(1+|\xi'|+\xi_n)^{2s}
\cdot
\xi_n^2
|\widehat{h}^{x',t}
\big(
\xi',-\xi'^2+\xi_n^2
\big)|^2
d\xi_n
d\xi'.
\end{align*}
Finally, making  change of variables $\xi'=\xi'$,  $\tau=-\xi'^2+\xi_n^2$, with the Jacobian 
$J=\frac{\partial (\xi', \tau)}{\partial (\xi', \xi_n)}=2\xi_n$,
we get
\begin{align}
\label{Bs-est-v0}
\hskip-0.1in
\|\psi_4V_0\|^2_{X^{s,b}}
\lesssim&
\int_{\rr^{n-1}}
\hskip-0.03in
\int_{\tau=-\xi'^2}^{-\xi'^2+1}
\hskip-0.08in
\big(
1
\hskip-0.01in+\hskip-0.01in
[\xi'^2
\hskip-0.01in+\hskip-0.01in
|\tau
\hskip-0.01in+\hskip-0.01in
\xi'^2|]^{\frac 12}
\big)^{2s}
|\tau
\hskip-0.01in+\hskip-0.01in
\xi'^2|^{\frac 12}
\big|
\widehat{h}^{x',t}
\big(
\xi',\tau
\big)
\big|^2
d\tau
d\xi'
\lesssim
\|h\|_{\overset{\dott}{\mathcal{B}^s}}^2
\lesssim
\|h\|_{\mathcal{B}^s}^2,
\end{align}
for
$s\in\rr$
and 
$b\ge 0$.
Note, in the second to last step we proceed like in the derivation of inequality \eqref{optimal-hom-est}.
This  implies  the desired inequality \eqref{V0-B-est} for $V_0$. Also, it gives us inequality \eqref{B-est-redu-pure-optimal} in Remark \ref{optimal-homogeneous-est} for $v_0$ once we prove Lemma \ref{near0-multiplier-lem}.
%We also  obtain the desired estimate \eqref{B-est-redu-pure} of $v_0$.
\,\,
$\square$

\noindent
{\bf Proof of Lemma \ref{near0-multiplier-lem}.} We begin with computing 
$\partial_{x_n}^k \widehat{\psi_4 V_0}^{x',t}$. 
From definition  \eqref{v0-extension}, taking the Fourier transform  with respect to $x'$ and $t$, and then differentiating $k$ times with respect to $x_n$, we obtain
\begin{equation}
\label{n-der-psi-V}
\p_{x_n}^k
\widehat{\psi_4(t)V_0}^{x',t}
(\xi',x_n,\tau)
\simeq
\int_{\xi_n=0}^1
\p_{x_n}^k
[e^{-\xi_n\varphi_1(x_n)}]
\widehat{\psi}_4(\tau+\xi'^2-\xi_n^2)
\xi_n
\widehat{h}^{x',t}
\big(
\xi',-\xi'^2+\xi_n^2
\big)
d\xi_n.
\end{equation}
Furthermore,
we split the $L^2_{x_n}$-norm of $\p_{x_n}^{k}
\widehat{\psi_4(t)V_0}^{x',t}
(\xi',x_n,\tau)$ for $x_n$ near  and away from $0$,
that is
\begin{align*}
\|
\p_{x_n}^k
\widehat{\psi_4(t)V_0}^{x',t}
(\xi',\cdot,\tau)\|_{L_{x_n}^2}^2
\doteq
\int_\rr
\big|
\p_{x_n}^k
\widehat{\psi_4(t)V_0}^{x',t}
(\xi',x_n,\tau)
\big|^2
dx_n
=
I_0+I_{-\infty}+I_{\infty},
\end{align*}
where
$
I_0
\doteq
\int_{-1}^1
|
\p_{x_n}^k
\widehat{\psi_4(t)V_0}^{x',t}
(\xi',x_n,\tau)
|^2
dx_n
$,
$
I_{\infty}
\doteq
\int_{1}^\infty
|
\p_{x_n}^k
\widehat{\psi_4(t)V_0}^{x',t}
(\xi',x_n,\tau)
|^2
dx_n
$
and
$
I_{-\infty}
\doteq
\int_{-\infty}^{-1}
|
\p_{x_n}^k
$
$
\widehat{\psi_4(t)V_0}^{x',t}
(\xi',x_n,\tau)
|^2
dx_n
$.
Since the estimate for $I_{-\infty}$ is similar to the estimate of $I_{\infty}$, here 
we only estimate the integrals $I_0$ and $I_\infty$.

\vskip0.05in\nin
{\it Estimation of $I_0$.} For this integral, taking the 
sup-norm for $x_n$ and using the Cauchy-Schwarz inequality for the integral of $\xi_n$ 
we get the desired estimate \eqref{near0-multiplier-est} for $I_0$.

\vskip0.05in\nin
{\it Estimation of $I_\infty$.} For $x_n>1$, we have 
$\varphi_1(x_n)=x_n$ and integral \eqref{n-der-psi-V} becomes
$$
\p_{x_n}^k
\widehat{\psi_4(t)V_0}^{x',t}
(\xi',x_n,\tau)
\simeq
\int_{\xi_n=0}^1
e^{-\xi_nx_n}
(-\xi_n)^k
\cdot
\widehat{\psi}_4(\tau+\xi'^2-\xi_n^2)
\xi_n
\widehat{h}^{x',t}
\big(
\xi',-\xi'^2+\xi_n^2
\big)
d\xi_n.
$$
Thus, we bound $I_{\infty}$ as follows
\begin{align*}
%\label{est-der-I1-term2-no2}
I_{\infty}
\lesssim&
\int_{1}^\infty
\Big(
\int_0^\infty
e^{- \xi_n x_n}
\cdot
\chi_{0<\xi_n<1}(\xi_n)
\big|
\widehat{\psi}_4(\tau+\xi'^2-\xi_n^2)
\xi_n
\widehat{h}^{x',t}
(
\xi',-\xi'^2+\xi_n^2
)
\big|
d\xi_n
\Big)^2
dx_n,
\end{align*}
and we are reduced to the  $L^2$-boundedness of 
the Laplace transform stated below
\cite{fhm2017,hardy1933}.
\begin{lemma}
[$L^2$ boundedness of the Laplace transform] 
\label{Laplace-lemma}
Suppose that $Q(\xi)\in L_\xi^2(0,\infty)$. Then, the map
$
Q(\xi)
\longmapsto
\int_0^\infty e^{-\xi y}Q(\xi)d\xi
$
is bounded from $L_\xi^2(0,\infty)$ into $L_y^2(0,\infty)$  with
\begin{align}
\label{L-T-inequ}
\Big\|\int_0^\infty e^{-\xi y}Q(\xi)d\xi\Big\|_{L_y^2(0,\infty)}
\leq
\sqrt{\pi}\|Q(\xi)\|_{L_\xi^2(0,\infty)}.
\end{align}
\end{lemma}
\nin
In fact, extending the $x_n$-integration from $(1,\infty)$ to $(0,\infty)$,
applying the estimate \eqref{L-T-inequ} 
with 
$
Q=\chi_{0<\xi_n<1}(\xi_n)
|\widehat{\psi}_4(\tau+\xi'^2-\xi_n^2)
\xi_n
\widehat{h}^{x',t}
(
\xi',-\xi'^2+\xi_n^2
)|,
$ 
we get
the desired inequality \eqref{near0-multiplier-est} for $I_{\infty}$.
This completes the proof of Lemma \ref{near0-multiplier-lem}.
\,\,$\Box$

\vskip0.05in
\noindent
\underline{\it Estimate \eqref{B-est-redu-pure} for $v_1$.}
Using the elementary but 
crucial identity
$\p_{x_n}e^{-\xi_n x_n}=-\xi_ne^{-\xi_n x_n}$ and
moving $\p_{x_n}$ outside the $d\xi_n$-integral, we rewrite $v_1$ as follows
\begin{align*}
%\label{v1-ext-fin}
v_1(x',x_n,t)
=
-\p_{x_n}
\int_{\xi'\in\rr^{n-1}} 
\int_{\xi_n=1}^\infty
e^{i\xi'\cdot x'-\xi_nx_n -i(\xi'^2-\xi_n^2)t}
\widehat{h}^{x',t}
\big(
\xi',-\xi'^2+\xi_n^2
\big)
d\xi_n
d\xi',
\quad
x_n>0.
\end{align*}
Since the integrand is not oscillatory in $\xi_n$, to estimate 
the Bourgain norm of $v_1$ we need to extend 
$e^{-\xi_n x_n}$ as a function of $x_n$ from $[0, \infty)$ 
to $\rr$ for each $\xi_n>0$.
%
%Like in \cite{et2016},
For this we use the one-sided  cut-off function $\rho(x)$, 
which is displayed 
 in Figure \ref{fig:pho}.

\begin{minipage}{0.6\linewidth}
\begin{center}
\begin{tikzpicture}[yscale=0.5, xscale=0.8]
%%%%%%%%%%%%%%%%%%%
%
%Variables defined
%
%%%%%%%%%%%%%%%%%%%%
\newcommand\X{0};
\newcommand\Y{0};
\newcommand\FX{11};
\newcommand\FY{11};
\newcommand\R{0.6};
\newcommand*{\TickSize}{2pt};
%%%%%%%%%%%%%%%%%%
%
%End
%
%%%%%%%%%%%%%%%%%%%%%
\draw[black,line width=1pt,-{Latex[black,length=2mm,width=2mm]}]
(-5,0)
--
(5,0)
node[above]
{\fontsize{\FX}{\FY}\bf \textcolor{black}{$x$}};

\draw[black,line width=1pt,-{Latex[black,length=2mm,width=2mm]}]
(0,0)
--
(0,3)
node[right]
{\fontsize{\FX}{\FY}\bf \textcolor{black}{$\rho$}};

\draw[line width=1pt, yscale=2,domain=-1.5:-4.3,smooth,variable=\x,red]  plot ({\x},{0});

\draw[line width=1pt, yscale=2,domain=0:4.3,smooth,variable=\x,red]  plot ({\x},{1});

\draw[smooth,line width=1pt, red]
(0,2)
to[out=-170,in=10]
(-1.5,0)
;

\draw[red]
(2,2.5)
node[]
{\fontsize{\FX}{\FY}$\rho(x)$}

(0,0)
node[yshift=-0.2cm]
{\fontsize{\FX}{\FY}$0$}

(-1.5,0)
node[yshift=-0.2cm]
{\fontsize{\FX}{\FY}$-1$};

\end{tikzpicture}
\end{center}
\vskip-0.08in
\captionof{figure}{}
\label{fig:pho}
\end{minipage}
\hskip-0.88in
\begin{minipage}{0.5\linewidth}
\begin{equation*}
%\label{rho-def}
\rho(x)
=
\begin{cases}
1,
&\quad
x\ge0,
\\
0,
&\quad
x\le -1.
\end{cases}
\end{equation*}
\end{minipage}

\nin
Now, we estimate the restricted Bourgain norm of $v_1$
by using its extension  to  $\rr^{n+1}$
\begin{align}
\label{v1-ext}
V_1(x',x_n,t)
=
-\p_{x_n}
\int_{\xi'\in\rr^{n-1}} 
\int_{\xi_n=1}^\infty
e^{i\xi'\cdot x'-\xi_nx_n -i(\xi'^2-\xi_n^2)t}
\cdot
\rho(\xi_nx_n)
\cdot
\widehat{h}^{x',t}
\big(
\xi',-\xi'^2+\xi_n^2
\big)
d\xi_n
d\xi'.
\end{align}
Thus,  to prove estimate \eqref{B-est-redu-pure} for $v_1$, it suffices to show that 
\begin{equation}
\label{V1-B-est}
\|V_1\|_{X^{s,b}(\rr^{n+1})}
\lesssim
\|h\|_{\mathcal{B}^s},
\quad
s>-3/2
\quad
\text{and}
\quad
0\le b<1/2.
\end{equation}
For this we compute the  full Fourier transform of $V_1$.
Defining the  Schwartz function 
%Schwarz function
\begin{equation}
\label{def-eta}
\eta(x_n)
\doteq
e^{-x_n} \rho(x_n),
\end{equation}
we rewrite  extension \eqref{v1-ext} as follows
\begin{align}
\label{v1-ext-1}
V_1(x',x_n,t)
=
-
\p_{x_n}
\int_{\xi'\in\rr^{n-1}} 
\int_{\xi_n=1}^\infty
e^{i\xi'\cdot x'-i(\xi'^2-\xi_n^2)t}
\cdot
\eta(\xi_nx_n)
\cdot
\widehat{h}^{x',t}
\big(
\xi',-\xi'^2+\xi_n^2
\big)
d\xi_n
d\xi'.
\end{align}
Furthermore, for the $\xi_n,\xi'$-integration in \eqref{v1-ext-1},
making the change of variables $\xi'=\xi'$, $\tau=-(\xi'^2-\xi_n^2)$ or $\xi_n=(\tau+\xi'^2)^{1/2}$, with Jacobian 
$J=\partial (\xi', \tau)/\partial (\xi', \xi_n)=2\xi_n$,
we get
\begin{align}
\label{eqn-V1-1}
V_1(x',x_n,t)
\simeq
\p_{x_n}
\int_{\rr^n}
e^{i\xi'\cdot x' +i\tau t}
\chi_{\tau+ \xi'^2\ge 1}(\xi',\tau)
\cdot
\eta
\Big(
(\tau+\xi'^2)^{1/2}x_n
\Big)
(\tau+\xi'^2)^{-\frac12}
\widehat{h}^{x',t}
\big(
\xi',\tau
\big)
d\tau
d\xi',
\end{align}
where $\chi_{\tau+\xi'^2\ge 1}(\xi',\tau)$  is the characteristic function on the domain $\{(\xi',\tau)\in\rr^n:\tau+ \xi'^2\ge 1\}$. 
Using formula \eqref{eqn-V1-1}, we find 
that the full Fourier transform of $V_1$ is 
\begin{align}
\label{v1-xt-FT}
\widehat{V}_1(\xi',\xi_n,\tau)
\simeq
\chi_{\tau+ \xi'^2\ge1}(\xi',\tau)
\cdot
\xi_n
F(\xi',\xi_n,\tau)
\frac{
\widehat{h}^{x',t}
\big(
\xi',\tau
\big)
}
{(\tau+\xi'^2)^{\frac12}},
\end{align}
where $F(\xi',\xi_n,\tau)$ is given by
\begin{align}
\label{F-def-reduced-ibvp}
F(\xi',\xi_n,\tau)
\doteq
\int_{\rr}
e^{-i\xi_n x_n}
\eta
\big(
(\tau+\xi'^2)^{1/2}x_n
\big)
dx_n.
\end{align}
Then, the Bourgain norm $\|V_1\|_{X^{s,b}}$ is 
equal to the quantity
\begin{align}
\label{X-norm-v1}
\|V_1\|_{X^{s,b}}^2
\hskip-0.03in
\simeq
\hskip-0.08in
\int_{\tau
\hskip-0.01in+\hskip-0.01in
 \xi'^2\ge1}
\Big[
\int_{\rr}
(1
\hskip-0.01in+\hskip-0.01in
\xi'^2
\hskip-0.01in+\hskip-0.01in
\xi_n^2)^{s}
(1
\hskip-0.01in+\hskip-0.01in
|\tau
\hskip-0.01in+\hskip-0.01in
\xi'^2
\hskip-0.01in+\hskip-0.01in
\xi_n^2|)^{2b}
|
\xi_n
F(\xi',\xi_n,\tau)
|^2
d\xi_n
\Big]
\Big|
\frac{
\widehat{h}^{x',t}
\big(
\xi',\tau
\big)
}
{(\tau+\xi'^2)^{\frac12}}
\Big|^2
d\tau
d\xi',
\end{align}
which we estimate by first bounding 
the $\xi_n$ integral in \eqref{X-norm-v1} by a multiplier 
 as follows.
\begin{lemma}
\label{v1-mult-lema}
Let $F$ be given by \eqref{F-def-reduced-ibvp}
and $\tau+\xi'^2\ge1$. 
If  $b\ge 0$,  $\varepsilon>0$, and $s\geq -3/2-\varepsilon$,
then
\begin{align}
\label{v1-mult-ine-1}
&\int_{\rr}
(1+\xi'^2+\xi_n^2)^{s}
(1+|\tau+\xi'^2+\xi_n^2|)^{2b}
|
\xi_n
F(\xi',\xi_n,\tau)
|^2
d\xi_n
\lesssim
(1+\xi'^2+|\tau+\xi'^2|)^{s}|\tau+\xi'^2|^{\frac{4b+1+2\varepsilon}{2}}.
\end{align}
\end{lemma}
\noindent
We shall prove Lemma \ref{v1-mult-lema} later. 
Now, we use it to prove inequality \eqref{V1-B-est}.
Since by our assumption $s>-3/2$, we can apply 
inequality \eqref{v1-mult-ine-1} for any $\varepsilon>0$
and $b \ge 0$ to get
\begin{align}
\label{Bs-est-v1}
\|V_1\|_{X^{s,b}(\rr^{n+1})}^2
\lesssim
\int_{\rr^n}
\chi_{\tau+ \xi'^2\ge1}(\xi',\tau)
\cdot
(1+\xi'^2+|\tau+\xi'^2|)^{s}|\tau+\xi'^2|^{\frac{4b-1+2\varepsilon}{2}}
|
\widehat{h}^{x',t}
\big(
\xi',\tau
\big)
|^2
d\tau d\xi',
\end{align}
which gives the desired inequality \eqref{V1-B-est} 
for $V_1$, if we can choose $b$ and $\varepsilon$ 
such that 
$\frac{4b-1+2\varepsilon}{2}
=
\frac12,
$
or 
$
 b
=
 \frac12-\frac12\varepsilon.
 $
 Since we also have $b\ge 0$ we get the following condition 
 for $b$ in the case of ibvp 
\begin{equation}
\label{b-ibvp-condition}
0\le b< 1/2,
 \end{equation}
 unlike in the ivp case where $\frac 12<b<1$.
 Furthermore, for any  $b$ satisfying condition \eqref{b-ibvp-condition} we can have
$\frac{4b-1+2\varepsilon}{2}
=
\frac12,
$
by choosing 
\begin{equation}
\label{epsilon-b-condition}
\varepsilon
=
1-2b.
 \end{equation}
 This completes the proof  estimates
  \eqref{B-est-redu-pure-optimal} and  \eqref{B-est-redu-pure} for $v_1$,  once we prove Lemma \ref{v1-mult-lema}.
\,\,
$\square$

\vskip.05in
\nin
{\bf Proof of Lemma \ref{v1-mult-lema}.} 
It is based on the  following decay estimate in 
$\xi_n$, which is a consequence of the fact that 
the function $\eta$, defined  by \eqref{def-eta},  is in 
the Schwartz space. 
\begin{lemma}
\label{F-bound-lem}
Let $F$ be given by \eqref{F-def-reduced-ibvp}.
For $\tau+\xi'^2> 1$
and  for any
$m\ge 0$,   we have
\begin{align}
\label{F-bound-ine}
|F(\xi',\xi_n,\tau)|
\leq
c_{\rho,m}
\cdot
\frac{|\tau+\xi'^2|^{\frac{m-1}{2}}}{|\xi_n|^m+|\tau+\xi'^2|^{\frac{m}{2}}},
\end{align}
where $c_{\rho ,m}$ is a constant only depending on $\rho$ (our spatial extension cutoff of $v_1$) and $m$.
\end{lemma}
Next, we  estimate the integrand
appearing in inequality \eqref{v1-mult-ine-1}, 
 which we denote as follows
\begin{align*}
%\label{I-def-vi}
I(\xi',\xi_n,\tau)
\doteq
(1+\xi'^2+\xi_n^2)^{s}
(1+|\tau+\xi'^2+\xi_n^2|)^{2b}
|
\xi_n
F(\xi',\xi_n,\tau)
|^2.
\end{align*}
Using decay estimate  \eqref{F-bound-ine}, we get
\begin{align*}
I(\xi',\xi_n,\tau)
\lesssim
(1+\xi'^2+\xi_n^2)^{s}
(1+|\tau+\xi'^2+\xi_n^2|)^{2b}
\xi_n^2
\frac{|\tau+\xi'^2|^{m-1}}{|\xi_n|^{2m}+|\tau+\xi'^2|^{m}},
\quad
m\ge0.
\end{align*}
Also,
since $\tau+\xi'^2\ge 1$ and $b\ge 0$,
separating $\xi_n^2$ and $\tau+\xi'^2$,  we get 
\begin{equation}
\label{L2-est-mult}
I(\xi',\xi_n,\tau)
\lesssim
(1+\xi'^2+\xi_n^2)^{s}
\big[
|\tau+\xi'^2|^{2b}
+
|\xi_n|^{4b}
\big]
\xi_n^2
\,
\frac{|\tau+\xi'^2|^{m-1}}{|\xi_n|^{2m}+|\tau+\xi'^2|^{m}},
\quad
m\ge 0.
\end{equation}
Now, we continue our estimation by considering  the following two cases.

\vskip.05in
\noindent
$\bullet$
Case 1: $s\ge 0$
\quad
$\bullet$
Case 2:
 $s<0$

\vskip0.05in
\noindent
{\bf  Case 1.}  To prove estimate \eqref{v1-mult-ine-1} in this case, we consider the following two subcases.

\vskip.05in
\noindent
$\bullet$
Subcase 1-1: $|\tau+\xi'^2|\ge \xi_n^2$
\quad
$\bullet$
Subcase 1-2:
 $|\tau+\xi'^2|< \xi_n^2$

\vskip0.05in
\noindent
\underline{Subcase 1-1.} Since $s\ge 0$ and $|\tau+\xi'^2|\ge \xi_n^2$, we have 
$
%\label{multiplier-bound-1-1}
(1+\xi'^2+\xi_n^2)^{s}
\lesssim
(1+\xi'^2+|\tau+\xi'^2|)^{s}
$
and
$
\big[
|\tau+\xi'^2|^{2b}
+
|\xi_n|^{4b}
\big]
\lesssim
|\tau+\xi'^2|^{2b},
$
which combined with 
inequality \eqref{L2-est-mult}  gives us the estimate for $I$ 
\begin{align}
\label{integrand-bnd}
I(\xi',\xi_n,\tau)
\lesssim
(1+\xi'^2+|\tau+\xi'^2|)^{s}
|\tau+\xi'^2|^{2b}
%\cdot
\xi_n^2
\,
\frac{|\tau+\xi'^2|^{m-1}}{|\xi_n|^{2m}+|\tau+\xi'^2|^{m}}.
\end{align}
Furthermore, choosing  
$m=3/2+\varepsilon$ 
where 
$
\varepsilon>0
$
and using $|\tau+\xi'^2|\ge 1$, from \eqref{integrand-bnd} we get
$$
I(\xi',\xi_n,\tau)
\lesssim
(1+\xi'^2+|\tau+\xi'^2|)^{s}|\tau+\xi'^2|^{2b}
\xi_n^2
\,
\frac{|\tau+\xi'^2|^{\frac12+\varepsilon}}{(|\xi_n|+1)^{3+2\varepsilon}}
\lesssim
(1+\xi'^2+|\tau+\xi'^2|)^{s}|\tau+\xi'^2|^{\frac{4b+1+2\varepsilon}{2}}
\,
\frac{1}{(|\xi_n|+1)^{1+2\varepsilon}}.
$$
Therefore, integrating $I$ with respect to $\xi_n$ we get
the desired estimate \eqref{v1-mult-ine-1} in this subcase.

\vskip0.05in
\noindent
\underline{Subcase 1-2.} For this, we consider the following two subcases.

\vskip.05in
\noindent
$\bullet$
Subcase 1-2-1: $|\xi'|\ge |\xi_n|$
\quad
$\bullet$
Subcase 1-2-2:
 $|\xi'|< |\xi_n|$

\vskip0.05in
\noindent
{\it Subcase 1-2-1.}
Since $s\ge 0$, $1\le |\tau+\xi'^2|< \xi_n^2$ and $|\xi'|\ge |\xi_n|$, we have 
$
(1+\xi'^2+\xi_n^2)^{s}
\lesssim
(1+\xi'^2)^{s}
\lesssim
(1+\xi'^2+|\tau+\xi'^2|)^s
$
and
$
\big[
|\tau+\xi'^2|^{2b}
+
|\xi_n|^{4b}
\big]
\lesssim
|\xi_n|^{4b}.
$
This combined with inequality \eqref{L2-est-mult}  gives us that
\begin{align}
\label{integrand-bnd-1-2-1}
I(\xi',\xi_n,\tau)
\lesssim
(1+\xi'^2+|\tau+\xi'^2|)^{s}
|\xi_n|^{4b}
%\cdot
\xi_n^2
\,
\frac{|\tau+\xi'^2|^{m-1}}{|\xi_n|^{2m}}.
\end{align}
Here, choosing $m=2b+3/2+\varepsilon$, from \eqref{integrand-bnd-1-2-1} we get
$$
I(\xi',\xi_n,\tau)
\lesssim
(1+\xi'^2+|\tau+\xi'^2|)^{s}|\xi_n|^{4b}
\xi_n^2
\,
\frac{|\tau+\xi'^2|^{2b+\frac12+\varepsilon}}{|\xi_n|^{4b+3+2\varepsilon}}
=
(1+\xi'^2+|\tau+\xi'^2|)^{s}|\tau+\xi'^2|^{\frac{4b+1+2\varepsilon}{2}}
\,
\frac{1}{|\xi_n|^{1+2\varepsilon}}.
$$
Now,  integrating $I$ with respect to $\xi_n$ for $|\xi_n|\ge 1$, we get the desired estimate \eqref{v1-mult-ine-1} in this subcase.

\vskip0.05in
\noindent
{\it Subcase 1-2-2.}
Since $s\ge 0$, $1\le |\tau+\xi'^2|< \xi_n^2$ and $|\xi'|< |\xi_n|$, we have 
$
(1+\xi'^2+\xi_n^2)^{s}
\lesssim
|\xi_n|^{2s}
$
and
$
\big[
|\tau+\xi'^2|^{2b}
+
|\xi_n|^{4b}
\big]
\lesssim
|\xi_n|^{4b}.
$
Using these inequalities, 
 from inequality \eqref{L2-est-mult}  we get 
\begin{align*}
%\label{integrand-bnd-1-2-2}
I(\xi',\xi_n,\tau)
\lesssim
|\xi_n|^{2s}
%\cdot
|\xi_n|^{4b}
%\cdot
\xi_n^2
\,
\frac{|\tau+\xi'^2|^{m-1}}{|\xi_n|^{2m}+|\tau+\xi'^2|^{m}}
\lesssim
|\tau+\xi'^2|^{\frac{2s+4b+1+2\varepsilon}{2}}
\,
\frac{1}{|\xi_n|^{1+2\varepsilon}},
\end{align*}
where in the last step we choose $m=s+2b+3/2+\varepsilon.$
Furthermore, since $s\ge 0$, we get 
$$
I(\xi',\xi_n,\tau)
\lesssim
(1+|\tau+\xi'^2|)^{s}
|\tau+\xi'^2|^{\frac{4b+1+2\varepsilon}{2}}
\,
\frac{1}{|\xi_n|^{1+2\varepsilon}}
\lesssim
(1+\xi'^2+|\tau+\xi'^2|)^{s}
|\tau+\xi'^2|^{\frac{4b+1+2\varepsilon}{2}}
\,
\frac{1}{|\xi_n|^{1+2\varepsilon}}.
$$
Hence, integrating $I$ with respect to  $\xi_n$ for $|\xi_n|>1$ we get
 the desired estimate \eqref{v1-mult-ine-1} in this subcase.

\vskip0.05in
\noindent
{\bf Subcase 2.} For this, we consider the following two subcases.

\vskip.05in
\noindent
$\bullet$
Subcase 2-1: $|\tau+\xi'^2|\le  \xi_n^2$
\quad
$\bullet$
Subcase 2-2:
 $|\tau+\xi'^2|> \xi_n^2$

\vskip0.05in
\noindent
\underline{Subcase 2-1.} Since $s< 0$, $b\ge 0$ and $1\le |\tau+\xi'^2|\le \xi_n^2$, we have 
$
(1+\xi'^2+\xi_n^2)^{s}
\lesssim
(1+\xi'^2+|\tau+\xi'^2|)^{s}
$
and
$
\big[
|\tau+\xi'^2|^{2b}
+
|\xi_n|^{4b}
\big]
\lesssim
|\xi_n|^{4b}.
$
Combining these inequalities with 
 \eqref{L2-est-mult},  we get 
\begin{align*}
%\label{integrand-bnd-2-1}
I(\xi',\xi_n,\tau)
\lesssim
(1+\xi'^2+|\tau+\xi'^2|)^{s}
%\cdot
|\xi_n|^{4b}
%\cdot
\xi_n^2
\,
\frac{|\tau+\xi'^2|^{m-1}}{|\xi_n|^{2m}}
\lesssim
(1+\xi'^2+|\tau+\xi'^2|)^{s}|\tau+\xi'^2|^{\frac{4b+1+2\varepsilon}{2}}
\,
\frac{1}{|\xi_n|^{1+2\varepsilon}},
\end{align*}
where in the last step we choose
 $m=2b+3/2+\varepsilon$.
Therefore, integrating $I$ with respect to $\xi_n$ over $|\xi_n|\ge 1$ we get
 the desired estimate \eqref{v1-mult-ine-1} in this subcase.

\vskip0.05in
\noindent
\underline{Subcase 2-2.} For this, we consider the following two subcases.

\vskip.05in
\noindent
$\bullet$
Subcase 2-2-1: $\xi'^2\ge |\tau+\xi'^2|$
\quad
$\bullet$
Subcase 2-2-2:
 $\xi'^2< |\tau+\xi'^2|$

\vskip0.05in
\noindent
{\it Subcase 2-2-1.}
Since $s< 0$, $b\ge0$, $|\tau+\xi'^2|> \xi_n^2$ and $\xi'^2\ge |\tau+\xi'^2|$, we have 
$
(1+\xi'^2+\xi_n^2)^{s}
\lesssim
(1+\xi'^2)^{s}
\lesssim
(1+\xi'^2+|\tau+\xi'^2|)^s
$
and
$
\big[
|\tau+\xi'^2|^{2b}
+
|\xi_n|^{4b}
\big]
\lesssim
|\tau+\xi'^2|^{2b},
$
which combined with 
 inequality \eqref{L2-est-mult}  gives 
\begin{align*}
%\label{integrand-bnd-2-2-1}
I(\xi',\xi_n,\tau)
\lesssim
(1+\xi'^2+|\tau+\xi'^2|)^{s}
%\cdot
|\tau+\xi'^2|^{2b}
%\cdot
\xi_n^2
\,
\frac{|\tau+\xi'^2|^{m-1}}{|\xi_n|^{2m}+|\tau+\xi'^2|^{m}}.
\end{align*}
Now, choosing $m=3/2+\varepsilon$ and using $|\tau+\xi'^2|\ge 1$, from the above estimate we get
$$
I(\xi',\xi_n,\tau)
\lesssim
(1+\xi'^2+|\tau+\xi'^2|)^{s}
|\tau+\xi'^2|^{2b}
\xi_n^2
\,
\frac{|\tau+\xi'^2|^{\frac12+\varepsilon}}{(|\xi_n|+1)^{3+2\varepsilon}}
\lesssim
(1+\xi'^2+|\tau+\xi'^2|)^{s}|\tau+\xi'^2|^{\frac{4b+1+2\varepsilon}{2}}
\,
\frac{1}{(|\xi_n|+1)^{1+2\varepsilon}}.
$$
Therefore, integrating $I$ with respect to  $\xi_n$ we get
 the desired estimate \eqref{v1-mult-ine-1} in this subcase.

\vskip0.05in
\noindent
{\it Subcase 2-2-2.}
Since $s< 0$, $b\ge0$ and $ |\tau+\xi'^2|> \xi_n^2$, we have 
$
(1+\xi'^2+\xi_n^2)^{s}
\lesssim
|\xi_n|^{2s}
$
and
$
\big[
|\tau+\xi'^2|^{2b}
+
|\xi_n|^{4b}
\big]
\lesssim
|\tau+\xi'^2|^{2b}.
$
Combining these inequalities  with estimate  \eqref{L2-est-mult}, we obtain
\begin{align*}
%\label{integrand-bnd-2-2-2}
I(\xi',\xi_n,\tau)
\lesssim
|\xi_n|^{2s}
%\cdot
|\tau+\xi'^2|^{2b}
%\cdot
\xi_n^2
\,
\frac{|\tau+\xi'^2|^{m-1}}{|\xi_n|^{2m}+|\tau+\xi'^2|^{m}}
\lesssim
|\tau+\xi'^2|^{\frac{2s+4b+1+2\varepsilon}{2}}
\,
\frac{1}{(|\xi_n|+1)^{1+2\varepsilon}},
\end{align*}
where in the last step we choose 
$m=s+3/2+\varepsilon$
 and we use $|\tau+\xi'^2|\ge 1$.
Furthermore, since $s< 0$, $\xi'^2< |\tau+\xi'^2|$ and $|\tau+\xi'^2|\ge 1$, from the above inequality we get 
$$
I(\xi',\xi_n,\tau)
\lesssim
(1+|\tau+\xi'^2|)^{s}
|\tau+\xi'^2|^{\frac{4b+1+2\varepsilon}{2}}
\,
\frac{1}{(|\xi_n|+1)^{1+2\varepsilon}}
\lesssim
(1+\xi'^2+|\tau+\xi'^2|)^{s}
|\tau+\xi'^2|^{\frac{4b+1+2\varepsilon}{2}}
\,
\frac{1}{(|\xi_n|+1)^{1+2\varepsilon}}.
$$
Therefore, integrating $I$ with respect to  $\xi_n$ we get
 the desired estimate \eqref{v1-mult-ine-1} in this subcase. 
This completes the proof of Lemma  \ref{v1-mult-lema}.
\,\,
$\square$

\vskip0.05in
\noindent
{\bf Temporal Bourgain spaces estimate \eqref{Y-B-est-redu-pure}.}
Using equation \eqref{sln-rpibvp-split}, we write $S[0,h;0]\simeq v_{r}+v_{i}$, where $v_r$, $v_i$ are defined in \eqref{vr-def}, \eqref{vi-def} respectively.
Estimate \eqref{Y-B-est-redu-pure} for $v_r$ is straightforward and follows a similar approach as estimate \eqref{B-est-redu-pure} for $v_r$. Hence, we will only prove  estimate \eqref{Y-B-est-redu-pure} for $v_i$ here.

To do this, using  equation \eqref{pure-sln-vi-split}, we rewrite $v_i=v_0+v_1$, where $v_0$ and $v_1$ are defined in \eqref{def-v0} and  \eqref{def-v1} respectively.

\vskip0.05in
\noindent
\underline{\it Estimate \eqref{Y-B-est-redu-pure} for $v_0$.}
To show this estimate, it suffices to prove the following inequality for the extension $V_0$, which is defined in \eqref{v0-extension}
\begin{equation*}
%\label{Y-V0-B-est}
\|\psi_4(t) V_0\|_{Y^{s,b}}^2
\doteq
\int_{\rr^{n+1}}
(1+|\tau|)^{s}
(1+|\tau+\xi^2|)^{2b}
|\widehat{\psi_4(t) V_0}(\xi,\tau)|^2
d\xi
d\tau
\lesssim
\|h\|_{\mathcal{B}^s}^2,
\quad
\,\,
s\in\rr
\,\,
\text{and}
\,\,
b\ge 0,
\end{equation*}
where $\psi_4=\psi(t/4)$ with $\psi$ being the time localizer given by \eqref{time-localizer}. 
Using inequality
$
(1+|\tau|)^{s}
\lesssim
(1+|\xi|)^{2s}(1+|\tau+\xi^2|)^{|s|},
$
we get
\begin{align}
\label{V0-Y-est}
\|\psi_4(t) V_0\|_{Y^{s,b}(\rr^{n+1})}^2
=
\int_{\rr^{n+1}}
(1+|\xi|)^{2s}
(1+|\tau+\xi^2|)^{|s|+2b}
|\widehat{\psi_4(t) V_0}(\xi,\tau)|^2
d\xi
d\tau
=
\|V_0\|_{X^{s,\frac12|s|+b}}^2.
\end{align}
Combining this with the space estimate \eqref{Bs-est-v0}
(where $b$ is replaced by $\frac12|s|+b$) we obtain 
the desired inequalities   \eqref{Y-B-est-redu-pure} and \eqref{Y-B-est-redu-pure-optimal} for $v_0$.

\vskip0.05in
\noindent
\underline{\it Estimate \eqref{Y-B-est-redu-pure} for $v_1$.}
Using the extension $V_1$, defined by \eqref{eqn-V1-1},
to prove estimate \eqref{Y-B-est-redu-pure} for $v_1$, it suffices to show that 
\begin{equation}
\label{Y-V1-B-est}
\|V_1\|_{Y^{s,b}(\rr^{n+1})}
\lesssim
\|h\|_{\mathcal{B}^s},
\quad
s\ge 0
\quad
\text{and}
\quad
0\le b<1/2.
\end{equation}
By the Fourier transform of $V_1$, given in \eqref{v1-xt-FT}, and by the definition of the $Y^{s,b}$-norm, we obtain
\begin{align}
\label{Y-norm-v1}
\|V_1\|_{Y^{s,b}}^2
\simeq
\int_{\tau+ \xi'^2\ge1}
(1+|\tau|)^{s}
\Big[
\int_{\rr}
(1+|\tau+\xi'^2+\xi_n^2|)^{2b}
|
\xi_n
F(\xi',\xi_n,\tau)
|^2
d\xi_n
\Big]
\Big|
\frac{
\widehat{h}^{x',t}
\big(
\xi',\tau
\big)
}
{(\tau+\xi'^2)^{\frac12}}
\Big|^2
d\tau
d\xi',
\end{align}
where the function $F$ is defined by \eqref{F-def-reduced-ibvp}.
Next, applying Lemma \ref{v1-mult-lema} with $s=0$, $b\ge 0$ and $\varepsilon>0$ we obtain
$
\int_{\rr}
(1+|\tau+\xi'^2+\xi_n^2|)^{2b}
|
\xi_n
F(\xi',\xi_n,\tau)
|^2
d\xi_n
\lesssim
|\tau+\xi'^2|^{\frac{4b+1+2\varepsilon}{2}},
$
which combined with  \eqref{Y-norm-v1} gives us
\begin{align}
\label{Y-norm-v1-no1}
\|V_1\|_{Y^{s,b}(\rr^{n+1})}^2
\lesssim
\int_{\rr^n}
\chi_{\tau+ \xi'^2\ge1}(\xi',\tau)
\cdot
(1+|\tau|)^{s}
%\cdot
|\tau+\xi'^2|^{\frac{4b-1+2\varepsilon}{2}}
%\cdot
\big|
\widehat{h}^{x',t}
(
\xi',\tau
)
\big|^2
d\tau
d\xi'.
\end{align}
Furthermore, using inequality $|\tau|\le |\tau+\xi'^2|+\xi'^2$,
and our assumption that $s\ge 0$,
 from  \eqref{Y-norm-v1-no1} we get
\begin{align}
\label{Y-norm-v1-no2}
\|V_1\|_{Y^{s,b}(\rr^{n+1})}^2
\lesssim
\int_{\tau+ \xi'^2\ge1}
(1+|\tau+\xi'^2|+\xi'^2)^{s}
%\cdot
|\tau+\xi'^2|^{\frac{4b-1+2\varepsilon}{2}}
%\cdot
\big|
\widehat{h}^{x',t}
(
\xi',\tau
)
\big|^2
d\tau
d\xi',
\,\,
s
\ge
0.
\end{align}
Finally,  choosing $\varepsilon=1-2b>0$ we
obtain inequality \eqref{Y-V1-B-est}. 
This also gives us the desired inequality \eqref{Y-B-est-redu-pure-optimal} in Remark \ref{optimal-homogeneous-est} for $v_1$.
We complete the proof of   Theorem \ref{pure-ibvp-thm}.
\,\,
$\square$
%

%
%
%
%%%%%%%%%%%%%%%%%%%%
%
%
% forced linear  ibvp estimates 
%
%
%%%%%%%%%%%%%%%%%%%% 
%
%
%
\section{Proof of forced linear  ibvp estimates (Theorem \ref{forced-linear-nls-thm})}
\label{sec:decompo}
\setcounter{equation}{0}
Here, we study the forced linear Schr\"odinger ibvp 
\eqref{lnls-ibvp}  using 
our knowledge of the reduced pure ibvp considered in 
the previous section. For this we decomposed it to three simpler 
problems.
%
%%%%%%%%%%%%%%%%%%%%%%%%%%%%%
%  
%
%    Decompositions
%
%
%
%%%%%%%%%%%%%%%%%%%%%%%%%%%%
%
 The first problem is the homogeneous initial value problem (ivp)
\begin{subequations}
\label{nls-ivp-homo}
\begin{align}
\label{nls-ivp-homo-eqn}
&iU _t+\Delta_{x}U
=
0,  
&& x\in \rr^n,\ t\in \rr,
\\
\label{nls-ivp-homo-ic}
&U (x,0)= U_0(x)\in H^s(\rr^n), 
\end{align}
\end{subequations}
with initial data  $U_0\in H^s(\rr^n)$ an extension 
 of the initial data $u_0\in H^s(\rr^{n-1}\times\rr^+)$ 
 of our forced linear ibvp \eqref{lnls-ibvp}, 
 such that
$
\|U_0\|_{H^s(\rr^n)} \leqslant 2\|u_0\|_{H^s(\rr^{n-1}\times\rr^+)}.
$
Its  solution is given by 
\begin{equation}
\label{nls-ivp-homo-sln}
U (x,t) 
= 
S\big[U_0; 0\big] (x, t) 
= 
\frac{1}{(2\pi)^n}
\int_{\rr^n}
e^{i\xi\cdot x}
e^{-i\xi^2t}
\widehat{U}_0(\xi)
d\xi,
\,\,
\text{where}
\,\,
\widehat{U}_0(\xi) 
=
\int_{\rr^n} e^{-i\xi\cdot x}\, U_0(x) dx,
\end{equation}
and satisfies  the  following estimates, which are proved in Section \ref{sec:homo-ivp}.
\begin{proposition}
[Homogeneous  ivp estimates]
\label{Livp-est}
We have
 the  space estimate in Bourgain spaces
\begin{align}
\label{ivp-homo-B-est}
\|\psi(t)\, S\big[U_0; 0\big]\|_{X^{s,b}(\rr^{n+1})}
\le
c_{s,b}
\|U_0\|_{H^{s}(\rr^n)},
&&
s\in\rr,\quad b\in\rr.
\end{align}
Also, 
we have the estimate in temporal Bourgain spaces
\begin{align}
\label{Y-ivp-homo-B-est}
\|\psi(t)
\, S\big[U_0; 0\big]\|_{Y^{s,b}(\rr^{n+1})}
\le
c_{s,b}
\|U_0\|_{H^{s}(\rr^n)},
&&
s\in\rr,\quad b\in\rr.
\end{align}
Furthermore, we have the
$\mathcal{B}^s$-norm estimate
\begin{align}
\label{ivp-homo-time-est-sharp}
\sup\limits_{x_n\in\rr}
\|
\psi(t)
S\big[U_0; 0\big](x_n)\|_{\mathcal{B}^s}
\le
c_s
\|\psi\|_{H^{|s|+2}}
\|U_0\|_{H^s(\rr^n)},
&&
s
\in
\rr.
\end{align}
\end{proposition}
The second problem is the following  inhomogeneous ivp with  initial data
zero 
\begin{subequations}
\label{nls-ivp-forced}
\begin{align}
\label{nls-ivp-forced-eqn}
&iW_t+  \Delta_{x}W
= 
w(x, t), 
&& x\in \rr^n,\ t\in \rr,
\\
\label{nls-ivp-forced-ic}
& W(x,0)
=
0,  
\end{align}
\end{subequations}
and forcing  $w$  an extension of the forcing $f$
appearing in our forced linear ibvp \eqref{lnls-ibvp}
 such that
\begin{equation}
\label{forcing-ext-bound}
\begin{cases}
\|w\|_{X^{s,-b}(\rr^{n+1})}
\le
2\|f\|_{X^{s,-b}_{\rr^{n-1}\times\rr^+\times(0,T)}},
&-\frac12\le s\le \frac12,
\\
\|w\|_{X^{s,-b}(\rr^{n+1})}+\|w\|_{Y^{s,-b}(\rr^{n+1})}
\hskip-0.03in
\le
\hskip-0.03in
2\|f\|_{X^{s,-b}_{\rr^{n-1}\times\rr^+\times(0,T)}}
\hskip-0.07in
+
2\|f\|_{Y^{s,-b}_{\rr^{n-1}\times\rr^+\times(0,T)}},
&s<-\frac12
\,\,
\text{or}
\,\,
s>\frac12.
\end{cases}
\end{equation}
The solution of this problem is given by the Duhamel representation
\begin{align}
\label{ivp-forcing-sln}
W(x,t) 
\doteq
S\big[0; w\big](x, t)
=&
\frac{-i}{(2\pi)^n}
\int_{\rr^n}
e^{i\xi\cdot x}
\Big[
\int_0^t
e^{-i\xi^2(t-t')}
\widehat{w}^{x}(\xi,t')
dt'
\Big]
d\xi
\\
\label{ivp-forcing-sln-d}
=&
-i
\int_{t'=0}^t  S\big[w(\cdot, t'); 0\big](x, t-t') dt',
\end{align}
where $\widehat{w}^{x}$ is the Fourier transform of $w$ with respect to $x$, 
and  $S\big[w(\cdot, t'); 0\big]$ in the Duhamel representation \eqref{ivp-forcing-sln-d} denotes the solution \eqref{nls-ivp-homo-sln} of homogeneous ivp \eqref{nls-ivp-homo} with $w(x, t')$ in place of the initial data $U_0(x)$.
We have the following estimates,
whose proof are provided in Section \ref{sec:forced-ivp}.
\begin{proposition}
[Inhomogeneous  ivp estimates]
\label{L-fivp-est}
The solution $S \big[0; w\big]$ 
satisfies 
the following
estimate
\begin{equation}
\label{forced-ivp-bourgain-est-nls}
\|\psi(t)S \big[0; w\big]\|_{X^{s,b}(\rr^{n+1}) }
\le
c_{s,b}
\|w\|_{X^{s,-b} (\rr^{n+1})},
\quad
s\in\rr,
\quad
0\le b<1/2.
\end{equation}
Also, it satisfies the following estimate in the temporal Bourgain spaces
\begin{align}
\label{Y-forced-ivp-bourgain-est-nls}
\|\psi(t)S \big[0; w\big]\|_{Y^{s,b}(\rr^{n+1}) }
\le
c_{s,b}
\big(\|w\|_{X^{s,-b} (\rr^{n+1})}
+
%c_{s,b}
\|w\|_{Y^{s,-b} (\rr^{n+1})}
\big),
\quad
s\in\rr,
\quad
0\le b<1/2.
\end{align}
Furthermore, 
we have  the $\mathcal{B}^s$-norm estimate
\begin{align}
\label{ivp-forced-est-tim-sharp}
\sup_{x_n\in\rr} \|\psi(t)S\big[0; w\big](x_n)\|_{\mathcal{B}^s} 
\le
\begin{cases}
c_{s,b}
\|w\|_{X^{s,-b}},
\quad
-\frac12
\le
 s  \le \frac12,
\,\,
0\le b<\frac12,
\\
c_{s,b}
\big(\|w\|_{X^{s,-b}}
+
%c_{s,b}
\|w\|_{Y^{s,-b}}
\big),
\quad
s\not\in[-\frac12,\frac12],
 \,\,
 0\le b<\frac12.
\end{cases}
\end{align}
\end{proposition}

Finally, the third problem in the decomposition of our 
forced linear ibvp   \eqref{lnls-ibvp} is the following  pure
ibvp having zero forcing, zero initial data, and boundary data
$g_0(x',t) - U (x',0,t)-W(x',0,t)$, that is 
\begin{subequations}
\label{nls-pure-ibvp}
\begin{align}
\label{nls-pure-ibvp-eqn}
&iu_t+\Delta_{x'}u+\p_{x_n}^2u= 0, 
\qquad
x'\in\rr^{n-1},\,x_n\in \rr^+, \ t\in (0,T), 
\\
\label{nls-pure-ibvp-ic}
&u(x',x_n,0)= 0,   
\\
\label{nls-pure-ibvp-bc}
&
u(x',0,t) = G_0(x',t) \doteq g_0(x',t) - U (x',0,t)
-W(x',0,t)
\in \mathcal{B}_T^s.
\end{align}
\end{subequations}
Note that, by our decomposition, the solution 
$
u(x',x_n, t) 
\doteq
S\big[0, G_0; 0\big](x',x_n, t)
$
of this problem, which is given  by  Fokas formula
\eqref{nd-lnls-utm-sln}, satisfies the relation
\begin{equation}
\label{UTM-deocomposition}
S\big[u_0, g_0; f\big]
=
S\big[U_0; 0\big]
+
S\big[0; w\big]
+
S\big[0, G_0;0\big],
\quad
x_n>0
\,\,
\text{and}
\,\,
0<t<T.
\end{equation}
Furthermore, it satisfies the following estimates.
\begin{proposition}
[Pure ibvp estimates]
\label{pure-ibvp-pro}
$S\big[0, G_0; 0\big](x, t)$ satisfies the Bourgain spaces estimate
\begin{align}
\label{pure-ibvp-B-est-sharp}
\|S\big[0, G_0; 0\big]\|_{X^{s,b}_{\rr^{n-1}\times\rr^+\times(0,T)}}
\le
c_{s,b}
\|G_0\|_{\mathcal{B}_T^s},
\quad
s>-\frac32, 
\,\,
\frac{2s-1}{4}\not\in \NN_0
\quad
\text{and}
\quad
0\le b<\frac12.
\end{align}
Also, we have the estimate in temporal Bourgain spaces
\begin{align}
\label{Y-pure-ibvp-B-est-sharp}
\|S\big[0, G_0; 0\big]\|_{Y^{s,b}_{\rr^{n-1}\times\rr^+\times(0,T)}}
\le
c_{s,b}
\|G_0\|_{\mathcal{B}_T^s},
\quad
s\ge 0, 
\,\,
\frac{2s-1}{4}\not\in \NN_0
\quad
\text{and}
\quad
0\le b<\frac12.
\end{align}
\end{proposition}

\noindent
{\bf Proof of Proposition \ref{pure-ibvp-pro}.} 
It is based on the estimates for the {\bf reduced} pure 
ibvp stated in Theorem \ref{pure-ibvp-thm} and proved
in the previous section.
However, for applying these estimates,  we need to {\bf extend} 
the  {\bf new} boundary data $G_0$, defined  by \eqref{nls-pure-ibvp-bc},
from $\rr^{n-1}\times (0,T)$ to $\rr^{n-1}\times \rr$ so that 
they become {\bf compactly supported in 
$\rr^{n-1}\times (0,2)$}. (Note that we have assumed  $T<1/2$.)
We do this next.
\begin{lemma}
\label{ext-thm-Bs-space}
Let $s>-\frac32$ and $\frac{2s-1}{4}\not\in \NN_0 $. 
If  $G_0\in \mathcal{B}^s_T$ 
and  satisfies the  condition
% at $t=0$:
%
\begin{align}
\label{zero-cond}
\p^j_{t}G_0(x',0)=0,\quad  j=0,1,2,\cdots,
\lfloor 
\frac{2s-1}{4}
\rfloor,
\end{align}
then we can extend $G_0(x',t)$ on $\rr^{n-1}\times (0, T)$
 to  $h(x',t)$ on $\rr^n$, such that the extension 
 $h$ satisfies
\begin{equation}
\label{Bs-ext-est}
\|h\|_{\mathcal{B}^s}
\lesssim
\|G_0\|_{\mathcal{B}_T^s}
\quad
\text{and}
\quad
\supp\,h\subset
\rr^{n-1}\times(0, 2).
\end{equation}
\end{lemma}

\vskip0.05in
To ensure that the boundary data $G_0$, defined by \eqref{nls-pure-ibvp-bc}, satisfies the condition \eqref{zero-cond}, we need to impose the following {\it compatibility conditions}
for $s>\frac 12$:
\begin{align}
\label{recursing-collec}
&g_0(x',0)
=
\phi_0
\doteq
u_0(x',0),
\\
&\p_t
g_0(x',0)
=
\phi_1
\doteq
i \Delta u_0 (x',0)\pm 
i|u_0(x',0)|^2u_0(x',0),
\\
\vdots
\nonumber
\\
\label{compa-condi}
&\p_t^k g_0(x',0)
=
\phi_k
\Big(u_0(x',0),\cdots,\Delta^{k}u_0(x', 0), \bar{u}_0(x', 0), \cdots,\Delta^{k}\bar{u}_0(x', 0)\Big),
\quad
k<\frac{2s-1}{4}.
\end{align}
Here, $\phi_k$, $k=0,1,2,\cdots$, is recursively determined by applying the operator $\partial_t$ to $\phi_{k-1}$ and then employing the nonlinear Schrödinger equation \eqref{NLS-eqn}.

The extension in $\mathcal{B}^s$ space is similar to the extension in Sobolev spaces $H^s$. In fact, for $s \geq 0$, we have
\begin{align*}
\|h\|_{\mathcal{B}^s}^2
\simeq
\|h\|_{X^{0,\frac{2s+1}{4}}}^2
+ 
\|h\|_{X^{s,\frac{1}{4}}}^2
\simeq
\int_{\rr}
%(1+\xi'^2)^0
\|
e^{i\xi'^2t}
\widehat{h}^{x'}
(\xi',t)
\|^2_{H_t^{\frac{2s+1}{4}}}
d\xi'
+
\int_{\rr}
(1+|\xi'|)^{2s}
\|
e^{i\xi'^2t}
\widehat{h}^{x'}
(\xi',t)
\|^2_{H_t^{\frac{1}{4}}}
d\xi'.
\end{align*}
Thus, the extension in the $\mathcal{B}^s$ space is analogous to the extension in the Sobolev spaces $H^s$, as described in Theorem 11.4 of \cite{lmbook}.
Also, for $-\frac{3}{2} < s < 0$, we can do the extension similarly.
Therefore, we omit the proof of Lemma \ref{ext-thm-Bs-space} here. 
Next, we use it to complete the proof  of Proposition \ref{pure-ibvp-pro}.

Using Lemma \ref{ext-thm-Bs-space}, we extend boundary data $G_0$ to $\rr^n$, such that 
the extension $h$ satisfies condition \eqref{Bs-ext-est}. Then, 
applying the reduced ibvp estimate  \eqref{B-est-redu-pure},
with boundary data $h$, we get 
\begin{align*}
%\label{ext-sln-bound-neg}
\|
S[0,G_0;0]
\|_{X^{s, b}_{\rr^{n-1}\times\rr^+\times (0,T)}} 
=
\|
S[0,h;0]
\|_{X^{s, b}_{\rr^{n-1}\times\rr^+\times (0,T)}} 
\le
\|
S[0,h;0]
\|_{X^{s, b}_{\rr^{n-1}\times\rr^+\times (0,2)}}
%\nonumber
\lesssim
\|h\|_{\mathcal{B}^s}
\lesssim
\|G_0\|_{\mathcal{B}_T^s}.
\end{align*}
This is the desired estimate \eqref{pure-ibvp-B-est-sharp}. Similarly, applying inequality \eqref{Y-B-est-redu-pure} with boundary data $h$ we get the desired estimate \eqref{Y-pure-ibvp-B-est-sharp}.
This completes the proof of Proposition 
\ref{pure-ibvp-pro}.
\,\,
$\square$

%%%%%%%%%%%%%%%%%%%%%%%%%%%%%%
%
%
%     Proof of linear estimates
%
%
%
%%%%%%%%%%%%%%%%%%%%%%%%%%%%
%
%
\vskip0.05in
Finally, combining the estimates for the three simpler problems above we obtain the  linear estimates for Schr\"odinger equation
(Theorem  \ref{forced-linear-nls-thm}).

\vskip0.05in
\nin
{\bf Proof of Theorem \ref{forced-linear-nls-thm}.} 
Here, we  prove only estimate \eqref{X-forced-linear-nls-est}.
The other estimates are similar.
Also, since the case $-\frac 12\le s< \frac 12$ is covered
by estimate \eqref{noY-forced-linear-nls-est},
we assume that $s>-\frac32$, $s\not\in [-\frac12,\frac12]$,  $\frac{2s-1}{4}\not\in \NN_0$.
Using the decomposition \eqref{UTM-deocomposition}, i.e. 
$$
S\big[u_0, g_0; f\big]
=
S\big[U_0; 0\big]
+
S\big[0; w\big]
+
S\big[0, G_0;0\big],
\quad
x_n>0
\,\,
\text{and}
\,\,
0<t<T,
$$
and applying the triangle inequality we have
\begin{align*}
\big\|
S\big[u_0, g_0; f\big] 
\big\|_{X^{s,b}_{\rr^{n-1}\times\rr^+\times(0,T)}}
\le&
\big\|
\psi(t)
S\big[U_0; 0\big]
\big\|_{X^{s,b}}
+
\big\|
\psi(t)
S\big[0; w\big]
\big\|_{X^{s,b}}
+
\big\|
S\big[0,  G_0; 0\big]
\big\|_{X^{s,b}_{\rr^{n-1}\times\rr^+\times(0,T)}}.
\end{align*}
Then, applying space estimates \eqref{ivp-homo-B-est},  \eqref{forced-ivp-bourgain-est-nls} and \eqref{pure-ibvp-B-est-sharp} we get 
\begin{align}
\label{sln-decom-est}
\big\|
S\big[u_0, g_0; f\big] 
\big\|_{X^{s,b}_{\rr^{n-1}\times\rr^+\times(0,T)}}
\lesssim
\big\|U_0\big\|_{H^s} 
+
\|w\|_{X^{s,-b}}
+\|G_0\|_{\mathcal{B}_T^s}.
\end{align}
Furthermore, using definition \eqref{nls-pure-ibvp-ic},  i.e. $ G_0(x',t) \doteq g_0(x',t) - U (x',0,t)-W(x',0,t)$, 
and $\mathcal{B}^s$-norm estimates \eqref{ivp-homo-time-est-sharp},  \eqref{ivp-forced-est-tim-sharp}  we obtain
\begin{align*}
\|G_0\|_{\mathcal{B}_T^s}
\le
\|g_0\|_{\mathcal{B}_T^s}
+
\|U (x',0,t)\|_{\mathcal{B}_T^s}
+
\|W (x',0,t)\|_{\mathcal{B}_T^s}
\lesssim
\|g_0\|_{\mathcal{B}_T^s}
+
\|U_0\|_{H^s}
+
\|w\|_{X^{s,-b}}+
\|w\|_{Y^{s,-b}}.
\end{align*}
Thus, from this inequality and estimate \eqref{sln-decom-est},
we get
\begin{align*}
\big\|
S\big[u_0, g_0; f\big] 
\big\|_{X^{s,b}_{\rr^{n-1}\times\rr^+\times(0,T)}}
\lesssim
\big\|U_0\big\|_{H^s}+\|g_0\|_{\mathcal{B}_T^s} 
+
\|w\|_{X^{s,-b}}+
\|w\|_{Y^{s,-b}}.
\end{align*}
Finally, using extension inequality for the initial data $u_0$,
$\|U_0\|_{H^s}\lesssim  \|u_0\|_{H^s(\rr^{n-1}\times\rr^+)}$
and the extension inequality \eqref{forcing-ext-bound} for the forcing $f$,
we get the desired  linear estimate 
\eqref{X-forced-linear-nls-est}.
\,\,
$\square$

%
%
%
%%%%%%%%%%%%%%%%%%%%
%
%
% Homogeneous linear ivp estimates
%
%
%%%%%%%%%%%%%%%%%%%% 
%
%
%
\section{Homogeneous ivp estimates (Proof of Proposition \ref{Livp-est}) }
\label{sec:homo-ivp}
\setcounter{equation}{0}
The proof of space estimate \eqref{ivp-homo-B-est} 
can be found in \cite{b1993-nls, tao-book}.
Concerning the  proof of temporal space estimate
 \eqref{Y-ivp-homo-B-est}, we reduce it to the space estimate \eqref{ivp-homo-B-est}.
In fact, using  inequality 
$
(1+|\tau|)^{s}
\lesssim
(1+|\xi|)^{2s}(1+|\tau+\xi^2|)^{|s|},
$
and applying estimate \eqref{ivp-homo-B-est}
 with $b$ replaced with $\frac12|s|+b$ we get 
$$
\|\psi(t) U\|_{Y^{s,b}}\lesssim \|\psi(t) U\|_{X^{s,\frac12|s|+b}}\lesssim \|U_0\|_{H^s}.
$$
Next, we focus on the proof  of the new
$\mathcal{B}^s$-norm estimate \eqref{ivp-homo-time-est-sharp}.

\vskip0.05in
\noindent
{\bf Proof of $\mathcal{B}^s$-norm estimate \eqref{ivp-homo-time-est-sharp}.}  
We recall  that the solution of the homogeneous ivp is given by
$$
U(x',x_n,t)
=
S\big[U_0; 0\big](x',x_n,t)
= 
\frac{1}{(2\pi)^n}
\int_{\rr^n}
e^{i\xi\cdot x}
e^{-i\xi^2t}
\widehat{U}_0(\xi)
d\xi.
$$
To calculate its $\mathcal{B}^s$-norm,
we shall make the change of variables $\xi_n^2\to \tau$ to find the Fourier transform of $U$ with respect to $x'$ and $t$.
For this, we need to split the $\xi_n$-integral
at $\xi_n=0$ and we write
$$
U(x',x_n,t)=I^++I^-,
$$
where $I^+$ and $I^-$ are defined as follows
\begin{align*}
 I^+(x',x_n,t)
\doteq
 \frac{1}{(2\pi)^n}
\int_{\rr^{n-1}}
\int_{\xi_n=0}^\infty
e^{i\xi'\cdot x'+i\xi_nx_n}
e^{-i(\xi'^2+\xi_n^2)t}
\widehat{U}_0(\xi',\xi_n)
d\xi_n
d\xi',
\\
I^-(x',x_n,t)
\doteq
 \frac{1}{(2\pi)^n}
\int_{\rr^{n-1}}
\int_{\xi_n=-\infty}^0
e^{i\xi'\cdot x'+i\xi_nx_n}
e^{-i(\xi'^2+\xi_n^2)t}
\widehat{U}_0(\xi',\xi_n)
d\xi_n
d\xi'.
\end{align*}
We shall prove $\mathcal{B}^s$-norm estimate \eqref{ivp-homo-time-est-sharp} for $I^+$. 
The estimate for $I^-$ is similar.
To deal with the singularity at $\xi_n=0$, we write
$I^+\simeq I_0+I_1$, where 
\begin{align}
\label{ivp-integral-0}
I_0
\doteq
\int_{\rr^{n-1}}
\int_{\xi_n=0}^1
e^{i\xi'\cdot x'+i\xi_nx_n}
e^{-i(\xi'^2+\xi_n^2)t}
\widehat{U}_0(\xi',\xi_n)
d\xi_n
d\xi',
\\
\label{ivp-integral-1}
I_1
\doteq
\int_{\rr^{n-1}}
\int_{\xi_n=1}^\infty
e^{i\xi'\cdot x'+i\xi_nx_n}
e^{-i(\xi'^2+\xi_n^2)t}
\widehat{U}_0(\xi',\xi_n)
d\xi_n
d\xi'.
\end{align}
\underline{\it Estimate \eqref{ivp-homo-time-est-sharp} for $I_0$.} 
This reads as follows
\begin{equation*}
%\label{psi-I0-est2}
\|\psi(t) I_0\|_{\mathcal{B}^s}
\lesssim
\|\psi\|_{H^{|s|+2}}
\|U_0\|_{H^s(\rr^n)},
\quad
s\in\rr,
\end{equation*}
where $\psi$ is the time localizer defined in \eqref{time-localizer}.
Since the Fourier transform  in  $(x', t)$  of $\psi(t) I_0$ is
\begin{align*}
\widehat{\psi(t)I}_0^{x',t}(\xi',x_n,\tau)
\simeq
\int_{\xi_n=0}^1
e^{i\xi_nx_n}
\widehat{\psi}(\tau+\xi'^2+\xi_n^2)
\widehat{U}_0(\xi',\xi_n)
d\xi_n,
\end{align*}
the  $\mathcal{B}^s$ norm of $\psi(t) I_0$, defined in \eqref{sharp-Bs}, is given by
\begin{align}
\label{psi-I0-est1}
\|\psi(t)I_0\|_{\mathcal{B}^s}^2
\doteq
\int_{\rr^n}
(1
+
\xi'^2
+
|\tau
+
\xi'^2|)^{s}
(1
+
|\tau
+
\xi'^2|)^{\frac12}
\big|
\widehat{\psi(t)I}_0^{x',t}(\xi',x_n,\tau)
\big|^2
d\tau
d\xi'
\simeq
\int_{\rr^{n-1}}
J_0(\xi')
d\xi',
\end{align}
where  $J_0(\xi')$ 
is defined as follows 
\begin{equation*}
%\label{homo-ivp-J0-def}
J_0(\xi')
\doteq
\int_{\rr}
\Big|
\int_{\xi_n=0}^1
e^{i\xi_nx_n}
(1+\xi'^2+|\tau+\xi'^2|)^{\frac{s}2}
(1+|\tau+\xi'^2|)^{\frac14}
\widehat{\psi}(\tau+\xi'^2+\xi_n^2)
\widehat{U}_0(\xi',\xi_n)
d\xi_n
\Big|^2
d\tau.
\end{equation*}
To estimate $J_0(\xi')$, we use the Minkowski inequality 
and we get 
\begin{align*}
%\label{homo-I0-minkow}
J_0(\xi')
\le
\Big(\int_{\xi_n=0}^1
\Big[
\int_{\rr}
(1+\xi'^2+|\tau+\xi'^2|)^{s}
(1+|\tau+\xi'^2|)^{\frac12}
\big|\widehat{\psi}(\tau+\xi'^2+\xi_n^2)
\widehat{U}_0(\xi',\xi_n)
\big|^2
d\tau
\Big]^{\frac12}
d\xi_n
\Big)^2.
\end{align*}
Furthermore, to separate   $|\tau+\xi'^2|$ and $\xi'^2$ 
we use the inequalities
\begin{align*}
%\label{homo-proof-mult-est1}
(1+\xi'^2+|\tau+\xi'^2|)^s
\lesssim&
(1+\xi'^2+
\xi_n^2)^s
(1+|\tau+\xi'^2+\xi_n^2|)^{|s|},
\quad
|\xi_n|
\le
1
\quad
\text{and}
\quad
s\in\rr,
\\
%\label{homo-proof-mult-est2}
(1+|\tau+\xi'^2|)^{1/2}
\le&
(1+\xi_n^2+|\tau+\xi'^2+\xi_n^2|)^{1/2}
\lesssim
(1+\xi_n^2)^{1/2}
(1+|\tau+\xi'^2+\xi_n^2|)^{1/2},
\end{align*}
and we bound $J_0$ as follows
\begin{align*}
%\label{homo-I0-minkow-1}
%\hskip-0.2in
J_0(\xi')
\lesssim
\hskip-0.03in
\Big(
\hskip-0.03in
\int_{0}^1
(1
\hskip-0.01in
+
\hskip-0.01in
\xi'^2
\hskip-0.01in
+
\hskip-0.01in
\xi_n^2)^{\frac{s}{2}}
(1
\hskip-0.01in
+
\hskip-0.01in
\xi_n^2)^{\frac14}
|\widehat{U}_0(\xi',\xi_n)|
\Big[
\hskip-0.03in
\int_\rr
\hskip-0.03in
(1
\hskip-0.01in
+
\hskip-0.01in
|\tau
\hskip-0.01in
+
\hskip-0.01in
\xi'^2
\hskip-0.01in
+
\hskip-0.01in
\xi_n^2|)^{|s|
\hskip-0.01in
+
\hskip-0.01in
\frac12}
\big|\widehat{\psi}(\tau
\hskip-0.01in
+
\hskip-0.01in
\xi'^2
\hskip-0.01in
+
\hskip-0.01in
\xi_n^2)
\big|^2
d\tau
\Big]^{\frac12}
d\xi_n
\Big)^2.
\end{align*}
In addition, utilizing the fact that
 $
 c_{\psi,s}
 \doteq
 \int_\rr
(1+|\tau+\xi'^2+\xi_n^2|)^{|s|+1/2}
\big|\widehat{\psi}(\tau+\xi'^2+\xi_n^2)
\big|^2
d\tau
\le
\|\psi\|_{H^{|s|+1}}^2
\lesssim 1$, 
from the last inequality we obtain the  bound
\begin{align}
\label{homo-I0-minkow-2}
J_0(\xi')
\lesssim
 c_{\psi,s}
\Big(\int_{\xi_n=0}^1
(1+\xi'^2+\xi_n^2)^{s/2}
(1+\xi_n^2)^{1/4}
|\widehat{U}_0(\xi',\xi_n)|
d\xi_n
\Big)^2.
\end{align}
Then, applying the Cauchy-Schwarz inequality for $\xi_n$ integral in  \eqref{homo-I0-minkow-2}, we get 
\begin{align}
\label{homo-I0-minkow-3}
J_0
\lesssim
\hskip-0.03in
c_{\psi,s}
\hskip-0.03in
\int_{0}^1
(1+\xi_n^2)^{1/2}
d\xi_n
\int_{0}^1
(1
\hskip-0.01in
+
\hskip-0.01in
\xi'^2
\hskip-0.01in
+
\hskip-0.01in
\xi_n^2)^{s}
|\widehat{U}_0(\xi',\xi_n)|^2
d\xi_n
\lesssim
\hskip-0.03in
c_{\psi,s}
\hskip-0.03in
\int_{0}^1
(1
\hskip-0.01in
+
\hskip-0.01in
\xi'^2
\hskip-0.01in
+
\hskip-0.01in
\xi_n^2)^{s}
|\widehat{U}_0(\xi',\xi_n)|^2
d\xi_n.
\end{align}
Finally, combining  \eqref{psi-I0-est1} with \eqref{homo-I0-minkow-3}, we get 
$
%\label{psi-I0-est2}
\|\psi I_0\|_{\mathcal{B}^s}
\hskip-0.01in
\lesssim
\hskip-0.01in
\|\psi\|_{H^{|s|+1}}
\hskip-0.01in
\|U_0\|_{H^s},
$
which gives  estimate \eqref{ivp-homo-time-est-sharp} for $I_0$.
\vskip0.1in
\noindent
\underline{\it Estimate \eqref{ivp-homo-time-est-sharp} for $I_1$.}  For $I_1$, defined by \eqref{ivp-integral-1}, to prove 
the desired estimate
\begin{equation}
\label{psi-I1-est2}
\|\psi(t) I_1\|_{\mathcal{B}^s}
\lesssim
\|\psi\|_{H^{|s|+2}}
\|U_0\|_{H^s(\rr^n)},
\quad
s\in\rr,
\end{equation}
first, we  make  the change of variables $\xi'=\xi'$ and $\tau=-(\xi'^2+\xi_n^2)$ or $\xi_n=\sqrt{-(\tau+\xi'^2)}$, with Jacobian 
$J=\partial (\xi', \tau)/\partial (\xi', \xi_n)=-2\xi_n$,
and we write it as follows
$$
I_1(x',x_n,t)
\simeq
\int_{\rr^{n-1}}
\int_{\tau=-\infty}^{-(\xi'^2+1)}
e^{i\xi'\cdot x'+i\sqrt{-(\tau+\xi'^2)} x_n}
e^{i\tau t}
\frac{1}{\sqrt{-(\tau+\xi'^2)}}
\widehat{U}_0
\big(\xi',\sqrt{-(\tau+\xi'^2)}
\big)
d\tau
d\xi'.
$$
Now we see that its Fourier transform with respect to $x'$ and $t$
is given by
\begin{align*}
\widehat{I_1}^{x',t}
(\xi',x_n,\tau)
=
\chi_{\tau<-(\xi'^2+1)}(\xi',\tau)
\cdot
e^{i\sqrt{-(\tau+\xi'^2)} x_n}
\frac{1}{\sqrt{-(\tau+\xi'^2)}}
\widehat{U}_0
\big(\xi',\sqrt{-(\tau+\xi'^2)}
\big),
\end{align*}
and using it in the definition of the $\mathcal{B}^s$-norm  we get
\begin{align*}
\|I_1\|_{\mathcal{B}^s}^2
\lesssim
\int_{\rr^n}
\chi_{\tau<-(\xi'^2+1)}(\xi',\tau)
(1+\xi'^2+|\tau+\xi'^2|)^{s}
(1+|\tau+\xi'^2|)^{\frac12}
\frac{1}{|\tau+\xi'^2|}
\big|\widehat{U}_0
(\xi',\sqrt{-(\tau+\xi'^2)})
\big|^2
d\tau
d\xi'.
\end{align*}
And, going back  to the variable $\xi_n$ via the change
 $\xi_n=\sqrt{-(\tau+\xi'^2)}$  we get 
\begin{align*}
\|I_1\|_{\mathcal{B}^s}^2
\lesssim
\int_{\rr^n}
\chi_{\xi_n> 1}(\xi_n)
\cdot
(1+\xi'^2+\xi_n^2)^{s}
\big|\widehat{U}_0
\big(\xi',\xi_n
\big)
\big|^2
d\xi_n
d\xi'
\lesssim
\|U_0\|_{H^s(\rr^n)}^2.
\end{align*}
The desired estimate \eqref{psi-I1-est2} for $I_1$ follows from the next multiplier lemma in 
boundary data space $\mathcal{B}^s$.

\begin{lemma}
\label{Bs-cut-lem}
Let $s\in \rr$ and $h\in \mathcal{B}^s$.
If $\psi(t)$ is a Schwartz function, then we have 
\begin{equation}
\label{Bs-cut-est}
\|\psi(t) h\|_{\mathcal{B}^s}
\lesssim
\|\psi\|_{H^{|s|+2}}
\|h\|_{\mathcal{B}^s},
\quad
s\in\rr.
\end{equation}
\end{lemma}
\noindent
The proof of inequality \eqref{Bs-cut-est} is similar 
to the corresponding multiplier lemma 
in  $X^{s,b}$ spaces, that is 
$
\|\psi(t) h\|_{X^{s,b}}
\lesssim
\|h\|_{X^{s,b}},
$
which can be found in \cite{tao-book} (see Lemma 2.11).
This completes the proof of Proposition \ref{Livp-est}.
\,\,
$\square$

%
%
%
%%%%%%%%%%%%%%%%%%%%
%
%
% Forced linear ivp estimates
%
%
%%%%%%%%%%%%%%%%%%%% 
%
%
%

%\newpage

\section{Inhomogeneous  ivp estimates
(Proof of Proposition \ref{L-fivp-est})
}
\label{sec:forced-ivp}
\setcounter{equation}{0}

In this section, we only prove the 
$\mathcal{B}^s$-norm estimate
 \eqref{ivp-forced-est-tim-sharp} for $S[0;w]$.
The proof of  estimate (\ref{forced-ivp-bourgain-est-nls}) 
in $X^{s,b}$ is similar to the corresponding estimate
in one-dimensional case
(see \cite{b1993-nls, b1993-kdv, tao-book}).
Additionally, the proof of the estimate \eqref{Y-forced-ivp-bourgain-est-nls} in $Y^{s,b}$ is analogous to the proof of  estimate \eqref{forced-ivp-bourgain-est-nls}.

\vskip.05in
\noindent
{\bf Proof of $\mathcal{B}^s$-norm estimate \eqref{ivp-forced-est-tim-sharp}.}
We start with recalling the Duhamel solution formula \eqref{ivp-forcing-sln}
$$
S\big[0; w\big](x, t)
=
\frac{-i}{(2\pi)^n}
\int_{\rr^n}
e^{i\xi\cdot x}
\Big[
\int_0^t
e^{-i\xi^2(t-t')}
\widehat{w}^{x}(\xi,t')
dt'
\Big]
d\xi.
$$
Using $\widehat{w}^{x}(\xi, t')
=
\frac{1}{2\pi}
\int_{\rr}
e^{i\tau t'}
\widehat{w} (\xi,\tau)d\tau$ and performing the $t'$ integration
$
 \int_0^t
e^{i(\tau+\xi^2)t'}   dt' 
=
-i \frac{e^{i(\tau+\xi^2)t} - 1}{\tau+\xi^2},
$
we  rewrite the solution formula  $\psi(t)S\big[0; w\big](x,t)$  as follows
\begin{subequations}
\label{f-nls-decomp}
\begin{align}
\label{w2-term-nls}
\psi(t)
S\big[0; w\big]
\simeq&
\psi(t)
\int_{\rr^{n+1}}
e^{i(\xi\cdot x+\tau t)}
\frac{1-\psi(\tau+\xi^2)}{\tau+\xi^2}
\widehat{w}(\xi,\tau)  d\tau  d\xi 
\\
\label{w3-term-nls}
-&
\psi(t)
\int_{\rr^{n+1}}
e^{i(\xi\cdot x-\xi^2t)}
\frac{1-\psi(\tau+\xi^2)}{\tau+\xi^2}
\widehat{w}(\xi,\tau)   d\tau  d\xi 
\\
\label{w4-term-nls}
+&
\psi(t)
\int_{\rr^{n+1}}
e^{i(\xi\cdot x-\xi^2t)}
 \frac{\psi(\tau+\xi^2)[e^{i(\tau+\xi^2)t} - 1]}{\tau+\xi^2}
\widehat{w}(\xi,\tau)   d\tau  d\xi. 
\end{align}
\end{subequations}
Next, we estimate the $\mathcal{B}^s$-norm of  terms
 \eqref{w2-term-nls}--\eqref{w4-term-nls} separately.
 We begin with the estimate for terms \eqref{w3-term-nls} and \eqref{w4-term-nls} that are based on the estimate for 
 the homogeneous ivp \eqref{nls-ivp-homo}.
\vskip.05in
\nin
\underline{\it Estimate \eqref{ivp-forced-est-tim-sharp} for  term \eqref{w3-term-nls}.} 
Rewriting  this term as the solution of homogeneous ivp \eqref{nls-ivp-homo} with initial data $U_0=F$ we get
$$
\eqref{w3-term-nls}
\simeq
\psi(t) S[F;0],
\quad
\text{where }
\quad
\widehat{F}(\xi)
\simeq
\int_{\rr}
\frac{1-\psi(\tau+\xi^2)}{\tau+\xi^2}
\widehat{w}(\xi,\tau) d\tau,
$$
and applying the $\mathcal{B}^s$-norm estimate \eqref{ivp-homo-time-est-sharp} for the homogeneous ivp  we have
\begin{align}
\label{w3-term-est}
\sup\limits_{x_n\in\rr}\|
\eqref{w3-term-nls}
\|_{\mathcal{B}^s}^2
\lesssim
\sup\limits_{x_n\in\rr}\|
\psi(t)
S[F;0]
\|_{\mathcal{B}^s}^2x
\lesssim
\|F\|_{H^{s}(\rr^n)}^2,
\,\,  s\in\rr.
\end{align}
Now, since $\psi(\tau+\xi^2)\le 1$ and 
$\psi(\tau+\xi^2)=1$ when $|\tau+\xi^2|\le \frac12$, 
we estimate  $\|F\|_{H^{s}(\rr^n)}$ as follows
\begin{align}
\label{inhom-termb}
\hskip-0.1in
\|F\|_{H^{s}(\rr^n)}^2
\simeq&
\int_{\rr^n}
(1+|\xi|)^{2s}
\Big|
\int_{\rr}
\frac{1-\psi(\tau+\xi^2)}{\tau+\xi^2}
\widehat{w}(\xi,\tau)   
d\tau
\Big|^2
d\xi
\lesssim
\int_{\rr^n}
(1+|\xi|)^{2s}
\Big(
\int_\rr
\frac{|\widehat{w}(\xi,\tau)| \, d\tau}{1+|\tau+\xi^2|}
\Big)^2
d\xi.
\end{align}
Furthermore, applying the Cauchy-Schwarz inequality for the $\tau$-integration with $0\le b<\frac12$, we have
\begin{align}
\label{2nd-bilinear-est-nls}
\Big(
\int_\rr
\frac{|\widehat{w}(\xi,\tau)| }{1+|\tau+\xi^2|}
d\tau
\Big)^2
\le
\int_\rr
\frac{1 }{(1+|\tau+\xi^2|)^{2(1-b)}}
d\tau
\int_\rr
\frac{|\widehat{w}(\xi,\tau)|^2 }{(1+|\tau+\xi^2|)^{2b}}
d\tau
\lesssim
\int_\rr
\frac{|\widehat{w}(\xi,\tau)|^2 }{(1+|\tau+\xi^2|)^{2b}}
d\tau.
\end{align}
Finally, combing  inequalities \eqref{w3-term-est}, \eqref{inhom-termb}  with \eqref{2nd-bilinear-est-nls}, we get 
$
\sup\limits_{x_n\in\rr}\|
\eqref{w3-term-nls}
\|_{\mathcal{B}^s}
\lesssim
\|w\|_{X^{s,-b}},
$
which gives the desired  estimate  \eqref{ivp-forced-est-tim-sharp} for \eqref{w3-term-nls}.

\vskip0.05in
\nin
\underline{\it Estimate \eqref{ivp-forced-est-tim-sharp} for  term \eqref{w4-term-nls}.} 
Using Taylor's expansion, we rewrite this term as follows
$$
\eqref{w4-term-nls}
\simeq
\sum\limits_{k=1}^\infty
\frac{1}{k!}t^k
\psi(t)S[c_k,0],
\,\,
\text{ where  }
\,\,
\widehat{c_k}(\xi)
\simeq
\int_{\rr}\psi(\tau+\xi^2)(\tau+\xi^2)^{k-1}
\widehat{w}(\xi,\tau)d\tau.
$$
Letting 
$\psi_k(t)\doteq t^k\psi(t)$ and 
using estimate \eqref{ivp-homo-time-est-sharp},  we get
\begin{align*}
\sup\limits_{x_n\in\rr}\|
\eqref{w4-term-nls}
\|_{\mathcal{B}^s}
\lesssim&
\sum_{k=1}^{\infty}\frac{1}{k!}
\sup\limits_{x_n\in\rr}
\|
\psi_k(t)
S[c_k;0]
\|_{\mathcal{B}^s}
\\
\lesssim&
\sum_{k=1}^{\infty}\frac{c_{\psi_k}}{k!}\Big( \int_{\rr^n} (1+|\xi|)^{2s} 
\Big|\int_{\rr}\psi(\tau+\xi^2)\cdot(\tau+\xi^2)^{k-1}
\widehat{w}(\xi,\tau)d\tau
\Big|^2 d\xi \Big)^{1/2},
\quad
s\in\rr,
\end{align*}
where $c_{\psi_k}\doteq \|\psi_k\|_{H^{|s|+2}}$  is given in estimate \eqref{ivp-homo-time-est-sharp}.
Also, since $\supp\, \psi \subset (-1,1)$ and  $|\psi(\tau+\xi^2)|\leq 1$
from the above inequality we obtain 
\begin{align*}
\sup\limits_{x_n\in\rr}\|
\eqref{w4-term-nls}
\|_{\mathcal{B}^s}
\lesssim&
\sum_{k=1}^{\infty}\frac{c_{\psi_k}}{k!}
\Big(\int_{\rr^n} (1+|\xi|)^{2s} \Big(\int_{|\tau+\xi^2|\leq 1} 
|\widehat{w}(\xi,\tau)|d\tau\Big)^2 d\xi \Big)^{1/2}
\\
\lesssim&
\Big(\int_{\rr^n} (1+|\xi|)^{2s} \Big(\int_{\rr} 
\frac{|\widehat{w}(\xi,\tau)|}{1+|\tau+\xi^2|}d\tau\Big)^2 d\xi
\Big)^{1/2}
\lesssim
\|w\|_{X^{s,-b}(\rr^{n+1})},
\end{align*}
where in the second step by using  
$c_{\psi_k}\lesssim k^m$, with $m=\lfloor |s|\rfloor+3$, we get 
$
\sum_{k=1}^{\infty}\frac{c_{\psi_k}}{k!}
\lesssim 1.
$
And the last relation  follows from  \eqref{2nd-bilinear-est-nls}.
This gives the desired  estimate \eqref{ivp-forced-est-tim-sharp} for term \eqref{w4-term-nls}. 

\vskip0.05in\nin
\underline{\it Estimate \eqref{ivp-forced-est-tim-sharp}  for  term \eqref{w2-term-nls}.}
We begin by writing 
$
\eqref{w2-term-nls}
\simeq
\psi(t)
h(x',x_n,t),
$
where
\begin{equation*}
%\label{w2-term-express}
h(x',x_n,t)
\hskip-0.03in
\doteq
\hskip-0.03in
\int_{\rr^{n+1}}
\hskip-0.03in
e^{i(\xi'\cdot x'+\xi_nx_n+\tau t)}
\frac{1-\psi(\tau+\xi'^2+\xi_n^2)}{\tau+\xi'^2+\xi_n^2}
\widehat{w}(\xi',\xi_n,\tau)  d\tau  d\xi'd\xi_n,
\end{equation*}
and using Lemma \ref{Bs-cut-lem}, we get 
$
\|\psi(t) h(x_n)\|_{\mathcal{B}^s}\lesssim \|h(x_n)\|_{\mathcal{B}^s},
$
from which we see that  for proving  estimate \eqref{ivp-forced-est-tim-sharp} for \eqref{w2-term-nls}, it suffices to show that 
\begin{align}
\label{sharp-Bs-h-est}
\sup\limits_{x_n\in\rr}
\|h(x_n)\|_{\mathcal{B}^s}
\lesssim
\begin{cases}
\|w\|_{X^{s,-b}(\rr^{n+1})},
&\quad
-\frac12\le s\le \frac12
\,\,
\text{and}
\,\,
 0\le b<\frac12,
\\
\|w\|_{X^{s,-b}(\rr^{n+1})}
+
\|w\|_{Y^{s,-b}(\rr^{n+1})},
&\quad
s<-\frac12
\,\,
\text{or}
\,\,
s>\frac12
\,\,
\text{and}
\,\,
 0\le b<\frac12.
 \end{cases}
\end{align}
Also, we  decompose $h$ as   $h=h^++h^-$, where 
\begin{align*}
%\label{inhom-h-pos-def}
h^+
\doteq
\int_{\rr^n}
e^{i(\xi'\cdot x'+\tau t)}
\Big(
\int_{\xi_n=0}^\infty
e^{i\xi_n x_n}
\frac{1-\psi(\tau+\xi'^2+\xi_n^2)}{\tau+\xi'^2+\xi_n^2}
\widehat{w}(\xi',\xi_n,\tau)d\xi_n
\Big)
d\tau  d\xi',
\\
%\label{inhom-h-neg-def}
h^-
\doteq
\int_{\rr^n}
e^{i(\xi'\cdot x'+\tau t)}
\Big(
\int_{\xi_n=-\infty}^0
e^{i\xi_n x_n}
\frac{1-\psi(\tau+\xi'^2+\xi_n^2)}{\tau+\xi'^2+\xi_n^2}
\widehat{w}(\xi',\xi_n,\tau)d\xi_n
\Big)
d\tau  d\xi',
\end{align*}
and prove estimate \eqref{sharp-Bs-h-est} for $h^+$,
since the proof of the estimate  for $h^-$ similar.
Furthermore, we  write $h^+=h_0+h_1$, where
\begin{align}
\label{inhom-h-0-def}
h_0
\doteq
\int_{\rr^n}
e^{i(\xi'\cdot x'+\tau t)}
\Big(
\int_{\xi_n=0}^1
e^{i\xi_n x_n}
\frac{1-\psi(\tau+\xi'^2+\xi_n^2)}{\tau+\xi'^2+\xi_n^2}
\widehat{w}(\xi',\xi_n,\tau)d\xi_n
\Big)
d\tau  d\xi',
\\
\label{inhom-h-1-def}
h_1
\doteq
\int_{\rr^n}
e^{i(\xi'\cdot x'+\tau t)}
\Big(
\int_{\xi_n=1}^\infty
e^{i\xi_n x_n}
\frac{1-\psi(\tau+\xi'^2+\xi_n^2)}{\tau+\xi'^2+\xi_n^2}
\widehat{w}(\xi',\xi_n,\tau)d\xi_n
\Big)
d\tau  d\xi',
\end{align}
and prove estimate \eqref{sharp-Bs-h-est} for both 
$h_0$ and $h_1$.

\vskip0.05in
\nin
{\it \bf Estimate \eqref{sharp-Bs-h-est} for $h_0$.} Taking the Fourier transform with respect to $x'$ and $t$, we get 
\begin{equation*}
%\label{inhom-h-0-FT}
\widehat{h}_0^{x',t}
(\xi',x_n,\tau)
\simeq
\int_{\xi_n=0}^1
e^{i\xi_n x_n}
\frac{1-\psi(\tau+\xi'^2+\xi_n^2)}{\tau+\xi'^2+\xi_n^2}
\widehat{w}(\xi',\xi_n,\tau)d\xi_n,
\end{equation*}
 and using the definition of $\mathcal{B}^s$ norm, defined in \eqref{sharp-Bs}, we have
\begin{align*}
%\label{inhom-h-0-Bs-est1}
\|h_0\|_{\mathcal{B}^s}^2
\simeq
\int_{\rr^n}
(1+\xi'^2+|\tau+\xi'^2|)^{s}
(1+|\tau+\xi'^2|)^{\frac12}
\Big|
\int_{0}^1
e^{i\xi_n x_n}
\frac{1-\psi(\tau+\xi'^2+\xi_n^2)}{\tau+\xi'^2+\xi_n^2}
\widehat{w}(\xi',\xi_n,\tau)d\xi_n
\Big|^2
d\xi'
d\tau.
\end{align*}
Now, applying  the Cauchy-Schwarz inequality in $\xi_n$-integration we get
\begin{align*}
\|h_0(x_n)\|_{\mathcal{B}^s}^2
\lesssim&
\int_{\rr^n}
\int_{\xi_n=0}^1
(1+\xi'^2+|\tau+\xi'^2|)^{s}
(1+|\tau+\xi'^2|)^{1/2}
\frac{|\widehat{w}(\xi',\xi_n,\tau)|^2}{(1+|\tau+\xi'^2+\xi_n^2|)^2}
d\xi_n
d\xi'
d\tau.
\end{align*}
Furthermore, employing the following two relations  
$$
1+\xi'^2+|\tau+\xi'^2|\simeq 1+\xi'^2+\xi_n^2+|\tau+\xi'^2+\xi_n^2|
\,\,
\text{and}
\,\,
1+|\tau+\xi'^2|\simeq 1+|\tau+\xi'^2+\xi_n^2|,
\quad
|\xi_n|
\le 1,
$$
 we obtain
\begin{align}
\label{inhom-h-0-Bs-est2}
\|h_0(x_n)\|_{\mathcal{B}^s}^2
\lesssim&
\int_{\rr^n}
\int_{\xi_n=0}^1
(1+\xi'^2+\xi_n^2+|\tau+\xi'^2+\xi_n^2|)^{s}
\frac{|\widehat{w}(\xi',\xi_n,\tau)|^2}{(1+|\tau+\xi'^2+\xi_n^2|)^{3/2}}
d\xi_n
d\xi'
d\tau.
\end{align}
At this point,  we consider the following two cases 
to separate $|\tau+\xi'^2+\xi_n^2|$ and $\xi'^2+\xi_n^2$.

\vskip0.05in
\noindent
$\bullet$ Case 1: $-\frac12\le s\le \frac12$
\qquad
$\bullet$ Case 2: $s<-\frac12$ or $s>\frac12$

\vskip0.05in
\noindent
{\it Case 1.} 
Using
$(1+\xi'^2+\xi_n^2+|\tau+\xi'^2+\xi_n^2|)^{s}
\lesssim
(1+\xi'^2+\xi_n^2)^s(1+|\tau+\xi'^2+\xi_n^2|)^{|s|}$, 
from \eqref{inhom-h-0-Bs-est2} we obtain
\begin{align*}
%\label{inhom-h-0-Bs-est3}
\|h_0(x_n)\|_{\mathcal{B}^s}^2
\lesssim&
\int_{\rr^n}
\int_{\xi_n=0}^1
(1+\xi'^2+\xi_n^2)^{s}
\frac{|\widehat{w}(\xi',\xi_n,\tau)|^2}{(1+|\tau+\xi'^2+\xi_n^2|)^{3/2-|s|}}
d\xi_n
d\xi'
d\tau
\lesssim
\|w\|_{X^{s,-\frac34+\frac12|s|}}^2,
\end{align*}
and since our assumption $|s|\le \frac12$
imply that
$
-\frac34+\frac12|s|\le -b,
$
when  $0\le b<\frac12$,
we get the
desired estimate \eqref{sharp-Bs-h-est} for $h_0$ in Case 1.

\vskip0.05in
\noindent
{\it Case 2.} Utilizing inequality $(1+\xi'^2+\xi_n^2+|\tau+\xi'^2+\xi_n^2|)^{s}
\lesssim
(1+\xi'^2+\xi_n^2)^s+(1+|\tau|)^{s}$, from
 \eqref{inhom-h-0-Bs-est2} we get 
\begin{align*}
%\label{inhom-h-0-Bs-est4}
\hskip-0.1in
\|h_0(x_n)\|_{\mathcal{B}^s}^2
\hskip-0.03in
\lesssim
\hskip-0.03in
\int_{\rr^n}
\hskip-0.03in
\int_{\xi_n=0}^1
\hskip-0.03in
\Big[
(1
\hskip-0.01in
+
\hskip-0.01in
\xi'^2
\hskip-0.01in
+
\hskip-0.01in
\xi_n^2)^s
\hskip-0.01in
+
\hskip-0.01in
(1
\hskip-0.01in
+
\hskip-0.01in
|\tau|)^{s}
\Big]
\frac{|\widehat{w}(\xi',\xi_n,\tau)|^2
d\xi_n
d\xi'
d\tau
}{(1+|\tau+\xi'^2+\xi_n^2|)^{3/2}}
\lesssim
\|w\|_{X^{s,-\frac34}}^2
\hskip-0.01in
+
\hskip-0.01in
\|w\|_{Y^{s,-\frac34}}^2.
\end{align*}
Then, using the assumption $0\le b<\frac12$ we get the
desired estimate \eqref{sharp-Bs-h-est} for $h_0$ in Case 2.

\vskip0.05in
\noindent
{\it \bf Estimate \eqref{sharp-Bs-h-est} for $h_1$.}
Since the Fourier transform in $x'$ and $t$ of  
$h_1$, which  defined in \eqref{inhom-h-1-def}, is given by
\begin{equation*}
%\label{inhom-h-1-FT}
\widehat{h}_1^{x',t}
(\xi',x_n,\tau)
\simeq
\int_{\xi_n=1}^\infty
e^{i\xi_n x_n}
\frac{1-\psi(\tau+\xi'^2+\xi_n^2)}{\tau+\xi'^2+\xi_n^2}
\widehat{w}(\xi',\xi_n,\tau)d\xi_n,
\end{equation*}
its  $\mathcal{B}^s$ norm \eqref{sharp-Bs} is given by
\begin{align*}
%\label{inhom-h-1-Bs-est1}
\hskip-0.1in
\|h_1\|_{\mathcal{B}^s}^2
\hskip-0.03in
\simeq
\hskip-0.03in
\int_{\rr^n}
\hskip-0.03in
(1
\hskip-0.01in
+
\hskip-0.01in
\xi'^2
\hskip-0.01in
+
\hskip-0.01in
|\tau
\hskip-0.01in
+
\hskip-0.01in
\xi'^2|)^{s}
(1
\hskip-0.01in
+
\hskip-0.01in
|\tau
\hskip-0.01in
+
\hskip-0.01in
\xi'^2|)^{\frac12}
\Big|
\int_{1}^\infty
\hskip-0.03in
e^{i\xi_n x_n}
\frac{1
\hskip-0.01in
-
\hskip-0.01in
\psi(\tau
\hskip-0.01in
+
\hskip-0.01in
\xi'^2
\hskip-0.01in
+
\hskip-0.01in
\xi_n^2)}{\tau+\xi'^2+\xi_n^2}
\widehat{w}(\xi',\xi_n,\tau)d\xi_n
\Big|^2
d\xi'
d\tau.
\end{align*}
And, making the change of variables $\eta_n=\xi_n^2$ or $\xi_n=\sqrt{\eta_n}$, we get 
$$
\int_{\xi_n=1}^\infty
e^{i\xi_n x_n}
\frac{1-\psi(\tau+\xi'^2+\xi_n^2)}{\tau+\xi'^2+\xi_n^2}
\widehat{w}(\xi',\xi_n,\tau)d\xi_n
=
\int_{\eta_n=1}^\infty
e^{i\sqrt{\eta_n} x_n}
\frac{1-\psi(\tau+\xi'^2+\eta_n)}{\tau+\xi'^2+\eta_n}
\widehat{w}(\xi',\sqrt{\eta_n},\tau)
\frac{1}{\sqrt{\eta_n}}
d\eta_n.
$$
Also, moving the absolute value inside the $\eta_n$-integral and 
taking into consideration that $\psi(x)\equiv 1$ when $|x|\le \frac12$, we obtain
\begin{align}
\label{inhom-h-1-Bs-est2}
\|h_1(x_n)\|_{\mathcal{B}^s}^2
\lesssim
\int_{\rr^n}
(1+\xi'^2+|\tau+\xi'^2|)^{s}
(1+|\tau+\xi'^2|)^{\frac12}
\Big(
\int_{\eta_n=1}^\infty
\frac{|\widehat{w}(\xi',\sqrt{\eta_n},\tau)|
d\eta_n
}{(1+|\tau+\xi'^2+\eta_n|)\cdot \sqrt{\eta_n}}
\Big)^2
d\xi'
d\tau.
\end{align}
Moreover, by comparing $\sqrt{\eta_n}$ with $\frac12|\tau+\xi'^2|$,
we rewrite  $\widehat{w}=\widehat{w}_1+\widehat{w}_2$, where 
\begin{align*}
%\label{inhom-w12-def}
\widehat{w}_1(\xi',\sqrt{\eta_n},\tau)
\doteq
\chi_{1<\eta_n\le \frac12|\tau+\xi'^2|}
(\xi',\eta_n, \tau)
\cdot
\widehat{w}
\quad
\text{and}
\quad
\widehat{w}_2(\xi',\sqrt{\eta_n},\tau)
\doteq
\chi_{\eta_n> \frac12|\tau+\xi'^2|}
(\xi',\eta_n, \tau)
\cdot
\widehat{w}.
\end{align*}
Then, since $|\widehat{w}|\le |\widehat{w}_1|+|\widehat{w}_2|$, 
from inequality \eqref{inhom-h-1-Bs-est2} we obtain $\|h_1(x_n)\|_{\mathcal{B}^s}^2\lesssim I_1+I_2$, where 
\begin{align}
\label{inhom-h-1-Bs-I12-def}
I_j
\doteq
\int_{\rr^n}
(1+\xi'^2+|\tau+\xi'^2|)^{s}
(1+|\tau+\xi'^2|)^{1/2}
\Big(
\int_{1}^\infty
\frac{|\widehat{w}_j(\xi',\sqrt{\eta_n},\tau)|\,\, d\eta_n}{(1+|\tau+\xi'^2+\eta_n|)\cdot \sqrt{\eta_n}}
\Big)^2
d\xi'
d\tau,
\,\,
j=1,2.
\end{align}
Therefore, to show estimate \eqref{sharp-Bs-h-est} for $h_1$, it suffices to show that
\begin{equation}
\label{inhom-h1-I-est}
I_j
\lesssim
\begin{cases}
\|w\|_{X^{s,-b}(\rr^{n+1})},
&\quad
-\frac12\le s\le \frac12
\,\,
\text{and}
\,\,
 0\le b<\frac12,
\\
\|w\|_{X^{s,-b}(\rr^{n+1})}
+
\|w\|_{Y^{s,-b}(\rr^{n+1})},
&\quad
s<-\frac12
\,\,
\text{or}
\,\,
s>\frac12
\,\,
\text{and}
\,\,
 0\le b<\frac12,
\end{cases}
\quad
j=1,2.
\end{equation}
{\it Estimate \eqref{inhom-h1-I-est} for $I_1$.} We begin by estimating the $\eta_n$-integral in \eqref{inhom-h-1-Bs-I12-def}. Using $|\tau+\xi'^2+\eta_n|\gtrsim \eta_n$, which follows from
 $
 1<\eta_n\le \frac12|\tau+\xi'^2|,
 $  and applying the
Cauchy-Schwarz inequality,
 for $\varepsilon>0$, we get 
\begin{align*}
\Big(
\int_{\eta_n=1}^\infty
\frac{|\widehat{w}_1(\xi',\sqrt{\eta_n},\tau)|}{(1+|\tau+\xi'^2+\eta_n|)\cdot \sqrt{\eta_n}}
d\eta_n
\Big)^2
\lesssim
\frac{1}{\varepsilon}
\int_1^\infty
\chi_{\eta_n\le \frac12|\tau+\xi'^2|}
(\xi',\eta_n, \tau)
\cdot
\frac{|\widehat{w}(\xi',\sqrt{\eta_n},\tau)|^2}{(1+|\tau+\xi'^2+\eta_n|)^{2-\varepsilon}}
d\eta_n.
\end{align*}
since 
$
\int_1^\infty \eta_n^{-(1+\varepsilon)} d\eta_n =1/\varepsilon.
$
This combined with definition \eqref{inhom-h-1-Bs-I12-def} gives 
\begin{align*}
%\label{inhom-h-1-I1-est}
I_1
\lesssim
\frac{1}{\varepsilon}
\int_{\rr^n}
(1+\xi'^2+|\tau+\xi'^2|)^{s}
(1+|\tau+\xi'^2|)^{\frac12}
\int_1^\infty
\chi_{\eta_n\le \frac12|\tau+\xi'^2|}
(\xi',\eta_n, \tau)
\frac{|\widehat{w}(\xi',\sqrt{\eta_n},\tau)|^2}{(1+|\tau+\xi'^2+\eta_n|)^{2-\varepsilon}}
d\eta_n
d\xi'
d\tau.
\end{align*}
Furthermore, making the change of variables $\xi_n=\sqrt{\eta_n}$, we get 
\begin{align}
\label{inhom-h-1-I1-est-1}
\hskip-0.1in
I_1
\lesssim
\hskip-0.05in
\frac{1}{\varepsilon}
\int_{\rr^n}
\int_1^\infty
\chi_{\xi_n^2\le \frac12|\tau+\xi'^2|}
(\xi',\xi_n, \tau)
\frac{
(1+\xi'^2+|\tau+\xi'^2|)^{s}
(1+|\tau+\xi'^2|)^{1/2}\xi_n
|\widehat{w}(\xi',\xi_n,\tau)|^2}{(1+|\tau+\xi'^2+\xi_n^2|)^{2-\varepsilon}}
d\xi_n
d\xi'
d\tau.
\end{align}
Now, we need the following result about the multiplier appearing in  estimate  \eqref{inhom-h-1-I1-est-1}.
\begin{lemma}
\label{inhom-I1-mut-lem}
Let $0\le b<\frac12$ and $0<\varepsilon\le 1-2b$. Then, for any $(\xi',\xi_n,\tau)\in\rr^{n+1}$  we have
\begin{align}
\label{inhom-I1-mut-est}
Q_1(\xi',\xi_n,\tau)
\doteq&
\chi_{1<\xi_n^2\le \frac12|\tau+\xi'^2|}
(\xi',\xi_n, \tau)
\frac{
(1+\xi'^2+|\tau+\xi'^2|)^{s}
(1+|\tau+\xi'^2|)^{1/2}\xi_n
}{(1+|\tau+\xi'^2+\xi_n^2|)^{2-\varepsilon}}
\nonumber
\\
\lesssim&
\begin{cases}
(1+\xi'^2+\xi_n^2)^{s}(1+|\tau+\xi'^2+\xi_n^2|)^{-2b},
\quad
&-\frac12\le s\le \frac12,
\\
\big[(1+\xi'^2+\xi_n^2)^{s}+(1+|\tau|)^{s}\big]
\cdot
(1+|\tau+\xi'^2+\xi_n^2|)^{-2b},
\quad
&s<-\frac12
\,\,
\text{or}
\,\,
s>\frac12.
\end{cases}
\end{align}
\end{lemma}
We prove Lemma \ref{inhom-I1-mut-lem} later. Now, combining it with inequality \eqref{inhom-h-1-I1-est-1}, we obtain the desired estimate 
\eqref{inhom-h1-I-est} for $I_1$.

\vskip0.1in
\noindent
{\it Estimate \eqref{inhom-h1-I-est} for $I_2$.} We begin by estimating 
the $\eta_n$-integral in definition \eqref{inhom-h-1-Bs-I12-def}. Using the Cauchy-Schwarz inequality, for any $\varepsilon>0$,  we get 
\begin{align*}
\Big(
\int_{\eta_n=1}^\infty
\frac{|\widehat{w}_2(\xi',\sqrt{\eta_n},\tau)|}{(1+|\tau+\xi'^2+\eta_n|)\cdot \sqrt{\eta_n}}
d\eta_n
\Big)^2
\lesssim
\frac{1}{\varepsilon}
\int_1^\infty
\chi_{\eta_n> \frac12|\tau+\xi'^2|}
(\xi',\eta_n, \tau)
\cdot
\frac{|\widehat{w}(\xi',\sqrt{\eta_n},\tau)|^2}{(1+|\tau+\xi'^2+\eta_n|)^{1-\varepsilon}\cdot \eta_n}
d\eta_n.
\end{align*}
Then, making the change of variables $\xi_n=\sqrt{\eta_n}$ or $\eta_n=\xi_n^2$, we obtain  
\begin{align*}
%\label{inhom-h1-CS-2}
\Big(
\int_{\eta_n=1}^\infty
\frac{|\widehat{w}_2(\xi',\sqrt{\eta_n},\tau)|}{(1+|\tau+\xi'^2+\eta_n|)\cdot \sqrt{\eta_n}}
d\eta_n
\Big)^2
\lesssim
\frac{1}{\varepsilon}
\int_{\xi_n=1}^\infty
\chi_{\xi_n^2> \frac12|\tau+\xi'^2|}
(\xi',\xi_n, \tau)
\cdot
\frac{|\widehat{w}(\xi',\xi_n,\tau)|^2}{(1+|\tau+\xi'^2+\xi_n^2|)^{1-\varepsilon}\cdot \xi_n}
d\xi_n,
\end{align*}
which combined with definition \eqref{inhom-h-1-Bs-I12-def} gives us that  
\begin{align}
\label{inhom-h-1-I2-est}
\hskip-0.1in
I_2
\lesssim
\frac{1}{\varepsilon}
\int_{\rr^n}
\int_1^\infty
\chi_{\xi_n^2> \frac12|\tau+\xi'^2|}
(\xi',\xi_n, \tau)
\cdot
\frac{
(1+\xi'^2+|\tau+\xi'^2|)^{s}
(1+|\tau+\xi'^2|)^{1/2}
|\widehat{w}(\xi',\xi_n,\tau)|^2}{(1+|\tau+\xi'^2+\xi_n^2|)^{1-\varepsilon}\xi_n}
d\xi_n
d\xi'
d\tau.
\end{align}
Now, we need the following result about the multiplier appearing in the estimate  \eqref{inhom-h-1-I2-est},
whose proof  is similar to that of Lemma \ref{inhom-I1-mut-lem}.
\begin{lemma}
\label{inhom-I2-mut-lem}
Let $0\le b<\frac12$ and $0<\varepsilon\le 1-2b$. Then, for any $(\xi',\xi_n,\tau)\in\rr^{n+1}$  we have
\begin{align}
\label{inhom-I2-mut-est}
Q_2(\xi',\xi_n,\tau)
\doteq&
\chi_{\xi_n>1}(\xi_n)
\chi_{\xi_n^2> \frac12|\tau+\xi'^2|}
(\xi',\xi_n, \tau)
\frac{
(1+\xi'^2+|\tau+\xi'^2|)^{s}
(1+|\tau+\xi'^2|)^{1/2}
}{(1+|\tau+\xi'^2+\xi_n^2|)^{1-\varepsilon}\xi_n}
\nonumber
\\
\lesssim&
\begin{cases}
(1+\xi'^2+\xi_n^2)^{s}(1+|\tau+\xi'^2+\xi_n^2|)^{-2b},
\quad
&-\frac12\le s\le\frac12,
\\
\big[(1+\xi'^2+\xi_n^2)^{s}+(1+|\tau|)^{s}\big] \cdot
(1+|\tau+\xi'^2+\xi_n^2|)^{-2b},
\quad
&s<-\frac12
\,\,
\text{or}
\,\,
s>\frac12.
\end{cases}
\end{align}
\end{lemma}
Combining it with inequality \eqref{inhom-h-1-I2-est} we obtain the desired estimate \eqref{inhom-h1-I-est} for $I_2$ ,
and complete the proof of $\mathcal{B}^s$-norm estimate \eqref{ivp-forced-est-tim-sharp}.
\,\,
$\square$

\vskip0.05in
\noindent
{\bf Proof of Lemma \ref{inhom-I1-mut-lem}.}
By our assumption $1< \xi_n^2\le \frac12 |\tau+\xi'^2|$, we get 
the relations
\begin{align}
\label{inhom-I1-mut-est1}
1
<
\xi_n^2
\lesssim
|\tau+\xi'^2|
\simeq
|\tau+\xi'^2+\xi_n^2|,
\end{align}
which we use in multiplier $Q_1$, defined in \eqref{inhom-I1-mut-est}, and obtain the bound
\begin{align}
\label{inhom-I1-mut-est2}
Q_1(\xi',\xi_n,\tau)
\lesssim
\chi_{1<\xi_n^2\le \frac12|\tau+\xi'^2|}
(\xi',\xi_n, \tau)
\frac{(1+\xi'^2+|\tau+\xi'^2|)^{s}\xi_n}
{(1+|\tau+\xi'^2+\xi_n^2|)^{\frac12+2b}},
\quad
\varepsilon\le 1-2b.
\end{align}
Next, we estimate the numerator $(1+\xi'^2+|\tau+\xi'^2|)^{s}\xi_n$
by  considering the following two microlocalizations.

\vskip0.05in
\noindent
$\bullet$ Microlocalization 1: $|\tau+\xi'^2|\le 10\xi'^2$
\qquad
$\bullet$ Microlocalization 2: $|\tau+\xi'^2|> 10\xi'^2$

\vskip0.05in
\noindent
\underline{Microlocalization 1.} In this situation we have
the relation $\xi_n^2\lesssim |\tau+\xi'^2|\lesssim \xi'^2$, 
and using it we get
\begin{equation}
\label{non-homo-te-I1-mult-1}
\chi_{1<\xi_n^2\le \frac12|\tau+\xi'^2|}
(\xi',\xi_n, \tau)
(1+\xi'^2+|\tau+\xi'^2|)^{s}
\simeq
(1+\xi'^2)^s
\simeq
(1+\xi'^2+\xi_n^2)^s,
\quad
s\in\rr,
\end{equation}
which combined with inequalities \eqref{inhom-I1-mut-est1} and \eqref{inhom-I1-mut-est2} gives estimate  \eqref{inhom-I1-mut-est}  in Microlocalization 1.

\vskip.05in
\noindent
\underline{Microlocalization 2.} 
Since 
 $\xi_n^2\le \frac12|\tau+\xi'^2|$ and $|\tau+\xi'^2|>10\xi'^2$, 
 we get the relations
\begin{equation}
\label{non-homo-te-I1-mult-2}
\xi'^2+\xi_n^2
\lesssim
|\tau+\xi'^2|
\simeq
|\tau+\xi'^2+\xi_n^2|
\simeq
|\tau|,
\end{equation}
and  by inequality \eqref{inhom-I1-mut-est2} to prove estimate \eqref{inhom-I1-mut-est} in this microlocalization, it suffices to show that 
\begin{align}
\label{Bs-Q1-mircro2-est1}
\chi_{1<\xi_n^2\le \frac12|\tau+\xi'^2|}
(\xi',\xi_n, \tau)
(1+\xi'^2+|\tau+\xi'^2|)^{s}\xi_n
\lesssim&
(1+|\tau|)^{s}
(1+|\tau+\xi'^2+\xi_n^2|)^{\frac12},
\quad
s\in\rr,
\end{align}
and
\begin{align}
\label{Bs-Q1-mircro2-est2}
\chi_{1<\xi_n^2\le \frac12|\tau+\xi'^2|}
(\xi',\xi_n, \tau)
(1+\xi'^2+|\tau+\xi'^2|)^{s}\xi_n
\lesssim&
(1+\xi'^2+\xi_n^2)^{s}
(1+|\tau+\xi'^2+\xi_n^2|)^{\frac12},
\quad
s\le \frac12.
\end{align}
Inequality \eqref{Bs-Q1-mircro2-est1} follows from relations \eqref{inhom-I1-mut-est1} and 
\eqref{non-homo-te-I1-mult-2}. 
Next, we prove inequality \eqref{Bs-Q1-mircro2-est2}.
For $s\le 0$, using assumption $\xi_n^2\le \frac12 |\tau+\xi'^2|$ and inequality \eqref{inhom-I1-mut-est1}, we have    estimate  \eqref{Bs-Q1-mircro2-est2}
\begin{equation}
\label{inhom-I1-mult-neg}
\chi_{1<\xi_n^2\le \frac12|\tau+\xi'^2|}
(\xi',\xi_n, \tau)
(1+\xi'^2+|\tau+\xi'^2|)^{s}\xi_n
\le
(1+\xi'^2+\xi_n^2)^{s}
(1+|\tau+\xi'^2+\xi_n^2|)^{\frac12}.
\end{equation}
For $0<s\le \frac12$, using relation \eqref{non-homo-te-I1-mult-2} 
we obtain
\begin{align*}
\frac{(1+\xi'^2+|\tau+\xi'^2|)^{s}\xi_n}
{(1+|\tau+\xi'^2+\xi_n^2|)^{\frac12}}
\lesssim
\frac{(1+\xi'^2+\xi_n^2)^{\frac12}}
{(1+|\tau+\xi'^2+\xi_n^2|)^{\frac12-s}}
=
\frac{(1+\xi'^2+\xi_n^2)^{s}(1+\xi'^2+\xi_n^2)^{\frac12-s}}
{(1+|\tau+\xi'^2+\xi_n^2|)^{\frac12-s}}
\lesssim
(1+\xi'^2+\xi_n^2)^{s},
\end{align*}
which  gives the desired estimate  \eqref{Bs-Q1-mircro2-est2}.
This completes the proof of Lemma \ref{inhom-I1-mut-lem}.
\,\,
$\square$

%\newpage
%
%
%%%%%%%%%%%%%%%%%%%%%
%
%	trilinear est in X^{s,b}
%
%%%%%%%%%%%%%%%%%%%%%
%
%
\section{Proof of  trilinear estimate}
\label{sec:tri-est}
\setcounter{equation}{0}
We begin by proving  the spatial trilinear estimates, which we restate as follows.
\begin{lemma}
\label{tri-est}
Let $n\ge 2$.
If $s>\frac{n}{2}-1$ and  $\max\{\frac38,\frac18n+\frac14-\frac14s\}< b'\le b<\frac12$, then we have 
\begin{equation}
\label{NLS-tri-est}
\|f\bar{g}h\|_{X^{s,-b}}
\lesssim
\|f\|_{X^{s,b'}}\|g\|_{X^{s,b'}}\|h\|_{X^{s,b'}},
\hspace{1em}
f, g, h \in X^{s,b'}(\rr^{n+1}).
\end{equation}
\end{lemma}
\noindent
{\bf Proof of Lemma \ref{tri-est}.}
Using  the  convenient notation in \cite{b1993-nls}
\begin{equation}
\label{cu-def}
c_u
\doteq
(1+|\xi|)^s
(1+|\tau+\xi^2|)^{b'}
|\widehat{u}(\xi,\tau)|,
\end{equation}
we can rewrite the Bourgain norm as an $L^2$ norm, that is
\begin{equation}
\label{Bourgain-L2-norms}
\|u\|_{X^{s,b'}}
=
\|c_u\|_{L^2_{\xi,\tau}}.
\end{equation}
Now, taking into consideration the identity
$
\widehat{\overline{g}}(\eta ,\tau_2) 
=
\overline{\widehat{g}}(-\eta ,-\tau_2),
$
we get the formula
\begin{align}
\label{convolution-triple}
\widehat{f\bar{g}h}(\xi,\tau)
=
\int_{\rr^{2(n+1)}} \widehat{f}(\zeta ,\tau_1) \overline{\widehat{g}}(-\eta ,-\tau_2) \widehat{h}(\xi-\zeta -\eta ,\tau-\tau_1-\tau_2) d\zeta  d\eta  d\tau_1 d\tau_2,
\end{align}
which we use together with notation \eqref{cu-def} 
to derive the following bound for the Bourgain norm
$\|f\bar{g}h\|_{X^{s,-b}}$ 
\begin{align*}
\|f\bar{g}h\|_{X^{s,-b}}^2
\le&
\int_{\rr^{n+1}}
\Big(
\int_{\rr^{2(n+1)}}
\frac{(1+|\xi|)^{s}}{(1+|\tau+\xi^2|)^{b}}
\frac{c_f(\zeta ,\tau_1)}{
(1+|\zeta |)^{s}
(1+|\tau_1+\zeta ^2|)^{b'}}
\frac{c_g(-\eta ,-\tau_2)}{
(1+|\eta |)^{s}
(1+|\tau_2-\eta ^2|)^{b'}}
\\
\times&
\frac{c_h(\xi-\zeta -\eta ,\tau-\tau_1-\tau_2)}{
(1+|\xi-\zeta -\eta |)^{s}
(1+|\tau-\tau_1-\tau_2+(\xi-\zeta -\eta )^2|)^{b'}}
d\eta d\tau_2
d\zeta d\tau_1
\Big)^2
d\xi
d\tau.
\end{align*}
Furthermore,
combining all multipliers, we arrive at the 
{\it $L^2$-formulation} of 
the trilinear estimate  
\eqref{NLS-tri-est} 
\begin{align*}
%\label{tril-est-l2}
\Big\|
\int_{\rr^{2(n+1)}}
\widetilde{Q}
\,\,
%(\xi,\tau,\zeta ,\tau_1,\eta ,\tau_2)
c_f(\zeta ,\tau_1)
c_g(-\eta ,-\tau_2)
c_h(\xi-\zeta -\eta ,\tau-\tau_1-\tau_2)
d\eta d\tau_2
d\zeta d\tau_1
\Big\|_{L^2}
\lesssim
\|
c_f
\|_{L^2}
\|
c_g
\|_{L^2}
\|
c_h
\|_{L^2},
\end{align*}
where  
$
\widetilde{Q}=\widetilde{Q}(\xi,\tau,\zeta ,\tau_1,\eta ,\tau_2)
$
is the combined multiplier given by
\begin{align*}
%\label{mult-tilde-Q-def}
\widetilde{Q}
\doteq&
\frac{
(1+|\xi|)^{s}
(1+|\zeta |)^{-s}
(1+|\eta |)^{-s}
(1+|\xi-\zeta -\eta |)^{-s}
}{
(1+|\tau+\xi^2|)^{b}
(1+|\tau_1+\zeta ^2|)^{b'}
(1+|\tau_2-\eta ^2|)^{b'}
(1+|\tau-\tau_1-\tau_2+(\xi-\zeta -\eta )^2|)^{b'}
}.
\end{align*}
Also, since $b'\le b$ we bound $\widetilde{Q}$ by the new multiplier $Q$ having only $b'$ exponents
in the denominator
\begin{equation}
\label{mult-Q-def}
Q
\doteq
\frac{
(1+|\xi|)^{s}
(1+|\zeta |)^{-s}
(1+|\eta |)^{-s}
(1+|\xi-\zeta -\eta |)^{-s}
}{
(1+|\tau+\xi^2|)^{b'}
(1+|\tau_1+\zeta ^2|)^{b'}
(1+|\tau_2-\eta ^2|)^{b'}
(1+|\tau-\tau_1-\tau_2+(\xi-\zeta -\eta )^2|)^{b'}
}.
\end{equation}
In addition, since
$
\|u(-\cdot)\|_{L^2}
=
\|u\|_{L^2},
$
we see that in order to prove trilinear estimate \eqref{NLS-tri-est}, it suffices to prove the  following more convenient $L^2$ inequality
\begin{align}
\label{tril-est-l2-no1}
\Big\|
\int_{\rr^{2(n+1)}}
Q
\,\,
%(\xi,\tau,\zeta ,\tau_1,\eta ,\tau_2)
c_f(\zeta ,\tau_1)
c_g(\eta ,\tau_2)
c_h(\xi-\zeta -\eta ,\tau-\tau_1-\tau_2)
d\eta d\tau_2
d\zeta d\tau_1
\Big\|_{L^2}
\lesssim
\|
c_f
\|_{L^2}
\|
c_g
\|_{L^2}
\|
c_h
\|_{L^2}.
\end{align}
We begin the proof of  \eqref{tril-est-l2-no1} by 
using duality and Fubini's theorem to get
\begin{align}
\label{u-tri-L2-duality}
&
\Big\|
\int_{\rr^{2(n+1)}}
Q
\,\,
c_{f}(\zeta ,\tau_1)
c_{g}(\eta ,\tau_2)
c_{h}(\xi-\zeta -\eta ,\tau-\tau_1-\tau_2)
d\eta d\tau_2
d\zeta d\tau_1
\Big\|_{L^2_{\xi,\tau}}
\nn
\\
=&
\sup\limits_{\|d\|_{L^2}=1}
\int_{\rr^{3(n+1)}}
Q
\,\,
d(\xi,\tau)
c_{f}(\zeta ,\tau_1)
c_{g}(\eta ,\tau_2)
c_{h}(\xi-\zeta -\eta ,\tau-\tau_1-\tau_2)
d\xi
d\tau
d\eta d\tau_2
d\zeta d\tau_1.
\end{align}
Moreover, by the symmetry of $\widehat{f}$
and $\widehat{h}$
in the  convolution
$\widehat{f}*
\widehat{\overline{g}} 
* \widehat{h}$
we may assume that 
\begin{equation}
\label{cov-assump}
|\zeta|
\le
|\xi-\zeta-\eta|.
\end{equation}
Therefore, to establish the $L^2$ trilinear estimate   \eqref{tril-est-l2-no1}, it suffices to prove the following  estimate
\begin{align}
\label{tri-L2-duality-split}
\hskip-0.1in
\sup\limits_{\|d\|_{L^2}=1}
\hskip-0.03in
\int_{A}
\hskip-0.03in
Q
d(\xi,\tau)
c_{f}(\zeta ,\tau_1)
c_{g}(\eta ,\tau_2)
c_{h}(\xi
\hskip-0.01in
-
\hskip-0.01in
\zeta
\hskip-0.01in
-
\hskip-0.01in
\eta ,\tau
\hskip-0.01in
-
\hskip-0.01in
\tau_1
\hskip-0.01in
-
\hskip-0.01in
\tau_2)
d\xi
d\tau
d\eta d\tau_2
d\zeta d\tau_1
%\nonumber
%\\
\hskip-0.03in
\lesssim
\hskip-0.03in
\|
c_{f}
\|_{L^2}
\|
c_{g}
\|_{L^2}
\|
c_{h}
\|_{L^2},
\end{align}
where the region $A$ is given by
\begin{equation}
\label{region-A-def}
A
\doteq
\{
(\xi,\tau,\zeta,\tau_1,\eta,\tau_2)\in\rr^{3(n+1)}:
|\zeta|
\le
|\xi-\zeta-\eta|
\}.
\end{equation}

Now, we will consider the following 
 two microlocalizations for $\eta$ appearing via $\widehat{\bar{g}}$.

\vskip0.05in
\noindent
{\bf Microlocalization I.} It is defined by the following sub-region of $A$
\begin{equation}
\label{region-A1-def}
A_1
\doteq
\{
(\xi,\tau,\zeta,\tau_1,\eta,\tau_2)\in \rr^{3(n+1)}:
|\zeta|
\le
|\xi-\zeta-\eta|
\quad
\text{and}
\quad
|\eta|
\le
|\xi-\zeta-\eta|
\}.
\end{equation}
{\bf Microlocalization II.} It is the remaining from $A_1$
 sub-region of $A$
\begin{equation}
\label{region-A2-def}
A_2
\doteq
\{
(\xi,\tau,\zeta,\tau_1,\eta,\tau_2)\in  \rr^{3(n+1)}:
|\zeta|
\le
|\xi-\zeta-\eta|
<
|\xi-\zeta-\eta|
\}.
\end{equation}
\vskip0.01in
\nin
{\bf Proof in microlocalization I.}  Since in $A_1$ we have $|\zeta| \le |\xi-\zeta-\eta|$ and $|\eta| \le |\xi-\zeta-\eta|$, applying the  triangle inequality, we get
$
|\xi|
\le
|\xi-\zeta-\eta|
+
|\zeta|
+
|\eta|
\le
3|\xi-\zeta-\eta|,
$
and  assuming that $s\ge 0$,  this implies that 
$
(1+|\xi|)^{s}\cdot  (1+|\xi-\zeta-\eta|)^{-s}\le 1.
$ 
Thus, we can bound the numerator of  $Q$ in \eqref{mult-Q-def} as follows
\begin{equation*}
\frac{(1+|\xi|)^{s}}
{
(1+|\zeta|)^{s}
(1+|\eta|)^{s}
(1+|\xi-\zeta-\eta|)^{s}
}
\lesssim
\frac{1}
{
(1+|\zeta|)^{s}
(1+|\eta|)^{s}
}
=
(1+|\zeta|)^{-s}
(1+|\eta|)^{-s},
\quad
s\ge0,
\end{equation*}
and using it we see that  in order to prove  trilinear estimate \eqref{tri-L2-duality-split}, it suffices to show that 
\begin{align}
\label{duality-Micro1}
&\sup\limits_{\|d\|_{L^2}=1}
\int_{\rr^{3(n+1)}}
%\hskip-0.1in
Q_1
d(\xi,\tau)
c_{f}(\zeta ,\tau_1)
c_{g}(\eta ,\tau_2)
c_{h}(\xi-\zeta -\eta ,\tau-\tau_1-\tau_2)
d\xi
d\tau
d\eta d\tau_2
d\zeta d\tau_1
\nonumber
\\
\lesssim&
\|
c_{f}
\|_{L^2}
\|
c_{g}
\|_{L^2}
\|
c_{h}
\|_{L^2},
\end{align}
where the new  multiplier  $Q_1
=
Q_1(\xi,\tau,\zeta ,\tau_1,\eta ,\tau_2)$ is defined as follows
\begin{align}
\label{mult-Q1-def}
Q_1
\doteq&
\frac{
(1+|\zeta |)^{-s}
(1+|\eta |)^{-s}
}{
(1+|\tau+\xi^2|)^{b'}
(1+|\tau_1+\zeta ^2|)^{b'}
(1+|\tau_2-\eta ^2|)^{b'}
(1+|\tau-\tau_1-\tau_2+(\xi-\zeta -\eta )^2|)^{b'}
}.
\end{align}
Making the change of variables $\widetilde{\xi}=\xi-\zeta $, $\widetilde{\tau}=\tau-\tau_1$ we obtain
\begin{align}
\label{duality-Micro1-no1}
\hskip-0.08in
\text{LHS of }
\eqref{duality-Micro1}
\hskip-0.03in
\lesssim
\hskip-0.03in
\sup\limits_{\|d\|_{L^2}=1}
\hskip-0.03in
\int_{\rr^{3(n+1)}}
\hskip-0.03in
Q_1
\,
d(\widetilde{\xi}
\hskip-0.01in
+
\hskip-0.01in
\zeta ,\widetilde{\tau}
\hskip-0.01in
+
\hskip-0.01in
\tau_1)
c_{f}(\zeta ,\tau_1)
c_{g}(\eta ,\tau_2)
c_{h}(\widetilde{\xi}
\hskip-0.01in
-
\hskip-0.01in
\eta ,\widetilde{\tau}
\hskip-0.01in
-
\hskip-0.01in
\tau_2)
d\widetilde{\xi}
d\widetilde{\tau}
d\eta d\tau_2
d\zeta d\tau_1,
\end{align}
where the multiplier $Q_1$ becomes $Q_1(\widetilde{\xi},\widetilde{\tau},\zeta ,\tau_1,\eta ,\tau_2)$ 
\begin{align}
\label{Q1-transform}
Q_1
=
\frac{
(1+|\zeta |)^{-s}
(1+|\eta |)^{-s}
}{
(1+|\widetilde{\tau}+\tau_1+(\widetilde{\xi}+\zeta )^2|)^{b'}
(1+|\tau_1+\zeta ^2|)^{b'}
(1+|\tau_2-\eta ^2|)^{b'}
(1+|\widetilde{\tau}-\tau_2+(\widetilde{\xi}-\eta)^2|)^{b'}
}.
\end{align}
Furthermore, following the idea in \cite{tao2001} we rewrite  the multiplier $Q_1$
 as the product of two multipliers 
 %$m_1$ and $m_2$
%
%
\begin{align}
\label{Q1-split}
Q_1
=
m_1(\widetilde{\xi},\widetilde{\tau},\zeta ,\tau_1)
m_2(\widetilde{\xi},\widetilde{\tau},\eta ,\tau_2),
\end{align}
where the multipliers $m_1$, $m_2$  are
\begin{align}
\label{m1-mult}
m_1(\widetilde{\xi},\widetilde{\tau},\zeta ,\tau_1)
\doteq&
\frac{
(1+|\zeta |)^{-s}
}{
(1+|\tau_1+\zeta ^2|)^{b'}
(1+|\widetilde{\tau}+\tau_1+(\widetilde{\xi}+\zeta )^2|)^{b'}
},
\\
\label{m2-mult}
m_2
(\widetilde{\xi},\widetilde{\tau},\eta ,\tau_2)
\doteq&
\frac{
(1+|\eta|)^{-s}
}{
(1+|\tau_2-\eta ^2|)^{b'}
(1+|\widetilde{\tau}-\tau_2+(\widetilde{\xi}-\eta )^2|)^{b'}
}.
\end{align}
Using \eqref{Q1-split}, Fubini, and
 Cauchy-Schwarz inequality in $d\widetilde{\xi} d\widetilde{\tau}$-integration, from \eqref{duality-Micro1-no1} we get 
\begin{subequations}
\label{tri-dual-1-est-3}
\begin{align}
\label{tri-dual-1-est-3-a}
\text{LHS of }
\eqref{duality-Micro1}
\lesssim
\sup\limits_{\|d\|_{L^2}=1}
\,
&
\Big\|
\int_{\rr^{n+1}}
m_1(\widetilde{\xi},\widetilde{\tau},\zeta ,\tau_1)
c_{f}(\zeta ,\tau_1)
d(\widetilde{\xi}+\zeta ,\widetilde{\tau}+\tau_1)
d\zeta  d\tau_1
\Big\|_{L^2_{\tilde{\xi},\tilde{\tau}}}
\\
\label{tri-dual-1-est-3-b}
\times&
\Big\|
\int_{\rr^{n+1}}
m_2
(\widetilde{\xi},\widetilde{\tau},\eta ,\tau_2)
c_{g}(\eta ,\tau_2)
c_{h}(\widetilde{\xi}-\eta ,\widetilde{\tau}-\tau_2)
d\eta  d\tau_2
\Big\|_{L^2_{\tilde{\xi},\tilde{\tau}}}.
\end{align}
\end{subequations}
%
%Now, for s and $b'$ satisfying the conditions given in  Lemma %\ref{tri-est}, 
Now, we  claim that, 
under the assumptions of Lemma \ref{tri-est},
the following two bilinear estimate hold
\begin{align}
&\label{m1-est}
\Big\|
\int_{\rr^{n+1}}
m_1(\xi,\tau,\zeta ,\tau_1)
c_{f}(\zeta ,\tau_1)
d(\xi+\zeta ,\tau+\tau_1)
d\zeta  d\tau_1
\Big\|_{L^2_{\xi,\tau}}
\lesssim
\|c_f\|_{L^2} \|d\|_{L^2},
\\
\label{m2-est}
&\Big\|
\int_{\rr^{n+1}}
m_2
(\xi,\tau,\eta ,\tau_2)
c_{g}(\eta ,\tau_2)
c_{h}(\xi-\eta,\tau-\tau_2)
d\eta  d\tau_2
\Big\|_{L^2_{\xi,\tau}}
\lesssim
\|c_g\|_{L^2} \|c_h\|_{L^2},
\end{align}
which used together with \eqref{tri-dual-1-est-3} 
yield the desired estimate \eqref{duality-Micro1}.
Thus, we complete the proof of $L^2$ trilinear estimate \eqref{tril-est-l2-no1} in Microlocalization I, once we prove 
the bilinear estimates \eqref{m1-est} and \eqref{m2-est}.

\vskip0.05in
\noindent
{\bf Proof of  bilinear estimate \eqref{m1-est}.} 
For this we consider the following two new microlocalizations.

\vskip0.05in
\noindent
\underline{Microlocalization 1.} 
In estimate \eqref{m1-est},
it restricts $d\zeta  d\tau_1$- integration in  the region 
\begin{equation}
\label{region-B1-def}
B_1
=
B_1(\xi,\tau)
\doteq
\{
(\zeta,\tau_1)\in \rr^{n+1}:
|\xi|
\ge
|\zeta|
\}.
\end{equation}
\underline{Microlocalization 2.} 
This restricts $d\zeta  d\tau_1$- integration in  the region 
\begin{equation}
\label{region-B2-def}
B_2
=
B_2(\xi,\tau)
\doteq
\{
(\zeta,\tau_1)\in \rr^{n+1}:
|\xi|
<
|\zeta|
\}.
\end{equation}
\underline{Proof in microlocalization 1.}
Then,  inequality \eqref{m1-est} reduces to the following one
\begin{align}
\label{m1-est-B1}
\Big\|
\int_{B_1}
m_1(\xi,\tau,\zeta ,\tau_1)
c_{f}(\zeta ,\tau_1)
d(\xi+\zeta ,\tau+\tau_1)
d\zeta  d\tau_1
\Big\|_{L^2_{\xi,\tau}}
\lesssim
\|c_f\|_{L^2} \|d\|_{L^2}.
\end{align}
Now, applying the Cauchy-Schwarz inequality 
for $d\zeta d\tau_1$-integration followed by taking the sup-norm for 
$
(\int_{\rr^{n+1}}m_1^2
d\zeta d\tau_1)^{1/2}, 
$
from  \eqref{m1-est-B1} we get
\begin{align*}
\text{LHS of\,\,}
\eqref{m1-est-B1}
\le
\Big\|
\Big(
\int_{B_1}
m_1^2
d\zeta d\tau_1 
\Big)^{1/2}
\Big\|_{L_{\xi,\tau}^\infty}
\|c_f\|_{L^2} \|d\|_{L^2}.
\end{align*}
Therefore, in order to prove  inequality \eqref{m1-est-B1}, it suffices to show the following result.
\begin{lemma}
\label{theta1-lem}
Let $n\ge 2$.
If $s>\frac{n}{2}-1$ and  $\max\{\frac38,\frac18n+\frac14-\frac14s\}< b'<\frac12$, then for $(\xi,\tau)\in\rr^{n+1}$  we have 
\begin{align}
\label{theta1-est}
\Theta_1(\xi,\tau)
\doteq
\int_{B_1}
m_1^2
d\tau_1
d\zeta 
=
\int_{\rr^{n+1}}
\frac{\chi_{|\xi|
\ge
|\zeta|}(\zeta)\cdot (1+|\zeta |)^{-2s}\, d\tau_1 d\zeta 
}{(1+|\tau_1+\zeta^{2}|)^{2{b'}}(1+|\tau+\tau_1+(\xi+\zeta)^2|)^{2{b'}}}
\lesssim
1.
\end{align}
\end{lemma}
\noindent
To prove this lemma, we need the following calculus result, which can be found in \cite{et2016}.
\begin{lemma}
\label{calc-lem-et}
If $\ell\ge \ell'\ge 0$ and $2\ell+2\ell'>1$, then 
\begin{equation}
\label{eq:calc_7}
\int_{\rr}\frac{dx}{(1+|x-a|)^{2\ell}(1+|x-c|)^{2\ell'}}
\lesssim
(1+|a-c|)^{-2\ell'}\phi_{\ell}(a-c),
\end{equation}
where
\begin{align*}
\phi_{\ell}(k)
\doteq
\begin{cases}
1,
\quad
&\ell>\frac12,
\\
\log(1+|k|),
\quad
&\ell=\frac12,
\\
(1+|k|)^{1-2\ell},
\quad
&\ell<\frac12.
\end{cases}
\end{align*}

\end{lemma}

\noindent
{\bf Proof of Lemma \ref{theta1-lem}.} 
Applying calculus inequality 
\eqref{eq:calc_7} with $x=\tau_1$, $a=-\zeta^2$ and $c=-\tau-(\xi+\zeta)^2$, for $1/4<\ell=\ell'={b'}<1/2$ we get 
\begin{align*}
%\label{theta-est-1}
\Theta_1(\xi,\tau)
\lesssim
\int_{\rr^n}
\frac{\chi_{|\xi|
\ge
|\zeta|}(\zeta)\cdot (1+|\zeta |)^{-2s} }{(1+|\tau+(\xi+\zeta)^2-\zeta^{2}|)^{4b'-1}}
d\zeta
=
\int_{\rr^n}
\frac{\chi_{|\xi|
\ge
|\zeta|}(\zeta)\cdot (1+|\zeta |)^{-2s}}{(1+|2\xi\cdot \zeta+\xi^2+\tau|)^{4b'-1}}
d\zeta
=
I_0
+
I_1,
\end{align*}
where the integrals $I_0$ and $I_1$ are given by
\begin{align*}
I_0(\xi,\tau)
\doteq
\int_{|\zeta|\le 1}
\frac{\chi_{|\xi|
\ge
|\zeta|}(\zeta)\cdot (1+|\zeta |)^{-2s}}{(1+|2\xi\cdot \zeta+\xi^2+\tau|)^{4b'-1}}
d\zeta,
\quad
I_1(\xi,\tau)
\doteq
\int_{|\zeta|> 1}
\frac{\chi_{|\xi|
\ge
|\zeta|}(\zeta)\cdot (1+|\zeta |)^{-2s}}{(1+|2\xi\cdot \zeta+\xi^2+\tau|)^{4b'-1}}
d\zeta.
\end{align*}
If $4b'-1\ge 0$  and $-2s\le 0$, then we have 
$
I_0(\xi,\tau)
\le
\int_{|\zeta|\le 1}
1
d\zeta
\lesssim 
1.
$
Next, we show that  integral $I_1(\xi,\tau)$ is bounded and we consider the following two cases.

\vskip0.05in
\noindent
Case 1: $|\tau|\ge 4\xi^2$
\quad
Case 2: $|\tau|< 4\xi^2$

\vskip0.05in
\noindent
\underline{Case 1.} 
Using assumptions
$
|\xi|\ge |\zeta|,
$ 
and
$|\tau|\ge 4\xi^2$,
we have 
$
|2\xi\cdot \zeta+\xi^2+\tau|
\ge
|\tau|
-
\xi^2
-
2|\xi||\zeta|
\ge
\xi^2,
$
and we bound $I_1$ as follows
$$
I_1(\xi,\tau)
\lesssim
\int_{|\zeta|> 1}
\frac{\chi_{|\xi|
\ge
|\zeta|}(\zeta)\cdot (1+|\zeta |)^{-2s}}{(1+\xi^2)^{4b'-1}}
d\zeta
\lesssim
\int_{\rr^n}
\frac{1}{(1+|\zeta|)^{2(4b'-1)+2s}}
d\zeta.
$$
This is bounded if $2(4b'-1)+2s>n$ or 
$
b'
>
\frac18n+\frac14-\frac14s.
$
Since $b'<\frac12$, we need $\frac18n+\frac14-\frac14s<\frac12$ or 
\begin{equation}
\label{s-critical}
s>\frac{n}{2}-1.
\end{equation}
This shows that $I_1$ is bounded in Case 1.

\vskip0.05in
\noindent
\underline{Case 2.} 
Since  
%$\max\limits_{j}|\xi_j|\ge \frac{1}{\sqrt{n}}|\xi|$, 
we can cover  $\rr^n$ with the conic regions
defined by
$
|\xi_j|
>
\frac{1}{n}|\xi|,
$
we may assume 
\begin{equation}
\label{zeta1-assump}
|\xi_1|
\gtrsim
|\xi|.
\end{equation}
Also, making the change of variables $(\zeta_1,\zeta_2,\cdots, \zeta_n)\to (\mu, \zeta_2,\cdots,\zeta_n)$, where $\mu=2\xi\cdot \zeta+\xi^2+\tau$, with the Jacobian
$
J
=
\Big|
\frac{\p(\mu, \zeta_2,\cdots,\zeta_n) }{\p (\zeta_1,\zeta_2,\cdots, \zeta_n)}
\Big|
=
2\xi_1
$
we obtain 
\begin{align}
\label{tri-microI-I1}
I_1(\xi,\tau)
\lesssim
\int_{|\zeta|> 1}
\frac{\chi_{|\xi|
\ge
|\zeta|}(\zeta)\cdot (1+|\zeta |)^{-2s}}{(1+|\mu|)^{4b'-1} |J|}
d\mu d\zeta_2\cdots d\zeta_n.
\end{align}
Since the integral $I_1$ is  over  $\{\zeta\in\rr^n: 1<|\zeta|\le |\xi|\}$, using  assumptions 
$
|\xi_1|
\gtrsim
|\xi|
$ and   $|\tau|< 4\xi^2$ we get 
$$
|\mu|
\le
|\tau|
+
\xi^2
+
2|\xi||\zeta|
\le
7\xi^2
\lesssim
\xi_1^2
\quad
\text{and}
\quad
|J|=2|\xi_1|\gtrsim |\xi|\gtrsim |\zeta|
>
1.
$$
Combining the above inequalities with estimate \eqref{tri-microI-I1} we have
\begin{align*}
I_1(\xi,\tau)
\lesssim
\int_{\rr^{n-1}}
\Big(
\int_{|\mu|\lesssim \xi_1^2}
\frac{1}{(1+|\mu|)^{4b'-1}}
d\mu
\Big)
\frac{\chi_{|\xi_1|\gtrsim 1+|\zeta'|} (\zeta')}{(1+|\zeta'|)^{2s}  |\xi_1| }
 d\zeta',
\end{align*}
where $\zeta'=(\zeta_2,\cdots,\zeta_n)$. Thus, for $\frac38\le b'<\frac12$, we obtain
\begin{align*}
I_1
\lesssim
\int_{\rr^{n-1}}
|\xi_1|^{2(2-4b')}
\frac{\chi_{|\xi_1|\gtrsim 1+|\zeta'|} (\zeta')}{(1+|\zeta'|)^{2s}  |\xi_1| }
 d\zeta'
=
\int_{\rr^{n-1}}
\frac{\chi_{|\xi_1|\gtrsim 1+|\zeta'|} (\zeta')}{(1+|\zeta'|)^{2s}  |\xi_1|^{8b'-3} }
d\zeta'
\lesssim
\int_{\rr^{n-1}}
\frac{1}{(1+|\zeta'|)^{2s+8b'-3} }
d\zeta',
\end{align*}
which is bounded if $2s+8b'-3>n-1$ or $
b'
>
\frac18n+\frac14-\frac14s.
$
Since $b'<\frac12$, we must have $\frac18n+\frac14-\frac14s<\frac12$ or 
$
s>\frac{n}{2}-1,
$
which is condition \eqref{s-critical}. This completes the proof of   Lemma \ref{theta1-lem}.
\,\,
$\square$

\vskip0.05in
\noindent
\underline{Proof in microlocalization 2.}
%involves $L^4$-norm estimate
%
Then  inequality \eqref{m1-est} becomes 
\begin{align}
\label{m1-est-B2}
\Big\|
\int_{B_2}
m_1(\xi,\tau,\zeta ,\tau_1)
c_{f}(\zeta ,\tau_1)
d(\xi+\zeta ,\tau+\tau_1)
d\zeta  d\tau_1
\Big\|_{L^2_{\xi,\tau}}
\lesssim
\|c_f\|_{L^2} \|d\|_{L^2}.
\end{align}
We prove this by writing $m_1$ as a product of two similar multipliers.
Since in $B_2$ we have $|\xi|<|\zeta|$ we get  $|\xi+\zeta|\le 2|\zeta|$ and
$
(1+|\zeta |)
\gtrsim
(1+|\zeta |)^{1/2}
(1+|\xi+\zeta |)^{1/2}.
$
Thus, we bound the multiplier $m_1$ as follows
$$
m_1(\xi,\tau,\zeta ,\tau_1)
\doteq
\frac{
(1+|\zeta |)^{-s}
}{
(1+|\tau_1+\zeta ^2|)^{b'}
(1+|\tau+\tau_1+(\xi+\zeta )^2|)^{b'}
}
\lesssim
\frac{
(1+|\zeta |)^{-s/2}
(1+|\xi+\zeta |)^{-s/2}
}{
(1+|\tau_1+\zeta ^2|)^{b'}
(1+|\tau+\tau_1+(\xi+\zeta )^2|)^{b'}
}.
$$
Also, we make change of variables $\zeta\to -\zeta$ and $\tau_1\to -\tau_1$ and we get 
\begin{align}
\label{m1-convo}
&\int_{B_2}
m_1(\xi,\tau,\zeta ,\tau_1)
c_{f}(\zeta ,\tau_1)
d(\xi+\zeta ,\tau+\tau_1)
d\zeta d\tau_1
\nonumber
\\
\lesssim&
\int_{\rr^{n+1}}
\frac{(1+|\zeta |)^{-s/2} c_{f}(-\zeta,-\tau_1)}{(1+|\tau_1-\zeta^2|)^{b'}}
\frac{(1+|\xi-\zeta |)^{-s/2} d(\xi-\zeta,\tau-\tau_1)}{(1+|\tau-\tau_1+(\xi-\zeta)^2|)^{b'}}
d\zeta d\tau_1
=
\widehat{F}
*
\widehat{G},
\end{align}
where the  functions $F$, $G$ (to be shown that are in $L^4$)
 are defined via the inverse Fourier transform
\begin{align*}
%\label{F-inverse}
F(x,t)
\doteq&
\frac{1}{(2\pi)^{n+1}}
\int_{\rr^{n+1}}
e^{i\xi\cdot x+i\tau t}
\frac{
(1+|\xi|)^{-s/2}
c_f(-\xi,-\tau)}{(1+|\tau-\xi^{2}|)^{b'}}
d\xi
d\tau,
\qquad
x\in\rr^n,
\,\,
t\in\rr,
\\
%\label{G-inverse}
G(x,t)
\doteq&
\frac{1}{(2\pi)^{n+1}}
\int_{\rr^{n+1}}
e^{i\xi\cdot x+i\tau t}
\frac{
(1+|\xi|)^{-s/2}
d(\xi,\tau)}{(1+|\tau+\xi^{2}|)^{b'}}
d\xi
d\tau,
\qquad
x\in\rr^n,
\,\,
t\in\rr.
\end{align*}
Now, following \cite{b1993-nls, hy2024-JMA-cNLS}, we apply Parseval's formula and then Cauchy-Schwarz inequality in  $dxdt$
to get
\begin{align}
\label{L2-general-est-3}
\|
\widehat{F}
*
\widehat{G}
\|_{L^2}
\simeq
\|
\widehat{F G}
\|_{L^2}
\simeq
\|
F G
\|_{L^2}
\le
\Big(
\int_{\rr^{n+1}}
|F|^4
dxdt
\Big)^{1/4}
\Big(
\int_{\rr^{n+1}}
|G|^4
dxdt
\Big)^{1/4}
=
\|F\|_{L^4}
\|G\|_{L^4}.
\end{align}
Combining inequality  \eqref{m1-convo} with relation \eqref{L2-general-est-3} we get
\begin{align}
\label{L2-general-est-4}
\Big\|
\int_{B_2}
m_1(\xi,\tau,\zeta ,\tau_1)
c_{f}(\zeta ,\tau_1)
d(\xi+\zeta ,\tau+\tau_1)
d\zeta d\tau_1
\Big\|_{L_{\xi,\tau}^2}
\lesssim
\|F\|_{L^4_{x,t}}
\|G\|_{L^4_{x,t}}.
\end{align}
Finally, we bound $\|F\|_{L^4_{x,t}}$ and $\|G\|_{L^4_{x,t}}$
by the corresponding $L^2$-norms by using the  following 
 estimates for Schr\"odinger  type equations
 (see \cite{a2018, b1993-nls, b1998, Ctbook, cdks2001, mp2015, s1977}).
% or  Strichartz type estimate.
%
\begin{lemma}
\label{stric-lem}
Let $n\ge 2$.
If $s>\frac{n}{2}-1$ and  $\max\{\frac14,\frac18n+\frac14-\frac14s\}< \rho<\frac12$, 
then for   $L^2$ function $f(\xi,\tau)$, $\xi\in\rr^n,\tau\in\rr$ we have 
\begin{align}
\label{stric-est-1}
&\Big\|
\int_{\rr^{n+1}}
e^{i\xi\cdot x+i\tau t}
\frac{
(1+|\xi|)^{-s/2}
f(\xi,\tau)}{(1+|\tau+\xi^{2}|)^{\rho}}
d\xi
d\tau
\Big\|_{L^4_{x,t}}
\lesssim
\|f\|_{L^2_{\xi,\tau}},
\\
\label{stric-est-2}
&\Big\|
\int_{\rr^{n+1}}
e^{i\xi\cdot x+i\tau t}
\frac{
(1+|\xi|)^{-s/2}
f(\xi,\tau)}{(1+|\tau-\xi^{2}|)^{\rho}}
d\xi
d\tau
\Big\|_{L^4_{x,t}}
\lesssim
\|f\|_{L^2_{\xi,\tau}}.
\end{align}
\end{lemma}
We prove  Lemma \ref{stric-lem} later.
Now, applying inequalities \eqref{stric-est-2} and \eqref{stric-est-1}  with $f=c_f(-\cdot,-\cdot)$ and  $f=d$ respectively,  we obtain
$
\|F\|_{L^4_{x,t}}
\lesssim
\|c_f(-\cdot,-\cdot)\|_{L^2_{\xi,\tau}}
=
\|c_f\|_{L^2_{\xi,\tau}}
$
and
$
\|G\|_{L^4_{x,t}}
\lesssim
\|d\|_{L^2_{\xi,\tau}}.
$
Therefore, from inequality \eqref{L2-general-est-4} we get the desired bilinear estimate
 \eqref{m1-est-B2} once we prove Lemma \ref{stric-lem}.
\,\,
$\square$

%
%%%%%%%%%%%%%%%%%%%%%
%
%
%
%	Proof of Strichartz estimate
%
%
%%%%%%%%%%%%%%%%%%%%%
%
%
\vskip0.1in
\nin
{\bf Proof of  Lemma \ref{stric-lem}.}
\label{sub-sec:Strichartz}
Since the second estimate follows from the first by a change 
of variables,
we only provide the proof of inequality \eqref{stric-est-1}.
Defining the function $F$ via inverse Fourier transform
\begin{align}
\label{f1-inverse}
F(x,t)
\doteq
\frac{1}{(2\pi)^{n+1}}
\int_{\rr^{n+1}}
e^{i\xi\cdot x+i\tau t}
\frac{
(1+|\xi|)^{-s/2}
f(\xi,\tau)}{(1+|\tau+\xi^{2}|)^{\rho}}
d\xi
d\tau,
\quad
x\in\rr^n,
\,\,
t\in\rr,
\end{align}
and we see that to prove inequality \eqref{stric-est-1}, it suffices to show  $\|F\|_{L^4_{x,t}}\lesssim \|f\|_{L^2_{\xi,\tau}}$.  
 Applying Parseval's formula we get
\begin{align}
\label{Bourgain-L4-express}
\|F\|_{L^4_{x,t}}^4
=
\int_{\rr^{n+1}}
|F
\cdot
F|^2
dxdt
=
\|F F\|_{L^2_{x,t}}^2
\simeq
\|
\widehat{F F}
\|_{L_{\xi,\tau}^2}^2
\simeq
\|
\widehat{F}
*
\widehat{F}
\|_{L_{\xi,\tau}^2}^2.
\end{align}
From the definition of $F$ in  \eqref{f1-inverse} we get
$
\widehat{F}(\xi,\tau)
=
\frac{
(1+|\xi|)^{-s/2}
f(\xi,\tau)}{(1+|\tau+\xi^{2}|)^{\rho}},
$
which gives
\begin{align*}
\widehat{F}*\widehat{F}
(\xi,\tau)
=
\int_{\rr^{n+1}}
\frac{
(1+|\zeta|)^{-s/2}
f(\zeta,\tau_1)}{(1+|\tau_1+\zeta^{2}|)^{\rho}}
\frac{
(1+|\xi-\zeta|)^{-s/2}
f(\xi-\zeta,\tau-\tau_1)}
{
(1+|\tau-\tau_1+(\xi-\zeta)^2|)^{\rho}
}
d\zeta d\tau_1,
\end{align*}
and therefore  we can write \eqref{Bourgain-L4-express} 
as follows  
\begin{align}
\label{trilinear-est-L2-form}
\|F\|_{L^4_{x,t}}^4
\simeq
\Big\|
\int_{\rr^{n+1}} 
Q(\xi,\tau,\zeta,\tau_1)
f(\zeta,\tau_1)
f(\xi-\zeta,\tau-\tau_1) 
d\zeta d\tau_1 \Big\|_{L^2_{\xi,\tau}}^2,
\end{align}
where $Q=Q(\xi,\tau,\zeta,\tau_1)$ is the multiplier defined as follows
\begin{equation*}
%\label{mult-Q-def}
Q(\xi,\tau,\zeta,\tau_1)
\doteq
\frac{(1+|\zeta|)^{-s/2}(1+|\xi-\zeta|)^{-s/2}}{(1+|\tau_1+\zeta^{2}|)^{\rho}(1+|\tau-\tau_1+(\xi-\zeta)^2|)^{\rho}}.
\end{equation*}
Now,   in \eqref{trilinear-est-L2-form}, 
applying the Cauchy-Schwarz inequality in 
 $d\zeta d\tau_1$ and taking the sup over $\xi, \tau$ 
we get
\begin{align*}
\|F\|_{L^4_{x,t}}^4
\lesssim
\Big\|
\Big(
\int_{\rr^{n+1}}
Q^2
d\zeta d\tau_1 
\Big)^{\frac12}
\Big\|_{L_{\xi,\tau}^\infty}^2
\|f\|_{L^2_{\xi,\tau}}^4.
\end{align*}
Therefore, in order to prove  inequality \eqref{stric-est-1}, it suffices to show the following result.
\begin{lemma}
\label{theta2-lem}
Let $n\ge 2$.
If $s>\frac{n}{2}-1$ and  $\max\{\frac14,\frac18n+\frac14-\frac14s\}< \rho<\frac12$, then for $(\xi,\tau)\in\rr^{n+1}$ we have 
\begin{align}
\label{theta2-est}
\Theta_2(\xi,\tau)
\doteq
\int_{\rr^{n+1}}
Q^2
d\tau_1
d\zeta 
=
\int_{\rr^{n+1}}
\frac{(1+|\zeta |)^{-s}(1+|\xi-\zeta |)^{-s}}{(1+|\tau_1+\zeta^{2}|)^{2{\rho}}(1+|\tau-\tau_1+(\xi-\zeta)^2|)^{2{\rho}}}
d\tau_1
d\zeta 
\lesssim
1.
\end{align}
\end{lemma}
\noindent
{\bf Proof of Lemma \ref{theta2-lem}.} 
Applying calculus inequality 
\eqref{eq:calc_7} with $x=\tau_1$, $a=-\zeta^2$ and $c=\tau+(\xi-\zeta)^2$, for $1/4<\ell=\ell'={\rho}<1/2$ we get 
\begin{align}
\label{theta2-est-1}
\Theta_2
\lesssim
\int_{\rr^n}
\frac{ (1+|\zeta |)^{-s} (1+|\xi-\zeta |)^{-s} }{(1+|\tau+(\xi-\zeta)^2+\zeta^{2}|)^{4\rho-1}}
d\zeta.
\end{align}
Now, we notice that the case $
|\xi-\zeta|
\le
|\zeta|
$
can be reduced to the case 
$
|\zeta|
\le
|\xi-\zeta|,
$
by
 making the change of variables $\tilde{\zeta}=\xi-\zeta$. 
Thus, we assume that 
$
|\zeta|
\le
|\xi-\zeta|.
$
Then, we have 
$
|\zeta|^2
\le
|\zeta||\xi-\zeta|,
$
which gives
 $
 (1+|\zeta |)^{-s} (1+|\xi-\zeta |)^{-s}\le (1+|\zeta |)^{-2s}
 $
 and we get the following bound  for $\Theta_2$ 
\begin{align*}
\Theta_2
\lesssim
\int_{\rr^n}
\frac{ (1+|\zeta |)^{-2s}}{(1+|2\zeta^2-2\xi\cdot \zeta+\xi^2+\tau|)^{4\rho-1}}
d\zeta
=
I_0
+
I_1,
\quad
s\ge 0,
\end{align*}
where the integrals $I_0(\xi,\tau)$ and $I_1(\xi,\tau)$ are defined by
\begin{align}
I_0
\doteq
\int_{|\zeta|\le 1}
\frac{ (1+|\zeta |)^{-2s}}{(1+|2\zeta^2-2\xi\cdot \zeta+\xi^2+\tau|)^{4\rho-1}}
d\zeta,
\quad
I_1
\doteq
\int_{|\zeta|> 1}
\frac{ (1+|\zeta |)^{-2s}}{(1+|2\zeta^2-2\xi\cdot \zeta+\xi^2+\tau|)^{4\rho-1}}
d\zeta.
\end{align}
If $4\rho-1\ge 0$  and $-2s\le 0$, then  
$
I_0(\xi,\tau)
\lesssim 
1.
$
Next, we estimate integral $I_1$ by using the decomposition
\begin{equation}
\label{tri-I1-split}
I_1(\xi,\tau)
=
I_{11}(\xi,\tau)
+
I_{11}(\xi,\tau),
\quad
\text{with}
\quad
I_{1j}
\doteq
\int_{D_j}
\frac{ (1+|\zeta |)^{-2s}}{(1+|2\zeta^2-2\xi\cdot \zeta+\xi^2+\tau|)^{4\rho-1}}
d\zeta,
\quad
j=1,2,
\end{equation}
where the regions $D_1$ and $D_2$ are defined as follows
\begin{align}
\label{region-D1-def}
D_1
\doteq
\{
\zeta\in \rr^{n}:
|\zeta|> 1
\,\,\,
\text{and}
\,\,\,
|2\zeta^2-2\xi\cdot \zeta+\xi^2+\tau|
\ge
|\zeta|^2
\},
\\
\label{region-D2-def}
D_2
\doteq
\{
\zeta\in \rr^{n}:
|\zeta|> 1
\,\,\,
\text{and}
\,\,\,
|2\zeta^2-2\xi\cdot \zeta+\xi^2+\tau|
<
|\zeta|^2
\}.
\end{align}
\underline{Estimate for $I_{11}$.} 
Since   
$
|2\zeta^2-2\xi\cdot \zeta+\xi^2+\tau|
\ge
|\zeta|^2,
$ 
we have
\begin{align*}
I_{11}(\xi,\tau)
\lesssim
\int_{D_1}
\frac{ (1+|\zeta |)^{-2s}}{(1+|\zeta|^2)^{4\rho-1}}
d\zeta
\lesssim
\int_{\rr^n}
\frac{ 1}{(1+|\zeta|)^{2(4\rho-1)+2s}}
d\zeta,
\end{align*}
which is bounded if $2(4\rho-1)+2s>n$ or 
$
\rho
>
\frac18n+\frac14-\frac14s.
$
Since $\rho<\frac12$, we must have $\frac18n+\frac14-\frac14s<\frac12$ or 
$
s>\frac{n}{2}-1,
$
which is condition \eqref{s-critical}.
This shows that 
 $I_{11}$ is bounded.

\vskip0.05in
\noindent
\underline{Estimate for $I_{12}$.} 
We prove it when  $n\ge 3$, since the case $n=2$ is similar.
Like in the proof of Lemma \ref{theta1-lem}, we 
may assume that
\begin{equation}
\label{zeta3-assump}
|\zeta_1|
>
\frac{1}{n}|\zeta|.
\end{equation}
Also, recalling the notation  $\zeta^2=\zeta_1^2+\cdots\zeta_n^2$, $\xi\cdot \zeta=\xi_1\zeta_1+\cdots+\xi_n\zeta_n$, and completing squares  we get
$$
2\zeta^2-2\xi\cdot \zeta+\xi^2+\tau
=
2(\zeta_{n-1}
-
\frac12\xi_{n-1})^2
+
2(\zeta_n-\frac12\xi_n)^2
+
2\zeta'^2
-
2\xi'\cdot \zeta'
+
\frac12\xi_{n-1}^2
+
\frac12\xi_n^2
+
\xi'^2
+
\tau,
$$
where $\zeta'=(\zeta_1,\cdots,\zeta_{n-2})$ and $\xi'=(\xi_1,\cdots,\xi_{n-2})$. Combining the  above identity we write $I_{12}$ as
follows
\begin{equation}
\label{I12-def}
I_{12}
=
\int_{\rr^{n-2}}
\tilde{I}_{12}
(\xi,\tau,\zeta')
d\zeta',
\end{equation}
where  the integral $\tilde{I}_{12}=\tilde{I}_{12}
(\xi,\tau,\zeta')$ is defined by
\begin{align*}
\tilde{I}_{12}
\doteq
\int_{\rr^2}
\frac{ 
\chi_{D_2}(\zeta)
\cdot
(1+|\zeta |)^{-2s}\quad  d\zeta_{n-1}d\zeta_n
}{(1+
|
2(\zeta_{n-1}
-
\frac12\xi_{n-1})^2
+
2(\zeta_n-\frac12\xi_n)^2
+
2\zeta'^2
-
2\xi'\cdot \zeta'
+
\frac12\xi_{n-1}^2
+
\frac12\xi_n^2
+
\xi'^2
+
\tau
|)^{4\rho-1}}.
\end{align*}
Next, we estimate $\tilde{I}_{12}$. Using polar coordinates $\zeta_{n-1}-\frac12\xi_{n-1}=r\cos\theta$, $\zeta_n-\frac12\xi_n=r\sin\theta$ we obtain
\begin{align*}
\tilde{I}_{12}
(\xi,\tau,\zeta')
=&
\int_{0}^\infty
\int_{0}^{2\pi}
\frac{ 
\chi_{D_2}(r,\zeta')
\cdot
\Big[1+
\big|
(r\cos\theta+\frac12\xi_{n-1})^2
+
(r\sin\theta+\frac12\xi_{n})^2
+\zeta'^2\big|^{\frac12}
\,
\Big]^{-2s}
\cdot
r
}{(1+
|
2r^2
+
2\zeta'^2
-
2\xi'\cdot \zeta'
+
\frac12\xi_{n-1}^2+\frac12\xi_n^2
+\xi'^2+\tau
|)^{4\rho-1}}
d\theta
dr
\\
\lesssim&
\int_{0}^\infty
\frac{ 
\chi_{D_2}(r,\zeta')
\cdot
(1+|\zeta'|)^{-2s}
\cdot
r
}{(1+
|
2r^2
+
2\zeta'^2
-
2\xi'\cdot \zeta'
+
\frac12\xi_{n-1}^2+\frac12\xi_n^2
+\xi'^2+\tau
|)^{4\rho-1}}
dr,
\end{align*}
where by   assumption \eqref{zeta3-assump},  the region $D_2$ given by \eqref{region-D2-def} becomes
$$
D_2
=
D_2(r, \zeta')
\doteq
\{
r\in\rr^+,
\,\,
\zeta'\in\rr^{n-2}:
|
2r^2
+
2\zeta'^2
-
2\xi'\cdot \zeta'
+
\frac12\xi_{n-1}^2+\frac12\xi_n^2
+\xi'^2+\tau
|
\lesssim
\zeta'^2
\}.
$$
Now, using the change of variables $\mu=\mu(r)=2r^2
+
2\zeta'^2
-
2\xi'\cdot \zeta'
+
\frac12\xi_{n-1}^2+\frac12\xi_n^2
+\xi'^2+\tau$, we get 
\begin{align}
\label{tilde-I12}
\tilde{I}_{12}
\lesssim&
\int_{|\mu|\lesssim \zeta'^2}
\frac{ 
(1+|\zeta'|)^{-2s}
\cdot
r
}{(1+
|
\mu
|)^{4\rho-1} \cdot r}
dr
\lesssim
(1+|\zeta'|)^{-2s}
(1+|\zeta'|)^{2(2-4\rho)}
=
\frac{1}{(1+|\zeta'|)^{2s-2(2-4\rho)}}
,
\quad
\rho<\frac12.
\end{align}
Finally, combining estimate \eqref{tilde-I12} with \eqref{I12-def}, we obtain
$
I_{12}
\lesssim
\int_{\rr^{n-2}}
\frac{1}{(1+|\zeta'|)^{2s-2(2-4\rho)}}
d\zeta',
$
which is bounded if $2s-2(2-4\rho)>n-2$ or  
$
\rho
>
\frac18n+\frac14-\frac14s.
$
Since $\rho<\frac12$, we must have 
$
s>\frac{n}{2}-1,
$
which is condition \eqref{s-critical}.
This shows that  $I_{12}$ is bounded  and completes the proof of   Lemma \ref{theta2-lem}.
\,\,
$\square$

\vskip0.05in
\noindent
{\bf Proof of bilinear estimate \eqref{m2-est}.} 
We reduce it to estimate \eqref{m1-est} for the $m_1$ multiplier,
which we have proved above.
For this, we make the change of variables $\widetilde{\eta}=-\eta$,  $\widetilde{\tau}_2=-\tau_2$ and we get
\begin{align*}
\int_{\rr^{n+1}}
\hskip-0.03in
m_2
(\xi,\tau,\eta ,\tau_2)
c_{g}(\eta ,\tau_2)
c_{h}(\xi-\eta,\tau-\tau_2)
d\eta  d\tau_2
=
\hskip-0.05in
\int_{\rr^{n+1}}
\hskip-0.03in
m_1(\xi,\tau,\widetilde{\eta},\widetilde{\tau}_2)
c_{g}(-\widetilde{\eta},-\widetilde{\tau}_2)
c_h(\xi+\widetilde{\eta},\tau+\widetilde{\tau}_2)
d\widetilde{\eta}d\widetilde{\tau}_2,
\end{align*}
where  the multipliers $m_1$, $m_2$ are defined in \eqref{m1-mult}, \eqref{m2-mult}. 
Now applying estimate  \eqref{m1-est} for $m_1$ we get
\begin{align*}
\text{LHS of }
\eqref{m2-est}
=
\big\|
\int_{\rr^{n+1}}
m_1(\xi,\tau,\widetilde{\eta},\widetilde{\tau}_2)
c_{g}(-\widetilde{\eta},-\widetilde{\tau}_2)
c_h(\xi+\widetilde{\eta},\tau+\widetilde{\tau}_2)
d\widetilde{\eta}d\widetilde{\tau}_2
\big\|_{L^2}
\overset{\eqref{m1-est}}{\lesssim}
\|c_g\|_{L^2}
\|c_h\|_{L^2},
\end{align*}
where we have used the fact
$
\|c_g(-\cdot,-\cdot)\|_{L^2}=\|c_g\|_{L^2}.
$
This completes the proof of inequality  \eqref{m2-est}.
\,\,
$\square$

%%%%%%%%%%%%%%%%%%%
%
%
% Microlocalization II
%
%
%%%%%%%%%%%%%%%%%%%
%
%
\vskip0.05in
\noindent
{\bf Proof of trilinear estimate \eqref{tri-L2-duality-split} in microlocalization II.} 
Like in the case of   Microlocalization I we reduce it 
to two bilinear estimate one with multiplier $m_1$ and the other
with a new multiplier $m_3$.
To accomplish this, in  region $A_2$, defined by \eqref{region-A2-def}, and where we have
$|\zeta| \le |\xi-\zeta-\eta|\le |\eta|$, we use inequality 
$
|\xi|
\le
|\xi-\zeta-\eta|
+
|\zeta|
+
|\eta|
\le
3|\eta| 
$
to simplify the numerator of  $Q$ defined in \eqref{mult-Q-def} as follows 
$$
\frac{(1+|\xi|)^{s}}
{
(1+|\zeta|)^{s}
(1+|\eta|)^{s}
(1+|\xi-\zeta-\eta|)^{s}
}
\lesssim
\frac{1}
{
(1+|\zeta|)^{s}
(1+|\xi-\zeta-\eta|)^{s}
}
=
(1+|\zeta|)^{-s}
(1+|\xi-\zeta-\eta|)^{-s},
\quad
s\ge 0.
$$
Thus, to prove  trilinear estimate \eqref{tri-L2-duality-split} in this microlocalization, it suffices to show that 
\begin{align}
\label{duality-Micro2}
&
\sup\limits_{\|d\|_{L^2}=1}
\int_{\rr^{3(n+1)}}
Q_2
\,\,
d(\xi,\tau)
c_{f}(\zeta ,\tau_1)
c_{g}(\eta ,\tau_2)
c_{h}(\xi-\zeta -\eta ,\tau-\tau_1-\tau_2)
d\xi
d\tau
d\eta d\tau_2
d\zeta d\tau_1
\nonumber
\\
\lesssim&
\|
c_{f}
\|_{L^2}
\|
c_{g}
\|_{L^2}
\|
c_{h}
\|_{L^2},
\end{align}
where the multiplier  $Q_2
=
Q_2(\xi,\tau,\zeta ,\tau_1,\eta ,\tau_2)$ is defined as follows
\begin{align}
\label{mult-Q2-def}
Q_2
\doteq&
\frac{
(1+|\zeta |)^{-s}
(1+|\xi-\zeta-\eta|)^{-s}
}{
(1+|\tau+\xi^2|)^{b'}
(1+|\tau_1+\zeta ^2|)^{b'}
(1+|\tau_2-\eta ^2|)^{b'}
(1+|\tau-\tau_1-\tau_2+(\xi-\zeta -\eta )^2|)^{b'}
}.
\end{align}
Working similarly to the derivation of inequality \eqref{tri-dual-1-est-3} in microlocalization I, we get  
\begin{subequations}
\label{tri-micro2-1}
\begin{align}
\label{tri-micro2-1-a}
\text{LHS of }
\eqref{duality-Micro2}
\lesssim
\sup\limits_{\|d\|_{L^2}=1}
\,
&
\Big\|
\int_{\rr^{n+1}}
m_1(\widetilde{\xi},\widetilde{\tau},\zeta ,\tau_1)
c_{f}(\zeta ,\tau_1)
d(\widetilde{\xi}+\zeta ,\widetilde{\tau}+\tau_1)
d\zeta  d\tau_1
\Big\|_{L^2_{\tilde{\xi},\tilde{\tau}}}
\\
\label{tri-micro2-1-b}
\times&
\Big\|
\int_{\rr^{n+1}}
m_3
(\widetilde{\xi},\widetilde{\tau},\eta ,\tau_2)
c_{g}(\eta ,\tau_2)
c_{h}(\widetilde{\xi}-\eta ,\widetilde{\tau}-\tau_2)
d\eta  d\tau_2
\Big\|_{L^2_{\tilde{\xi},\tilde{\tau}}},
\end{align}
\end{subequations}
where the multiplier $m_1$ is defined by \eqref{m1-mult}
and the new multiplier $m_3$ is defined as follows
\begin{align}
\label{m3-mult}
m_3
(\widetilde{\xi},\widetilde{\tau},\eta ,\tau_2)
\doteq&
\frac{
(1+|\widetilde{\xi}-\eta|)^{-s}
}{
(1+|\tau_2-\eta ^2|)^{b'}
(1+|\widetilde{\tau}-\tau_2+(\widetilde{\xi}-\eta )^2|)^{b'}
}.
\end{align}
Also,  making the change of variables  $\widetilde{\xi}-\eta=\widetilde{\eta}$ and $\widetilde{\tau}-\tau_2=\widetilde{\tau}_2$, and letting    $\xi=-\widetilde{\xi}$ and $\tau=-\widetilde{\tau}$, we get 
\begin{align*}
&\int_{\rr^{n+1}}
\hskip-0.05in
m_3
(\widetilde{\xi},\widetilde{\tau},\eta ,\tau_2)
c_{g}(\eta ,\tau_2)
c_{h}(\widetilde{\xi}-\eta ,
\hskip-0.02in
\widetilde{\tau}-\tau_2)
d\eta  d\tau_2
=
\hskip-0.08in
\int_{\rr^{n+1}}
\hskip-0.05in
m_1(\xi,\tau,\widetilde{\eta},\widetilde{\tau}_2)
c_{h}(\widetilde{\eta} ,\widetilde{\tau}_2)
c_g(
\hskip-0.02in
-
\hskip-0.02in
(\xi
\hskip-0.01in
+
\hskip-0.01in
\widetilde{\eta}),
\hskip-0.03in
-(\tau
\hskip-0.01in
+
\hskip-0.01in
\widetilde{\tau}_2))
d\widetilde{\eta}  d\widetilde{\tau}_2,
\end{align*}
in which $m_3$ is replaced by $m_1$.
Combining this with  inequality \eqref{tri-micro2-1} we get
\begin{align*}
\text{LHS of }
\eqref{duality-Micro2}
\lesssim
\sup\limits_{\|d\|_{L^2}=1}
\,
&
\Big\|
\int_{\rr^{n+1}}
m_1(\widetilde{\xi},\widetilde{\tau},\zeta ,\tau_1)
c_{f}(\zeta ,\tau_1)
d(\widetilde{\xi}+\zeta ,\widetilde{\tau}+\tau_1)
d\zeta  d\tau_1
\Big\|_{L^2_{\tilde{\xi},\tilde{\tau}}}
\\
\times&
\Big\|
\int_{\rr^{n+1}}
m_1(\xi,\tau,\widetilde{\eta},\widetilde{\tau}_2)
c_{h}(\widetilde{\eta} ,\widetilde{\tau}_2)
c_g(-(\xi+\widetilde{\eta}),-(\tau+\widetilde{\tau}_2))
d\widetilde{\eta}  d\widetilde{\tau}_2
\Big\|_{L^2_{\xi,\tau}}.
\end{align*}
Finally, applying inequality \eqref{m1-est} for $m_1$ twice and using 
$
\|c_g(-\cdot,-\cdot)\|_{L^2}=\|c_g\|_{L^2}
$ 
we obtain the desired  inequality \eqref{duality-Micro2}.
This completes the proof of spatial trilinear estimate \eqref{NLS-tri-est}.
\,\,
$\square$

%
%
%%%%%%%%%%%%%%%%%%%%%
%
%
%
%	Proof of trilinear estimate in Y space
%
%
%%%%%%%%%%%%%%%%%%%%%
%
%
\section{Proof of temporal  trilinear estimates}
\label{sec:tri-est-Y}
\setcounter{equation}{0}
Here, we prove the NLS trilinear estimates 
in temporal Bourgain spaces, that is 
estimates \eqref{trilinear-est-Y}, which 
we restate in the following lemma.
\begin{lemma}
\label{trilinest-Y}
Let $n\ge 2$.
If $s>\frac{n}{2}-1$ and  $\max\{\frac38,\frac18n+\frac14-\frac14s\}< b'\le b<\frac12$, then we have 
\begin{equation}
\label{NLS-tri-est-Y}
\|f\bar{g}h\|_{Y^{s,-b}}
\lesssim
(\|f\|_{X^{s,b'}}+\|f\|_{Y^{s,b'}})
(\|g\|_{X^{s,b'}}+\|g\|_{Y^{s,b'}})
(\|h\|_{X^{s,b'}}+\|h\|_{Y^{s,b'}}),
\hspace{1em}
f, g, h \in X^{s,b'}\cap Y^{s,b'}.
\end{equation}
\end{lemma}

\vskip0.05in
\noindent
{\bf Proof of Lemma \ref{trilinest-Y}.}
Using   the following notation, 
which is similar to the one used in $X^{s,b}$-norm,
\begin{equation}
\label{cu-def-Y}
c_u
\doteq
[(1+|\xi|)^s+(1+|\tau|)^{s/2}]
(1+|\tau+\xi^2|)^{b'}
|\widehat{u}(\xi,\tau)|,
\end{equation}
we write the sum  norm 
$
\|u\|_{X^{s,b'}}+\|u\|_{Y^{s,b'}}
$
as an $L^2$ norm, that is
\begin{equation}
\label{Bourgain-L2-norms-Y}
\|u\|_{X^{s,b'}}
+
\|u\|_{Y^{s,b'}}
\simeq
\|c_u\|_{L^2_{\xi,\tau}}.
\end{equation}
Using this notation, and working like in the derivation of 
the $L^2$ formulation  \eqref{tri-L2-duality-split},
we arrive at the following  $L^2$ inequality,
which implies trilinear estimate  \eqref{NLS-tri-est-Y},
\begin{align}
\label{tri-L2-duality-split-Y}
&
\sup\limits_{\|d\|_{L^2}=1}
\int_{A}
Q
\,\,
d(\xi,\tau)
c_{f}(\zeta ,\tau_1)
c_{g}(\eta ,\tau_2)
c_{h}(\xi-\zeta -\eta ,\tau-\tau_1-\tau_2)
d\xi
d\tau
d\eta d\tau_2
d\zeta d\tau_1
\nonumber
\\
\lesssim&
\|
c_{f}
\|_{L^2}
\|
c_{g}
\|_{L^2}
\|
c_{h}
\|_{L^2},
\end{align}
where the region $A$, by the symmetry of $\widehat{f}$
and $\widehat{h}$
in the  convolution
$\widehat{f}*
\widehat{\overline{g}} 
* \widehat{h}$, 
is given as follows
\begin{equation}
\label{region-A-def-Y}
A
\doteq
\{
(\xi,\tau,\zeta,\tau_1,\eta,\tau_2)\in\rr^{3(n+1)}:
|\tau_1|
\le
|\tau-\tau_1-\tau_2|
\},
\end{equation}
and the multiplier $Q
=
Q(\xi,\tau,\zeta ,\tau_1,\eta ,\tau_2)$ is defined as follows
\begin{subequations}
\label{Y-mult-Q-def}
\begin{align}
\label{Y-mult-Q-def-no1}
Q
\doteq&
\frac{(1+|\tau|)^{s/2}}
{
[(1+|\zeta |)^{s}+(1+|\tau_1|)^{\frac{s}{2}}]
[(1+|\eta |)^{s}+(1+|\tau_2|)^{\frac{s}{2}}]
[(1+|\xi-\zeta -\eta |)^{s}+
(1+|\tau-\tau_1-\tau_2|)^{\frac{s}{2}}]
}
\\
\label{Y-mult-Q-def-no2}
\times&
\frac{
1
}{
(1+|\tau+\xi^2|)^{b'}
(1+|\tau_1+\zeta ^2|)^{b'}
(1+|\tau_2-\eta ^2|)^{b'}
(1+|\tau-\tau_1-\tau_2+(\xi-\zeta -\eta )^2|)^{b'}
}.
\end{align}
\end{subequations}

Next, we prove temporal inequality \eqref{tri-L2-duality-split-Y},
by considering the following two ``temporal" microlocalizations,
in analogy to spatial  microlocalizations used in the proof
of our spatial trilinear estimate.
 
\vskip0.05in
\noindent
{\bf Microlocalization I.} We restrict the integration in \eqref{tri-L2-duality-split-Y} to the region $A_1$, which is defined by
\begin{equation}
\label{region-A1-def-Y}
A_1
\doteq
\{
(\xi,\tau,\zeta,\tau_1,\eta,\tau_2)\in\rr^{3(n+1)}:
|\tau_1|
\le
|\tau-\tau_1-\tau_2|
\quad
\text{and}
\quad
|\tau_2|
\le
|\tau-\tau_1-\tau_2|
\}.
\end{equation}
{\bf Microlocalization II.} We restrict to the region $A_2$, which is defined as follows
\begin{equation}
\label{region-A2-def-Y}
A_2
\doteq
\{
(\xi,\tau,\zeta,\tau_1,\eta,\tau_2)\in\rr^{3(n+1)}:
|\tau_1|
\le
|\tau-\tau_1-\tau_2|
<
|\tau_2|
\}.
\end{equation}
\underline{Proof in microlocalization I.} In region $A_1$, since $|\tau_1|
\le
|\tau-\tau_1-\tau_2|$ and $|\tau_2|
\le
|\tau-\tau_1-\tau_2|$, 
we get 
$
|\tau|
\le
|\tau-\tau_1-\tau_2|
+
|\tau_1|
+
|\tau_2|
\le
3|\tau-\tau_1-\tau_2|.
$
Using this we obtain
$
(1+|\tau|)^{s/2}
[(1+|\xi-\zeta -\eta |)^{s}+
(1+|\tau-\tau_1-\tau_2|)^{s/2}]^{-1}
\lesssim
1
$
and we bound the first fraction (involving $s$) in the multiplier $Q$ as follows
\begin{align}
\label{Y-mult-Q1-est}
\eqref{Y-mult-Q-def-no1}
\lesssim&
\frac{1}
{
[(1+|\zeta |)^{s}+(1+|\tau_1|)^{s/2}]
[(1+|\eta |)^{s}+(1+|\tau_2|)^{s/2}]
}
\le
(1+|\zeta|)^{-s}
(1+|\eta|)^{-s},
\quad
s\ge0.
\end{align}
Combining inequality \eqref{Y-mult-Q1-est} with definition \eqref{Y-mult-Q-def}, we obtain 
$Q\lesssim Q_1$,  
where  $Q_1$  is the multiplier defined in \eqref{mult-Q1-def}. 
Therefore, estimate \eqref{tri-L2-duality-split-Y} reduces to  
$L^2$ estimate \eqref{duality-Micro1}.

\vskip0.05in
\noindent
\underline{Proof in  microlocalization II.} 
Then, we have $|\tau_1|
\le
|\tau-\tau_1-\tau_2|< |\tau_2|$. Working like
in  \eqref{Y-mult-Q1-est}, we get
\begin{align}
\label{Y-mult-Q2-est}
\eqref{Y-mult-Q-def-no1}
\lesssim
(1+|\zeta|)^{-s}
(1+|\xi-\zeta -\eta |)^{-s},
\quad
s\ge0.
\end{align}
Combining inequality \eqref{Y-mult-Q2-est} with definition \eqref{Y-mult-Q-def}, we obtain 
$Q\lesssim Q_2$, 
where $Q_2$ is defined in \eqref{mult-Q2-def}. 
Therefore, we reduce  estimate \eqref{tri-L2-duality-split-Y} to  \eqref{duality-Micro2}. This
completes the proof of trilinear estimate \eqref{NLS-tri-est-Y}.
\,\,
$\square$

%
%%%%%%%%%%%%%%%%%%%%%%
%
%
% Optimality
%
%
%%%%%%%%%%%%%%%%%%%%%%  
%
\section{optimality of  trilinear estimates}
\label{sec:tri-counter}
\setcounter{equation}{0}

In this section,  we conclude our study of NLS trilinear estimates
by proving their optimality.
More precisely, we have the following result.

\begin{lemma}
\label{counter-example}
Let $n\ge 2$.
If $s\le \frac{n}{2}-1$, then the trilinear estimate \eqref{NLS-tri-est} with $b'=b$ fails for any  $b<\frac12$.
\end{lemma}

\vskip0.05in
\noindent
{\bf Proof of Lemma \ref{counter-example}.}
Using  notation \eqref{cu-def},
the trilinear estimate \eqref{NLS-tri-est} with $b'=b$
reads
\begin{align}
\label{tril-est-l2-counter}
\Big\|
\int_{\rr^{2(n+1)}}
Q
\,
%(\xi,\tau,\zeta ,\tau_1,\eta ,\tau_2)
c_f(\zeta ,\tau_1)
c_g(\eta ,\tau_2)
c_h(\xi-\zeta -\eta ,\tau-\tau_1-\tau_2)
d\eta d\tau_2
d\zeta d\tau_1
\Big\|_{L^2_{\xi,\tau}}
\lesssim
\|
c_f
\|_{L^2}
\|
c_g
\|_{L^2}
\|
c_h
\|_{L^2},
\end{align}
where the multiplier $Q=Q(\xi,\tau,\zeta ,\tau_1,\eta ,\tau_2)$ is defined as follows
\begin{equation}
\label{mult-Q-counter}
Q
\doteq
\frac{
(1+|\xi|)^{s}
(1+|\zeta |)^{-s}
(1+|\eta |)^{-s}
(1+|\xi-\zeta -\eta |)^{-s}
}{
(1+|\tau+\xi^2|)^{b}
(1+|\tau_1+\zeta ^2|)^{b}
(1+|\tau_2-\eta ^2|)^{b}
(1+|\tau-\tau_1-\tau_2+(\xi-\zeta -\eta )^2|)^{b}
}.
\end{equation}
We will show  that if $b<\frac12$ and $L^2$
trilinear estimate \eqref{tril-est-l2-counter} holds, then   $s>\frac{n}{2}-1$, $n\ge 2$.  To show this, for any
given  $N\in \zz^+$ we  define the three $L^2$ functions  $c_{f}$, $c_g$ and $c_h$ as follows
\begin{equation}
\label{cf-counter-example-function}
c_{f}(\zeta,\tau_1)
\doteq
\chi_{A_1}(\zeta,\tau_1),
\quad
c_{g}(\eta,\tau_2)
\doteq
\chi_{A_2}(\eta,\tau_2)
\quad
\text{and}
\quad
c_{h}(\mu,\tau_3)
\doteq
\chi_{A_3}(\mu,\tau_3),
\end{equation}
where  $\chi_{A_j}(\cdot)$ are the characteristic functions of the regions $A_j$, $j=1,2,3$, which  are 
defined by
\begin{align}
\label{A1-counter-example-domain}
A_1
\doteq&
\Big\{
(\zeta,\tau_1)
\in
\rr^{n+1}
:
|\tau_1+\zeta^2|
\le
1,
\,\,
N
\le
\zeta_k
\le
2N,
\,\,
k=1,2,3,\cdots,n
\Big\},
\\
\label{A2-counter-example-domain}
A_2
\doteq&
\Big\{
(\eta,\tau_2)
\in
\rr^{n+1}
:
|\tau_2-\eta^2|
\le
1,
\,\,
N
\le
\eta_k
\le
2N,
\,\,
k=1,2,3,\cdots,n
\Big\},
\\
\label{A3-counter-example-domain}
A_3
\doteq&
\Big\{
(\mu,\tau_3)
\in
\rr^{n+1}
:
|\tau_3|
\le
47nN^2,
\,\,
N
\le
\mu_k
\le
4N,
\,\,
k=1,2,3,\cdots,n
\Big\}.
\end{align}
Also, for the left hand-side of inequality \eqref{tril-est-l2-counter} we will need  the following region $A_0$ 
\begin{align}
\label{A0-counter-example-domain}
A_0
\doteq
\Big\{
(\xi,\tau)
\in
\rr^{n+1}
:
|\tau+\xi^2|
\le
1,
\,\,
5N
\le
\xi_k
\le
6N,
\,\,
k=1,2,3,\cdots,n
\Big\}.
\end{align}
Now, by straightforward computation we get
\begin{align}
\label{area-A0}
&
\|\chi_{A_0}\|_{L^2}^2
=
\int_{\rr^{n+1}}
\chi_{A_0}(\xi,\tau)
d\xi
d\tau
=
\int_{\xi_n=5N}^{6N}
\cdots
\int_{\xi_1=5N}^{6N}
\int_{\tau=-1-\xi^2}^{1-\xi^2}
1
\,
d\tau
d\xi_1
\cdots
d\xi_n
=
2N^{n},
\\
\label{L2-norm-cf}
&\|c_{f}\|_{L^2}^2
=
\|\chi_{A_1}\|_{L^2}^2
=
\|c_{g}\|_{L^2}^2
=
\|\chi_{A_2}\|_{L^2}^2
\simeq
N^{n}
\quad
\text{and}
\quad
\|c_{h}\|_{L^2}^2
=
\|\chi_{A_3}\|_{L^2}^2
\simeq 
N^{n+2}.
\end{align}
Next, we define the integrand function  $\Theta=\Theta(\xi,\tau)$
\begin{align}
\label{counter-theta}
\Theta
\doteq
\int_{\rr^{2(n+1)}}
Q(\xi,\tau,\zeta ,\tau_1,\eta ,\tau_2)
%(\xi,\tau,\zeta ,\tau_1,\eta ,\tau_2)
c_f(\zeta ,\tau_1)
c_g(\eta ,\tau_2)
c_h(\xi-\zeta -\eta ,\tau-\tau_1-\tau_2)
d\eta d\tau_2
d\zeta d\tau_1,
\end{align}
for which we shall prove that the following key estimate 
\begin{align}
\label{count-theta-bound}
\| 
\Theta
\|_{L^2_{\xi,\tau}}
\gtrsim
N^{-2s-2b+\frac52n}.
\end{align}
However, before proving it, we use it to complete the proof of 
Lemma \ref{counter-example}.
Combining  inequality \eqref{L2-norm-cf} with \eqref{count-theta-bound},  we see that if  the trilinear estimate \eqref{tril-est-l2-counter} holds, then we must  have
$
N^{-2s-2b+\frac52n}
\lesssim
N^{\frac32n+1}.
$
For this to hold when $N\gg 1$, we must have  $-2s-2b+\frac52n\le \frac32n+1$ or
$
\frac12n-s-\frac12\le b.
$
And since $b<\frac12$, we must have $\frac12n-s-\frac12<\frac12$ or
\begin{equation}
\label{2bar-counter-s-cond}
s
>
\frac{n}{2}-1,
\quad
n\ge 2.
\end{equation}
This completes the proof of  Lemma \ref{counter-example}, once we prove inequality \eqref{count-theta-bound}.
\,\,
$\square$

\vskip0.05in
\noindent
{\bf Proof of inequality \eqref{count-theta-bound}.}
We first show that if $(\xi,\tau)\in A_0$, $(\zeta,\tau_1)\in A_1$ and $(\eta,\tau_2)\in A_2$, then
\begin{equation}
\label{count-translate}
(\mu,\tau_3)
\in
A_3,
\quad
\text{where}
\,\,
\mu=\xi-\zeta-\eta
\,\,\,
\text{and}
\,\,\,
\tau_3
=
\tau-\tau_1-\tau_2.
\end{equation}
In fact, then we have $5N\le \xi_k\le 6N$, $N\le \zeta_k\le 2N$ and $N\le \eta_k\le 2N$. Therefore,
$
\mu_k
\doteq
\xi_k-\zeta_k-\eta_k
\ge 5N-2N-2N=N
$ 
and $\mu_k\le 6N-N-N=4N$. Using these two inequalities, we get 
$$
N
\le
\mu_k
\le
4N,
\quad
k=1,2, 3,\cdots,n.
$$
Furthermore, combining the above inequalities with $|\tau+\xi^2|\le 1$, $|\tau_1+\zeta^2|\le 1$ and $|\tau_2-\eta^2|\le 1$  we obtain
\begin{align*}
|\tau_3|
\le&
|\tau|+|\tau_1|+|\tau_2|
\le
|\tau+\xi^2|+\xi^2
+
|\tau_1+\zeta^2|+\zeta^2
+
|\tau_2-\eta^2|+\eta^2
\le
47nN^2,
\end{align*}
which gives 
$
(\mu,\tau_3)\in A_3.
$

Now, using these we bound $Q$ from below.
Since $(\xi,\tau)\in A_0$, $(\zeta,\tau_1)\in A_1$  and $(\eta,\tau_2)\in A_2$, we have
\begin{align*}
|\tau-\tau_1-\tau_2+(\xi-\zeta-\eta)^2|
=&
|(\tau+\xi^2)-(\tau_1+\zeta^2)-(\tau_2-\eta^2)-[\xi^2-\zeta^2+\eta^2-(\xi-\zeta-\eta)^2]|
\\
\le&
|\tau+\xi^2|
+
|\tau_1+\zeta^2|
+
|\tau_2-\eta^2|
+
\xi^2+\zeta^2+\eta^2+(\xi-\zeta-\eta)^2
\lesssim
N^2,
\end{align*}
and we can bound the multiplier $Q$, defined in \eqref{mult-Q-counter},  as follows
$$
Q(\xi,\tau,\zeta ,\tau_1,\eta ,\tau_2)
\gtrsim
\frac{(1+N)^s(1+N)^{-s}(1+N)^{-s}(1+N)^{-s}}{
(1+1)^b(1+1)^b(1+1)^b(1+N^2)^b}
\gtrsim
N^{-2s-2b}.
$$
Finally, combining this inequality  with relations \eqref{area-A0} and \eqref{L2-norm-cf}  we get
\begin{align*}
\| 
\Theta
\|_{L^2_{\xi,\tau}}
\ge&
N^{-2s-2b}
\Big\|
\chi_{A_0}(\cdot,\cdot)
\int_{\rr^{n+1}}
\chi_{A_1}(\zeta,\tau_1)
d\zeta d\tau_1
\int_{\rr^{n+1}}
\chi_{A_2}(\eta,\tau_2)
d\eta d\tau_2
\Big\|_{L^2_{\xi,\tau}}
\simeq
N^{-2s-2b+\frac52n},
\end{align*}
which is  the desired inequality \eqref{count-theta-bound}.
\,\,
$\square$

%
%%%%%%%%%%%%%%%%%%
%
%
% 
%   Proof of WP Theorem 
%  
%
%
%%%%%%%%%%%%%%%%%%%
%
\section{ Proof of well-posedness}
\label{sec:proof-wp}
\setcounter{equation}{0}
In this section, using the linear and trilinear estimates derived
earlier we  prove the well-posedness Theorems \ref{wp-thm-small} and  \ref{wp-thm} via a fixed point argument.

\vskip0.05in
\noindent
{\bf Proof of Theorem \ref{wp-thm-small}.}
Let  $n\ge 2$  and $s>\frac{n}{2}-1$.
For given initial data $u_0\in H^s(\rr^{n-1}\times\rr^+)$ 
and boundary data $g_0\in \mathcal{B}_T^s$, 
we define an iteration map by replacing in the Fokas solution 
formula
  \eqref{nd-lnls-utm-sln}
 the forcing by the nonlinearities $\mp |u|^2u$.
More precisely, for any input $u$  we define the following map
\begin{align}
\label{iteration-map-T-loc-2}
u(x, t)
=
\Phi( u)
\doteq
S
\Big[
u_0,g_0;
\mp |u|^2u
\Big],
\quad
(x, t)
\in \rr^{n-1} \times \rr^+ \times (0, T),\, \,
 T<1/2.
\end{align}
Also, for the data we assume the smallness condition
\begin{equation}
\label{smallness} 
\|u_0\|_{H_x^s(\rr^{n-1}\times\rr^+)}
+
\|g_0\|_{\mathcal{B}_T^s}
\le 
\frac{1}{32C^2},
\quad
\text{with}
\quad
C\doteq1+2c_1+32c_1c_2,
\end{equation}
where $c_1$  
is the constant appearing in the linear estimates 
(Theorem \ref{forced-linear-nls-thm}) 
and $c_2$ is the constant appearing  in the trilinear estimates
(Theorem \ref{Trilinear-estimate-thm}). 
Finally,  in our solution space $X^{s,b}_{\rr^{n-1}\times\rr^+\times(0,T)}
\cap
Y^{s,b}_{\rr^{n-1}\times\rr^+\times(0,T)}
$
with the sum norm
$
|||u|||
\doteq
\|u\|_{X^{s,b}_{\rr^{n-1}\times\rr^+\times(0,T)}}
+
\|u\|_{Y^{s,b}_{\rr^{n-1}\times\rr^+\times(0,T)}}
$
we choose a closed ball $B=B(r)$
\begin{align}
\label{r-ball-small}
B
\doteq
\big\{
u:
u\in X^{s,b}_{\rr^{n-1}\times\rr^+\times(0,T)}
\cap
Y^{s,b}_{\rr^{n-1}\times\rr^+\times(0,T)}
\,\,
\text{and}
\,\,
|||u|||
\le
r\big\}.
\end{align} 
Next, we choose $r$ so that  map $\Phi$ goes from $B(r)$
into  $B(r)$ and is a contraction on the closed ball $B(r)$.

First, we show that $\Phi$ is into $B(r)$.
For $u\in B(r)$, using linear estimates 
\eqref{X-forced-linear-nls-est}, \eqref{Y-forced-linear-nls-est}, 
we  estimate $|||\Phi(u)|||$ as follows
\begin{align}
\label{forced-linear-nls-est-small}
|||\Phi(u)|||
\le&
2c_1\big[
\|u_0\|_{H_x^s(\rr^{n-1}\times\rr^+)}
+
\|g_0\|_{\mathcal{B}_T^s}
+
\|
|u|^2u
\|_{X^{s,-b}_{\rr^{n-1}\times\rr^+\times(0,T)}}
+
\|
|u|^2u
\|_{Y^{s,-b}_{\rr^{n-1}\times\rr^+\times(0,T)}}
 \big]
\nonumber
\\
\le&
2c_1\big[
\|u_0\|_{H_x^s(\rr^{n-1}\times\rr^+)}
+
\|g_0\|_{\mathcal{B}_T^s}
+
\|
|\tilde{u}|^2\tilde{u}
\|_{X^{s,-b}}
+
\|
|\tilde{u}|^2\tilde{u}
\|_{Y^{s,-b}}
 \big],
\end{align}
where $\tilde{u}$ is the extension of $u$ from $\rr^{n-1}\times\rr^+\times(0,T)$ to $\rr^{n+1}$, 
which  satisfies 
$$
%\label{u-extension}
\|\tilde u\|_{X^{s,b}}+
\|\tilde u\|_{Y^{s,b}}
\le
2
\big(
\|u\|_{X^{s,b}_{\rr^{n-1}\times\rr^+\times(0,T)}}
+
\|u\|_{Y^{s,b}_{\rr^{n-1}\times\rr^+\times(0,T)}}
\big).
$$
Next, using trilinear estimates \eqref{trilinear-est-B} and \eqref{trilinear-est-Y} from inequality \eqref{forced-linear-nls-est-small} we obtain
\begin{align}
\label{forced-linear-nls-est-small-1}
|||\Phi(u)|||
\le
2c_1\big[
\|u_0\|_{H_x^s(\rr^{n-1}\times\rr^+)}
+
\|g_0\|_{\mathcal{B}_T^s}
+
16c_2
\big(
\|u\|_{X^{s,b}_{\rr^{n-1}\times\rr^+\times(0,T)}}
+
\|u\|_{Y^{s,b}_{\rr^{n-1}\times\rr^+\times(0,T)}}
\big)^3
\big].
\end{align}
Thus, for $u\in B(r)$ we have 
\begin{equation}
\label{small-onto-ine1}
|||\Phi(u)|||
\le
C
\big(
\|u_0\|_{H_x^s(\rr^{n-1}\times\rr^+)}
+
\|g_0\|_{\mathcal{B}_T^s}
\big)
+
C
r^3,
\end{equation}
where $C$ is the constant defined in  \eqref{smallness}.
Now, choosing 
\begin{equation}
\label{small-r-def}
r
\doteq
2C
\big(
\|u_0\|_{H_x^s(\rr^{n-1}\times\rr^+)}
+
\|g_0\|_{\mathcal{B}_T^s}
\big),
\end{equation}
and using the smallness condition \eqref{smallness}, we get
\begin{equation}
\label{small-r-cond}
r
\le
2C
\cdot
\frac1{32C^2}
\le
\frac{1}{16C}.
\end{equation}
Using above inequalities, from \eqref{small-onto-ine1} we obtain
$$
|||\Phi(u)|||
\le
\frac12r+C\cdot \frac{1}{(16C)^2}r
\le
r,
\quad
u\in B(r),
$$
which
 shows that $\Phi$ is  from $B(r)$ into $B(r)$.

To show that  $\Phi$ is contraction in $B(r)$,
for $u,v\in B(r)$, 
using  linear estimates 
\eqref{X-forced-linear-nls-est} and \eqref{Y-forced-linear-nls-est}
we get 
\begin{align}
\label{forced-linear-nls-est-small-contr}
|||\Phi(u)-\Phi(v)|||
\le
2c_1\big[
\|
|\tilde u|^2\tilde u-|\tilde v|^2\tilde v
\|_{X^{s,-b}}
+
\|
|\tilde u|^2\tilde u-|\tilde v|^2\tilde v
\|_{Y^{s,-b}}
 \big],
\end{align}
where $\tilde u$, $\tilde v$ extend $u$, $v$ from $\rr^{n-1}\times\rr^+\times(0,T)$ to $\rr^{n+1}$ respectively,
and  satisfy the conditions 
\begin{equation}
\label{uv-ext-cond}
\|\tilde u\|_{X^{s,b}}+
\|\tilde u\|_{Y^{s,b}}
\le
2r
\quad
\text{and}
\quad
\|\tilde v\|_{X^{s,b}}+
\|\tilde v\|_{Y^{s,b}}
\le
6r.
\end{equation}
These extensions are made as follows.
First, we extend  $u$ from $\rr^{n-1}\times\rr^+\times(0,T)$ to $\rr^{n+1}$ such that extension  $\tilde{u}$
satisfies 
$$
%\label{u-extension}
\|\tilde u\|_{X^{s,b}}+
\|\tilde u\|_{Y^{s,b}}
\le
2
\big(
\|u\|_{X^{s,b}_{\rr^{n-1}\times\rr^+\times(0,T)}}
+
\|u\|_{Y^{s,b}_{\rr^{n-1}\times\rr^+\times(0,T)}}
\big)
\le
2r.
$$
Also, we extend $w\doteq u-v$ from 
$\rr^{n-1}\times\rr^+\times(0,T)$ to $\rr^{n+1}$ such that the extension $\tilde{w}$
satisfies 
$$
%\label{u-extension}
\|\tilde w\|_{X^{s,b}}+
\|\tilde w\|_{Y^{s,b}}
=
\|\tilde u-\tilde v\|_{X^{s,b}}+
\|\tilde u-\tilde v\|_{Y^{s,b}}
\le
2
\big(
\|w\|_{X^{s,b}_{\rr^{n-1}\times\rr^+\times(0,T)}}
+
\|w\|_{Y^{s,b}_{\rr^{n-1}\times\rr^+\times(0,T)}}
\big)
\le
4r
$$
Furthermore, we notice that $\tilde v\doteq \tilde u-\tilde w$ extends $v$ from 
$\rr^{n-1}\times\rr^+\times(0,T)$ to $\rr^{n+1}$  and it satisfies
\begin{align*}
\|\tilde v\|_{X^{s,b}}+
\|\tilde v\|_{Y^{s,b}}
\le
(\|\tilde u\|_{X^{s,b}}+
\|\tilde u\|_{Y^{s,b}})
+
(\|\tilde w\|_{X^{s,b}}+
\|\tilde w\|_{Y^{s,b}})
\le
6r.
\end{align*}
Thus, we get the extensions $\tilde{u}$ and $\tilde{v}$, satisfying conditions \eqref{uv-ext-cond}. 
Now,  applying the identity
$
|\tilde u|^2\tilde u-|\tilde v|^2\tilde v
=
\bar{\tilde u}(\tilde u-\tilde v)(\tilde u+\tilde v)
+
\overline{\tilde u-\tilde v}\tilde v^2
$
and using \eqref{uv-ext-cond},
from inequality \eqref{forced-linear-nls-est-small-contr} we obtain
\begin{align}
\label{forced-linear-nls-est-small-contr-no1}
|||\Phi(u)-\Phi(v)|||
\le
104 c_1 r^2
\big(
\|u-v\|_{X^{s,b}_{\rr^{n-1}\times\rr^+\times(0,T)}}+
\|u-v\|_{Y^{s,b}_{\rr^{n-1}\times\rr^+\times(0,T)}}
\big).
\end{align}
Combining this with  inequality \eqref{small-r-cond}, for $u,v\in B(r)  $ we get
\begin{align*}
|||\Phi(u)-\Phi(v)|||
\le
\frac12
\big(
\|u-v\|_{X^{s,b}_{\rr^{n-1}\times\rr^+\times(0,T)}}+
\|u-v\|_{Y^{s,b}_{\rr^{n-1}\times\rr^+\times(0,T)}}
\big)
=
\frac12
|||u-v|||,
\end{align*}
which shows that $\Phi$ is contraction in the ball $B$. 

The proofs for the Lip-continuity of the data-to-solution map
and of the uniqueness of solution are similar to the ones presented in \cite{bop1998} for the well-posedness  of initial value problems.
\,\,
 $\square$

\vskip0.05in
\noindent
{\bf Proof of Theorem \ref{wp-thm}.}
To deal with large data, we modify  iteration map \eqref{iteration-map-T-loc-2} as follows
\begin{align}
\label{iteration-map-T-loc}
u(x, t)
=
\Phi_{T^*}( u)
\doteq
S
\Big[
u_0,g_0;
\mp\psi_{2T^*}
|u|^2u
\Big],
\quad
(x, t)
\in \rr^{n-1} \times \rr^+ \times (0, T),\, \,
T^*\le  T<1/2,
\end{align}
where   $\psi_{T^*}(t)=\psi(t/{T^*})$, which makes $T^*$
the lifespan of the solution to NLS ibvp \eqref{NLS-ibvp}.
Defining the ball in spatial Bourgain spaces
\begin{equation}
\label{r-ball}
B_2(r)
\doteq
\Big\{
u:
u\in X^{s,b}_{\rr^{n-1}\times\rr^+\times(0,T)}
\quad
\text{and}
\quad
\|u\|_{X^{s,b}_{\rr^{n-1}\times\rr^+\times(0,T)}}
\le
r\Big\}
\subseteq X^{s,b}_{\rr^{n-1}\times\rr^+\times(0,T)},
\end{equation} 
and using  linear estimate \eqref{noY-forced-linear-nls-est} and trilinear estimate \eqref{trilinear-est-B}, 
for $u$ in $B_2(r)$, we obtain 
\begin{align}
\label{onto-map-est-no1}
\|\Phi_{T^*}(u)\|_{X^{s,b}_{\rr^{n-1}\times\rr^+\times(0,T)}}
\le
c_1
\Big(
\|
u_0
\|_{H^s(\rr^{n-1}\times\rr^+)}
+
\|g_0\|_{\mathcal{B}_T^s}
+
c_2
\big\|
\psi_{2{T^*}}(\cdot)
\tilde{u}
\big\|_{X^{s,b'}}
\big\|
\tilde{u}
\big\|_{X^{s,b'}}^2
\Big),
\end{align}
where $\tilde{u}$ is the extension of $u$ from $\rr^{n-1}\times\rr^+\times(0,T)$ to $\rr^{n+1}$, 
satisfying
$
%\label{u-extension}
\|\tilde u\|_{X^{s,b}}
\le
2
\|u\|_{X^{s,b}_{\rr^{n-1}\times\rr^+\times(0,T)}}.
$

Now, we use following multipliers estimate,
 whose  proof can be found in  
\cite{tao-book}.
\begin{lemma}
\label{tao-lemma}
Let $\eta(t)$ be a function in the Schwartz space $\mathcal{S}(\rr)$. If 
$
-\frac12
<
b'
\le 
b
<
\frac12,
$
then for any $0<{T^*}\le 1$ we have
\begin{align}
\label{tao-est}
&\|\eta(t/{T^*})u\|_{X^{s,b'}}
\le
c_3(\eta,b',b)
\,\,
{T^*}^{b-b'}\|u\|_{X^{s,b}}.
\end{align}
\end{lemma}
Applying  inequality \eqref{tao-est} with
$
%\label{b-b'-choice-wp-2}
b
=
\frac12-\frac12\beta_{2}
$
and
$
b'
=
\frac12-\beta_{2},
$
where $\beta_{2}=\frac{s}{8}$ with $0<s<\frac12$, is defined in \eqref{beta-choice}, we get 
$
\|
\psi_{2{T^*}}
\tilde{u}
\|_{X^{s,b'}}
\le
c_3
{T^*}^{\beta_{2}/2}
\big\|
\tilde{u}
\big\|_{X^{s,b}},
$
and
combining this   with inequality  \eqref{onto-map-est-no1}
we obtain
\begin{align}
\label{onto-map-est-no2}
\|\Phi_{T^*}(u)\|_{X^{s,b}_{\rr^{n-1}\times\rr^+\times(0,T)}}
\le
c_4
\Big(
\|
u_0
\|_{H^{s}(\rr^{n-1}\times\rr^+)}
+
\|g_0\|_{\mathcal{B}_T^s}
+
{T^*}^{\beta_{2}/2}
\big\|
u
\big\|_{X^{s,b}_{\rr^{n-1}\times\rr^+\times(0,T)}}^3
\Big),
\end{align}
where $c_4\doteq c_1+8c_1c_2c_3$.
Finally, using inequality \eqref{onto-map-est-no2}, we see that
$\Phi_{T^*}$  goes into  $B_2(r)$, if $T^*$ satisfies
\begin{align}
\label{onto-map-condition-nls}
&c_4
\Big(
\|
u_0
\|_{H^{s}(\rr^{n-1}\times\rr^+)}
+
\|g_0\|_{\mathcal{B}_T^s}
\Big)
+
c_4
{T^*}^{\beta_{2}/2}
r^3
\le 
r.
\end{align}

Working similarly, for  $u,v\in B_2(r)$ we obtain
\begin{align}
\label{contraction-map-fin-est-1}
\|\Phi_{T^*}(u)-\Phi_{T^*}(v)\|_{X^{s,b}_{\rr^{n-1}\times\rr^+\times(0,T)}}
\le
108c_4r^2
{T^*}^{
\beta_{2}/2
}
\|
u-v
\|_{X^{s,b}_{\rr^{n-1}\times\rr^+\times(0,T)}},
\end{align}
which makes  $\Phi_{T^*}$ a contraction by choosing
$
r
=
2c_4
\big(
\|
u_0
\|_{H^{s}(\rr^{n-1}\times\rr^+)}
+
\|g_0\|_{\mathcal{B}_T^s}
\big)
$
and 
$T^*$ as in  \eqref{mnls-lifespan}.

The proofs of Lip-continuity and uniqueness is
similar to those for initial value problems \cite{bop1998}.
\,\,
$\square$

%
%
%
%%%%%%%%%%%%%%%%%%%%
%
%
%     UTM
%
%
%%%%%%%%%%%%%%%%%%%% 
%
%
\section{Fokas solution formula for linear Schr\"odinger equation}
\label{sec:Fokas-solution}
\setcounter{equation}{0}
We solve forced linear  ibvp \eqref{lnls-ibvp} on the 
half-space by reducing it  to the half-line case. For this, we take
 the Fourier transform with respect 
to $x' \in \rr^{n-1}$ and letting 
\begin{equation}
\label{Fourier-line}
u_1(\xi',x_n,t)
=
\widehat{u}^{x'}
(\xi',x_n,t)
\doteq
\int_{\rr^{n-1}}
e^{-i\xi'\cdot x'}
u(x',x_n,t)
dx',
\quad
\xi'
\in
\rr^{n-1},
\end{equation}
we reduce the forced linear  ibvp \eqref{lnls-ibvp} 
to the following one
\begin{subequations}
\label{lnls-ibvp-reduce}
\begin{align}
\label{lnls-eqn-reduce}
&i\partial_t u_1
-
\xi'^2 u_1
+
\partial^2_{x_n} u_1
=
\widehat{f}^{x'}(\xi',x_n,t),
\hskip-0.7in
&& 
(\xi', x_n, t)\in \rr^{n-1} \times \rr^+ \times (0, T),
\\
\label{lnls-ic-reduce}
&u_1(\xi',x_n,0) 
= 
\widehat{u}_0^{x'} (\xi',x_n),
\hskip-0.7in
&& (\xi', x_n)\in \rr^{n-1} \times \rr^+,
\\
\label{lnls-bc-reduce}
&u_1(\xi',0, t) 
= 
\widehat{g}_0^{x'}(\xi',t),
\hskip-0.7in
&& (\xi', t)\in \rr^{n-1} \times (0, T)
.
\end{align}
\end{subequations}
%
%
%
%%%%%%%%%%%%%%%%%%
%
%     Fourier transform on half-line  
%
%%%%%%%%%%%%%%%%%%
%
Since the solution of ibvp \eqref{lnls-ibvp} is the inverse Fourier transform of $u_1$ with respect to $\xi'$, in order to solve 
\eqref{lnls-ibvp}, it suffices to solve \eqref{lnls-ibvp-reduce} for any $\xi'\in\rr^{n-1}$. Now, using the integrating factor $e^{i\xi'^2t}$ and letting 
\begin{equation}
\label{u2-def-ls}
u_2(x_n,t)
=
u_2(\xi',x_n,t)
\doteq
e^{i\xi'^2t}\cdot u_1(\xi',x_n,t),
\quad
\text{or}
\quad
i\p_t u_2
=
e^{i\xi'^2t}
\cdot
(i\p_t u_1-\xi'^2 u_1),
\quad
\forall
\,
\xi'
\in
\rr^{n-1},
\end{equation}
and  regarding $\xi'$ as a parameter,
we reduce   ibvp \eqref{lnls-ibvp-reduce} to the following 
one on the half-line $x_n>0$
\begin{subequations}
\label{lnls-half-line}
\begin{align}
\label{lnls-half-line-reduce}
&i\partial_t u_2
+
\partial^2_{x_n} u_2
=
f_2(x_n,t),
\hskip-0.7in
&& 
(x_n, t)\in  \rr^+ \times (0, T),
\\
\label{lnls-half-line-reduce}
&u_2(x_n,0) 
= 
u_{2,0}(x_n),
\hskip-0.7in
&& x_n\in  \rr^+,
\\
\label{lnls-half-line-reduce}
&u_2(0,t)
= 
g_{2,0}(t),
\hskip-0.7in
&& t\in  (0, T),
\end{align}
\end{subequations}
where the forcing, initial data, and boundary data are defined as follows
\begin{align*}
%\label{u2-forcing}
&f_2(x_n,t)
=
f_2(\xi',x_n,t)
\doteq
e^{i\xi'^2t}
\cdot
\widehat{f}^{x'}(\xi',x_n,t)
=
e^{i\xi'^2t}
\cdot
\int_{\rr^{n-1}}
e^{-i\xi'\cdot x'}
f(x',x_n,t)
dx',
\\
%\label{u2-ic}
&u_{2,0}(x_n)
=
u_{2,0}(\xi', x_n)
\doteq
\widehat{u}_0^{x'} (\xi',x_n)
=
\int_{\rr^{n-1}}
e^{-i\xi'\cdot x'}
u_0(x',x_n)
dx',
\\
%\label{u2-bc}
&g_{2,0}(t)
=
g_{2,0}(\xi',t)
\doteq
e^{i\xi'^2t}
\cdot
\widehat{g}_0^{x'}(\xi',t)
=
e^{i\xi'^2t}
\int_{\rr^{n-1}}
e^{-i\xi'\cdot x'}
g_0(x',t)
dx'.
\end{align*}
Now, using Fokas method \cite{f2008}, we
 solve the linear Schr\"odinger  ibvp \eqref{lnls-half-line} and obtain the formula
\begin{align*}
%\label{u2-utm}
u_2(x_n,t)
=&
\frac{1}{2\pi}\int_{\rr}e^{i\xi_n x_n-i\xi_n^{2}t}[\widehat{u}_{2,0}(\xi_n)-iF_2(\xi_n,t)]d\xi_n
-
\frac{1}{2\pi}\int_{\p D^+}e^{i\xi_n x_n-i\xi_n^{2}t}[\widehat{u}_{2,0}(-\xi_n)-iF_2(-\xi_n,t)]d\xi_n
\nonumber
\\
+&
\frac{1}{\pi}\int_{\p D^+}e^{i\xi_n x_n-i\xi_n^{2}t}\xi_n\tilde{g}_{2,0}(-\xi_n^{2},T)d\xi_n,
\end{align*}
where the domain $D^+$ is given by Figure \ref{nls-domain},  $\widehat{u}_{2,0}(-\xi_n)$ is the half-line Fourier transform with respect to $x_n$, $\tilde{g}_{2,0}$ is the time transform of $g_{2,0}$ and $F_2(\xi_n,t)$ is the spatio-temporal transform of $f_2$, i.e.
\begin{align*}
%\label{half-line-FT}
&\widehat {u}_{2,0}(\xi_n)
=
\widehat {u}_{2,0}(\xi',\xi_n)
\doteq
\int_0^\infty
e^{-ix_n\xi_n}
u_{2,0}(x_n)
dx_n,
\\
%\label{time-transform}
&
\tilde g_{2,0}(-\xi_n^2,T)
=
\tilde g_{2,0}(\xi',-\xi_n^2,T)
\doteq
\int_0^Te^{i\xi_n^2 t} g_{2,0}(t)dt,
\\
%\label{space-time-transform}
&F_2(\xi_n,t)
=
F_2(\xi',\xi_n,t)
\doteq
\int_0^t e^{i\xi_n^2\tau}\widehat f_2(\xi_n,\tau)d\tau
=
\int_0^t e^{i\xi_n^2\tau}\int_0^\infty e^{-i\xi_n x_n}f_2(x_n,\tau)dx_nd\tau,
\quad
t\in (0,T).
\end{align*}
Finally, using definition \eqref{Fourier-line} and \eqref{u2-def-ls}, we get  $\widehat{u}^{x'}
(\xi',x_n,t)=u_1(\xi',x_n,t)=e^{-i\xi'^2t} u_2(\xi',x_n,t)$. Then, taking the inverse Fourier transform with respect to $\xi'$, we obtain
the Fokas solution formula \eqref{nd-lnls-utm-sln}.

%%%%%%%%%%%%%%%%%%%%%%
%
%
%		Acknowledgements
%
%
%
%%%%%%%%%%%%%%%%%%%%%%  
%
%
\vskip0.1in
\noindent
{\bf Acknowledgements.} The first author was partially supported by a grant from the 
Simons Foundation (\#524469 to Alex Himonas).  
The authors would like to thank the referee for constructive 
input.

\vskip0.15in
\noindent
{\bf \Large  Declarations}
 
 \vskip0.05in
 \noindent
{\bf \large Ethical Approval. }
Not applicable

 \vskip0.05in
 \noindent
{\bf \large  Competing interests. }
Both authors have no interests of a financial or personal nature.

 \vskip0.05in
 \noindent
{\bf \large  Authors' contributions. }
Both authors wrote and reviewed the manuscript.

 \vskip0.05in
 \noindent
{\bf \large  Funding. }
This work was partially supported by a grant from the 
Simons Foundation (\#524469 to Alex Himonas).
 
 \vskip0.05in
 \noindent
{\bf \large  Availability of data and materials. }
 Data sharing not applicable to this article as no datasets were used
 in this study.

%
%
%%%%%%%%%%%%%%%%%%%%%%  
%
%
%. Bibliography  
%
%
%
%%%%%%%%%%%%%%%%%%%%%  
%
%
%

%

%
%%%%%%%%%%%
%
% Author's address
%
%%%%%%%%%%%
%
\vspace{1mm}
\noindent
A. Alexandrou Himonas  \hfill Fangchi Yan\\
Department of Mathematics  \hfill Department of Mathematics\\
University of Notre Dame  \hfill Virginia Polytechnic Institute and State University\\
Notre Dame, IN 46556  \hfill Blacksburg, VA 24061 \\
E-mail: \textit{himonas.1$@$nd.edu}  \hfill E-mail: \textit{fyan1@alumni.nd.edu}

\end{document}